\documentclass[a4paper,12pt]{report}

\usepackage{stdclsdv}
\usepackage{tocloft}
\usepackage{bm}
\usepackage{amsmath}
\usepackage{amsfonts}
\usepackage{amssymb}
\usepackage{latexsym}
\usepackage{amscd}
\usepackage{ifthen}
\usepackage{graphicx}
\usepackage[mathscr]{eucal}

\usepackage{verbatim}
\usepackage{cases}
\usepackage[tips,graph,matrix]{xy}

\newcommand{\bu}{\bullet}
\newcommand{\mc}[1]{\mathcal{#1}}
\newcommand{\C}{\mathbb{C}}
\newcommand{\Z}{\mathbb{Z}}

\newcommand{\N}{\mathbb{N}}

\renewcommand{\H}{\mathscr{H}}
\newcommand{\G}{\varGamma}
\newcommand{\Gh}{\hat{\G}}
\newcommand{\A}{\mc{A}}
\newcommand{\Cc}{\mc{C}}
\newcommand{\R}{\mathscr{R}}
\renewcommand{\S}{\mc{S}}

\newcommand{\g}{\mathfrak{g}}
\newcommand{\f}{\mathfrak{f}}
\newcommand{\m}{\mathfrak{m}}
\renewcommand{\k}{\Bbbk}

\newcommand{\D}{\mc{D}}
\newcommand{\B}{\mc{B}}

\newcommand{\X}{\mc{X}}

\newcommand{\ot}{\otimes}
\newcommand{\rar}{\longrightarrow}
\renewcommand{\c}{\circ}
\renewcommand{\_}[1]{_{\scsc (#1)}}
\renewcommand{\^}[1]{^{\scsc (#1)}}
\newcommand{\ydcat}[1]{{}_{\scsc #1}^{\scsc #1}\mathscr{YD}}

\newcommand{\scsc}{\scriptscriptstyle}

\renewcommand{\th}{\theta}
\newcommand{\alp}{\alpha}
\newcommand{\si}{\sigma}

\newcommand{\lam}{\lambda_{ij}}
\renewcommand{\a}{a_{ij}}
\newcommand{\uu}{{\bm u}}

\newcommand{\xdig}{\cup} %?or \times ??
\newcommand{\dtab}[7]{$#1$&$#2$&$#3$&$#4$&$#5$&$#6$&$#7$\\}
\newcommand{\dtabb}[6]{$#1$&$#2$&$#3$&$#4$&$#5$&$#6$\rule[-0.8em]{0pt}{0.5em}\\}
\newcommand{\extab}[3]{(4.#1)&\parbox{.3\textwidth}{\centering\includegraphics[scale=1]{#2}}&$#3$\\}

\DeclareMathOperator{\gr}{\bm gr}
\DeclareMathOperator{\del}{\Delta}
\DeclareMathOperator{\ep}{\varepsilon}

\DeclareMathOperator{\id}{id}
\DeclareMathOperator{\an}{\mathscr S}
\DeclareMathOperator{\ad}{ad}
\DeclareMathOperator{\ord}{ord}
\DeclareMathOperator{\co}{co}
\DeclareMathOperator{\son}{\mathfrak{so}(n)}
\DeclareMathOperator{\sln}{\mathfrak{sl}(n)}
\DeclareMathOperator{\spn}{\mathfrak{sp}(n)}

\newcommand{\beer}[3][leer]{\ifthenelse{\equal{#1}{leer}}{\begin{#2}}{\begin{#2}[#1]} \begin{description} \item[{}] #3 \ifthenelse{\equal{#2}{exp}}{\mbox{\hspace{2em}} \hspace{\fill} $\bigtriangleup$}{} \end{description}\end{#2}}
\newcommand{\prf}[2][leer]{\begin{description} \ifthenelse{\equal{#1}{leer}}{\item[Proof:]}{\item[Proof #1:]} #2 \hspace{5em} \hspace*{\fill} \bf{qed.} \end{description} }
\newcommand{\rem}[1]{\begin{description} \item[Remark:] {#1} \end{description}}

\newcommand{\liti}[6]{#1, \emph{#2}, #3 {\bf #4} (#5), #6.}%1Name,2Titel,3Journal,4Number,5Year,6Pages
\newcommand{\litibookbf}[5]{#1, \emph{#2}, #3 {\bf #4}, #5.}
\newcommand{\litibook}[4]{#1, \emph{#2}, #3, #4.}

\newcommand{\lkb}{\cdots}
\renewcommand{\l}{\lambda}
\newcommand{\q}{\bm q}

\newcommand{\qibin}[2]{\binom{#1}{#2}_{\hspace{-1ex}q_i\hspace{0ex}}}
\newcommand{\abild}[1]{\begin{center}\includegraphics{#1}\end{center}}
\newcommand{\bbild}[2][1]{\begin{center}\includegraphics[scale=#1]{#2}\end{center}}

\newcommand{\floatbild}[4][1]{\begin{figure}\centering\includegraphics[scale=#1]{#2}\caption{#3}\label{#4}\end{figure}}
\newcommand{\pth}[1]{#1^{\text{th}}}
\newcommand{\uqg}{U_q(\g)}
\newcommand{\tth}{\tilde{\theta}}
\newcommand{\ga}{\gamma}
\newcommand{\de}{\delta}

\newcommand{\uf}{\mathfrak u}
\newcommand{\Uf}{\mathfrak U}

\newcommand{\U}{\mathcal U}

\newcommand{\K}{\mathcal K}
\newcommand{\Y}{\Upsilon}

\renewcommand{\k}{\Bbbk}

\begin{document}
\sloppy
%\author{Daniel Didt}
%\title{Linkable Dynkin diagrams and Quasi-isomorphisms for finite dimensional pointed Hopf algebras}
%\maketitle
\begin{titlepage}
\phantom{Hallo}\vspace{7ex}
\centering

{\Huge Linkable Dynkin diagrams\\\vspace{1ex}and\\\vspace{1ex} Quasi-isomorphisms\\\vspace{1ex} for finite dimensional pointed\\\vspace{1ex} Hopf algebras\\}
\vspace{19ex}
{\Large Dissertation}\\
\vspace{5ex}
{\large an der Fakult\"at f\"ur Mathematik, Informatik und Statistik der Ludwig-Maximilians-Universit\"at M\"unchen\\
\vspace{8ex}
eingereicht von\\}
\vspace{5ex}
{\LARGE Daniel Didt}\\
\vspace{5ex}
M\"unchen, den 8.~November 2002

\end{titlepage}
%\documentclass[a4paper,12pt]{report}

%\begin{document}
%\sloppy

%\begin{titlepage}
%\setlength{\parindent}{0pt}
\thispagestyle{empty}
\phantom{Hallo}
\vspace{85ex}
%\vfill
\parbox{10cm}{1. Gutachter: Herr Prof. Dr. H.-J.~Schneider\\
2. Gutachter: Herr Priv.-Doz. Dr. P.~Schauenburg}\\

Tag des Rigorosums: 10.~Februar 2003\\

%\end{titlepage}

%\end{document}
%%%%%%%%%%%%%%%%%%%%%%%%%%%%%%%%%%%%%%%%%%%%%%%%%%%%%%%%%%%%%%%%%

%\include{deckblatt}
\setcounter{page}{1}
\tableofcontents
\begingroup
\renewcommand*{\addvspace}[1]{}
\listoffigures
%\listoftables
\endgroup
\chapter{Introduction}

Hopf algebras are named after Heinz Hopf, who introduced such objects in 1941 \cite{Hop} to settle a question in homology theory which was posed to him by Cartan.
The first textbook on Hopf algebras \cite{Swe} came out in 1969 and in spite of many interesting results, the number of people studying this field was small.
 
This changed dramatically with the invention of Quantum Groups in the mid-80s. Suddenly, starting with examples by Drinfeld and Jimbo \cite{Dri,Jim}, there was a vast class of
non-commutative and non-cocommutative Hopf algebras coming from deformations of the enveloping algebras of semisimple Lie algebras \cite{FRT}. The original ideas came from the physical theory of integrable systems and there were lots of attempts to apply the new theory to develop a new quantum or $q$-physics. For this, quantum spaces had to be constructed and symmetries quantized \cite{CSSW1,CSSW2,MajMey}. Differential structures for the new spaces had to be defined \cite[Part IV]{KS} and even experimental evidence was looked for \cite{ACM}. Kreimer found a Hopf algebra that can be used to explain the renormalization process of quantum field theories in mathematical terms \cite{CK}. After a decade of fruitful research, most of the physics community however, started converting to String theory. 

This highly active period brought lots of new notions and constructions, explicit computations and some fundamental structural results to the theory of Hopf algebras, and a major movement to classify finite dimensional Hopf algebras was started.
Nichols and Zoeller proved a freeness theorem \cite{NZ}, Zhu showed that a Hopf algebra of dimension $p$ is necessarily the group algebra of the group with $p$ elements \cite{Zhu}. Recently, Ng \cite{Ng} was able to prove that the only Hopf algebras of dimension $p^2$ are group algebras and the Taft algebras \cite{Taf} introduced in 1971. Many low-dimensional Hopf algebras have been classified. For the case of semisimple Hopf algebras, which contains all group algebras, considerable progress was made by translating lots of proofs and results from group theory into this situation. 
For an overview we suggest \cite{Mon-semi,Som}. Etingof and Gelaki were successful in classifying all finite dimensional triangular Hopf algebras \cite{EG-clas}.

The representation theory of Quantum groups turned out to be closely related to the classical theory of Lie algebras in the cases when the deformation parameter $q$ is not a root of unity. In an attempt to examine the case when $q$ \emph{is} a root of unity, Lusztig found an important class of finite dimensional Hopf algebras \cite{Lus-1,Lus-2}. The representation theory of these new examples is related to that of semisimple groups over a field of positive characteristic and of Kac-Moody algebras. Lusztig was even able to interpret these new Quantum groups as kernels of a ``Quantum Frobenius'' map, which coined the term Frobenius-Lusztig kernels. 

The area this thesis is concerned with is pointed Hopf algebras, which includes all the newly found quantized enveloping algebras of Lie algebras and the finite dimensional Frobenius Lusztig kernels. For pointed Hopf algebras, the so-called coradical which for semisimple Hopf algebras is the whole algebra, is just a group algebra. Substantial results in this case were established with the help of the lifting method of Andruskiewitsch and Schneider \cite{AS-Poin}.

Good introductions to Hopf algebras and related topics can be found in any of the textbooks \cite{CP,Kas,KS,Lus-book,Maj-book,Mon}, the survey article \cite{And} on finite dimensional Hopf algebras, and \cite{MSRI-proc}. The proceedings \cite{Wars-proc} feature a nice series of introductory lectures on various aspects of non-commutative geometry, including new developments on generalisations of the theory of Quantum groups to non-compact groups.\\

In this thesis we want to contribute to some classification results for pointed Hopf algebras with abelian coradical found recently by Andruskiewitsch and Schneider \cite{AS-p3,AS-p17,AS-Poin,AS-char}. Their lifting method produces new classes of Hopf algebras. These algebras are constructed from a linking datum consisting of a group, a Dynkin diagram, some linking parameters and a number of group elements and characters fulfilling certain compatibility conditions. These conditions are rather implicit and hence an explicit description of these Hopf algebras is often not easy. In this work we treat various aspects of such a description in detail.

 One of our main contributions is the clarification of the  concept of linking. Based on the original work \cite{AS-p17}, we first introduce some suitable terminology, Definitions \ref{link.diag}-\ref{link.alg}. Then we give an easily applicable criterion, Theorem \ref{mainthm}, that helps in deciding which linkings can produce finite dimensional Hopf algebras and what possible restrictions have to be imposed on the coradical. This involves simply counting certain objects in graphs and computing the so-called genus from this data. We extend this result to treat affine Dynkin diagrams as well, Theorem \ref{affinethm}. Examples of ``exotic'' linkings are given in Figure \ref{exotics}.  Some exceptional cases that usually have to be excluded from classification results come from setups we call self-linkings. We present the prototypes of Hopf algebras arising from such situations in Section \ref{rank2comps}. The new Hopf algebras derived from the diagram $B_2,$ which we compute using a Computer algebra program, are given in Figure \ref{B2rel}.

Another open question concerns the compatibility of the groups and the Dynkin diagrams in a linking datum. Although a general answer seems out of reach, we are able to contribute an answer for the groups $(\Z/(p))^2$ in Theorem \ref{groupthm}. We prove that apart from a few exceptions, all diagrams with at most four vertices can be used for the construction of finite dimensional pointed Hopf algebras with these groups as the coradical. 

Finally, the last major topic of this thesis is the investigation of the relation between the new Hopf algebras constructed by the lifting method. It turns out that different linking parameters lead to quasi-isomorphic Hopf algebras, Theorem \ref{linkthm}. All Hopf algebras that arise from the lifting method using only Dynkin diagrams of type $A_n$ display the same behaviour, Theorem \ref{mixedthm}. This means that all the finite dimensional pointed Hopf algebras constructed in this way, which only differ in their choice of parameters are 2-cocycle deformations of each other. Our proof should be easily adaptable to the Hopf algebras associated with the other types of finite Dynkin diagrams, once all parameters have been determined for these algebras explicitly. This raises the hope that Masuoka's conjecture in \cite{Mas} can be saved in spite of the counter-example in \cite{EG} by specializing it slightly (page \pageref{conjec}).\\

In Chapter \ref{chap.basics}, to fix notation and present some important terminology, we introduce some basic definitions and results. In Chapter \ref{chap.lift} we give an overview of the lifting method and some important applications. Here we adapt the presentation to suit our needs and cover only the aspects related to this thesis. Chapter \ref{chap.link} then deals extensively with all the aspects concerning linkings. The result for the group realization is presented in Chapter \ref{chap.group} and the last chapter contains all the results connected to quasi-isomorphisms. The programs used for determining the relations of the algebras from Section \ref{rank2comps} are listed in the Appendix together with some documentation.\\

I would like to thank Prof. H.-J.~Schneider for his guidance, the referee of \cite{D-link} for some useful remarks, my family for moral support and the Graduiertenkolleg ``Mathematik im Bereich ihrer Wechselwirkung mit der Physik'' for providing me with the scientific environment and financial support that enabled me to carry out this research.
%%% Local Variables: 
%%% mode: latex
%%% TeX-master: "main"
%%% End: 

\chapter{Basics}\label{chap.basics}

In this chapter we want to present most of the basic ingredients that will be used in the rest of this work. Much of the material here is meant only as a quick reference for displaying our conventions, and helps to present this work in a self contained fashion. There are numerous textbooks and nice expositions treating the various sections in much more depth, giving motivations and historical comments. We will point to some references in  appropriate places and suggest \cite{Mon} for the first two sections.

%%%%%%%%%%%%%%%%%%%%%%%%%%%%%%%%%%%%%%%%%%%%%%%%%%%%%%%%%%%%%%%%%%%%%%%%%%%%%%%%

\section{Coalgebras}

A Hopf algebra, the main object of this work, is first an associative algebra over a base field, that we will denote by $\k$\footnote{We assume, unless stated otherwise, that $\k$ is of characteristic 0 and algebraically closed.}. So it is a $\k$-vector space together with a multiplication and a unit. But at the same time it is a coassociative coalgebra, which is a dualized version of an associative algebra. 
\beer{defi}{A coassociative \emph{coalgebra} is a $\k$-vector space $\Cc$ together with two $\k$-linear maps  $$\del : \Cc\rar \Cc\ot \Cc\qquad\text{and}\qquad\ep : \Cc\rar\k,$$ called the \emph{comultiplication} and the \emph{counit}. The coassociativity constraint requires $$(\del\ot\id)\c\del=(\id\ot\del)\c\del,$$ and the counit has to fulfill $$(\ep\ot\id)\c\del=\id=(\id\ot\ep)\c\del.$$}

To be able to do calculations in coalgebras more easily, a certain convention of notation is now widely used. It is based on an original version by Sweedler and Heyneman and helps to denote the comultiplication. Applying $\del$ to an element $c$ of a coalgebra leads to an element in the tensor product of the coalgebra with itself. This tensor product element is normally a sum of simple tensors. To facilitate notation, one leaves out the summation sign and indicates the tensor components with sub-indices in brackets. This leads to the following notation
$$\del(c)=c\_1\ot c\_2.$$

Coalgebras with the property $\tau\c\del=\del,$ where $\tau$ denotes the flip operator that simply exchanges tensor factors, are called \emph{cocommutative}. Elements $g$ with $\del(g)=g\ot g$ and $\ep(g)=1$ are called \emph{group-like} and elements $x$ with $\del(x)=g\ot x+x\ot h,$ where $g$ and $h$ are group-like, are called \emph{$(g,h)$-primitive}. The set of group-like elements of a coalgebra $\Cc$ is denoted by $G(\Cc)$ and the set of $(g,h)$-primitives by $P_{g,h}(\Cc)$.

Standard examples of coalgebras are group algebras where every element of the group is considered group-like, and the enveloping algebras of Lie algebras. Here the elements of the Lie algebra are $(1,1)$-primitive and the counit evaluates to zero. Also, the dual of every finite dimensional algebra is a coalgebra. The dual of a cocommutative coalgebra is a commutative algebra.

We note that a group-like element spans a one dimensional subcoalgebra and $\k G(\Cc)$ is another example of a subcoalgebra. We call a coalgebra \emph{simple} if it contains no proper subcoalgebras and \emph{cosemisimple} if it is a direct sum of simple subcoalgebras. The sum of all the simple subcoalgebras of a coalgebra $\Cc$ is denoted by $\Cc_0$ and called the \emph{coradical}. Hence $\Cc$ is cosemisimple iff $\Cc=\Cc_0$.

A first fundamental fact about coalgebras is that every simple subcoalgebra is finite dimensional. If, moreover, every simple subcoalgebra is one dimensional, then the coalgebra is called \emph{pointed}. In this case the coradical is necessarily the group coalgebra of the group-like elements.

For every coalgebra $\Cc$ we define inductively for $n\ge 1,$
$$\Cc_n:=\{x\in \Cc : \del(x)\in \Cc_0\ot \Cc + \Cc\ot \Cc_{n-1}\}.$$
According to \cite[Theorem 5.2.2]{Mon} we have for every $n\ge 0,$
\begin{itemize}
\item $\Cc_n\subseteq \Cc_{n+1}$ and $\Cc=\bigcup_{n\ge 0}\Cc_n,$
\item $\del \Cc_n\subseteq\sum^n_{i=0}\Cc_i\ot \Cc_{n-i}.$
\end{itemize}
These properties are exactly the ones defining a coalgebra filtration. So we see that the coradical $\Cc_0$ is the bottom piece of such a filtration and all $\Cc_n$ are subcoalgebras of $\Cc$. We call this filtration the \emph{coradical filtration}. Moreover, it can be proved that the lowest term of \emph{any} coalgebra filtration contains the coradical.

To every filtered coalgebra one can associate a graded coalgebra $\gr\Cc$ by setting $\gr\Cc:=\oplus_{n\ge 0} \Cc(n),$ where
$$  \Cc(n):=\Cc_n/\Cc_{n-1},\quad\text{for }n\ge 1,\qquad\text{and}\quad  \Cc(0):=\Cc_0,$$ and extending the structure maps from $\Cc_0$ in a natural way.
For a coalgebra to be graded we have to have 
\begin{itemize}
\item $\Cc=\oplus_{n\ge 0} \Cc(n)\quad$ and
\item $\del \Cc(n)\subseteq\sum_{i=0}^n \Cc(i)\ot \Cc(n-i),\qquad \ep|_{\Cc(n)}=0$ for $n>0.$
\end{itemize}
Here $\Cc(n)$ is not usually a subcoalgebra.
%%%%%%%%%%%%%%%%%%%%%%%%%%%%%%%%%%%%%%%%%%%%%%%%%%%%%%%%%%%%%%%%%%%%%%%%%%%%%%%%%%%%%%%%%%%%%%%%%%%%

\section{Hopf algebras}
\beer{defi}{A \emph{Hopf algebra} $\H$ is 
\begin{itemize}
\item an associative algebra with unit $1$,
\item a coassociative coalgebra with a comultiplication $\del$ and a counit $\ep$, which are both algebra maps, i.e. satisfying 
\begin{equation}\label{delalg} \del(ab)=\del(a)\del(b),\qquad\ep(ab)=\ep(a)\ep(b),\end{equation}
\item equipped with a linear map $\an$, called the antipode, from $\H$ to itself, fulfilling $$\an(a\_1)a\_2=\ep(a)1=a\_1\an(a\_2).$$
\end{itemize}
}
Here the tensor product $\H\ot\H$ is considered as an algebra with component-wise multiplication. A basic property of the antipode is that $\an$ is an algebra antihomomorphism.

As examples of Hopf algebras we again have group algebras where the antipode on a group element $g$ is defined by $\an(g):=g^{-1},$ and enveloping algebras of Lie algebras where we set $\an(x):=-x$ for the elements of the Lie algebra. A slightly more interesting class are the \emph{Taft algebras} $T(\xi).$ For $\xi\in\k$, a root of 1 of order $N$, we set
$$T(\xi):=\k<g,x : gx=\xi xg, g^N=1, x^N=0>.$$
This is a Hopf algebra where the co-structures are determined by the comultiplication on the generators
$$\del(g):=g\ot g,\qquad\del(x):=g\ot x+x\ot 1.$$

A Hopf algebra is \emph{cosemisimple} or \emph{pointed}, if its underlying coalgebra is so. 

The associated graded coalgebra of a coalgebra filtration of the Hopf algebra $\H,$ is again a Hopf algebra if the filtration is a Hopf filtration. For this we need also that $\H_n\H_m\subseteq \H_{n+m}$ and $\an(A_n)\subseteq A_n,$ for all $n,m\ge 0.$ In the case of the coradical filtration this condition is equivalent to the coradical $\H_0$ being a Hopf subalgebra of $\H$. Therefore, for a pointed Hopf algebra the graded coalgebra associated with the coradical filtration is a Hopf algebra, because the coradical is the group algebra of the group-like elements and this is a Hopf subalgebra.

Pointed Hopf algebras comprise a large class of Hopf algebras. Apart from group algebras and enveloping algebras of Lie algebras, every cocommutative Hopf algebra over an algebraically closed field is pointed.
% This can be seen by considering any simple subcoalgebra $D$. It is finite dimensional so the dual $D^*$ is a commutative and simple algebra over $\k$, and thus isomorphic to $\k$. Hence $D$ is also one dimensional.
In addition, when the base field has characteristic 0, the classic Cartier-Kostant-Milnor-Moore theorem states that any cocommutative Hopf algebra is just the semi-direct product of a group algebra and the enveloping algebra of a Lie algebra.

An interesting aspect of the recent work on pointed Hopf algebras is the somewhat converse statement that large classes of pointed Hopf algebras can be obtained from a group algebra and a \emph{deformed} version of the enveloping algebra of a Lie algebra.

%%%%%%%%%%%%%%%%%%%%%%%%%%%%%%%%%%%%%%%%%%%%%%%%%%%%%%%%%%%%%%%%%%%%%%

\section{Yetter-Drinfeld modules}

The exposition in this section follows closely \cite[Section 2]{AS-Poin}.

A (left) \emph{module} over an algebra $\A$ is a $\k$-vector space $M$ and an action $.: \A\ot M\rar M,$ such that for $a,b\in\A$ and $m\in M$ we have $$(ab).m=a.(b.m)\qquad\text{and}\qquad 1.m=m\:.$$
Analogously, the dual notion is defined. A (left) \emph{comodule} for a coalgebra $\Cc$ is a $\k$-vector space $M$ and a coaction $\rho: M\rar \Cc\ot M,$ such that we have $$(\del\ot\id)\c\rho=(\id\ot\rho)\c\rho\qquad\text{and}\qquad(\ep\ot\id)\c\rho=\id.$$
We extend the Sweedler notation to comodules by writing $$\rho(m)=m\_{-1}\ot m\_0.$$ The negative indices stand for the coalgebra components and the (0) index always denotes the comodule component.

A (left) \emph{Yetter-Drinfeld module} $V$ over a Hopf algebra $\H$ is simultaneously a module and a comodule over $\H$, where the action and coaction fulfill the following compatibility condition:
$$ \rho(h.v)=h\_1v\_{-1}\an(h\_3)\ot h\_2.v\_0,\qquad v\in V, h\in\H.$$

The category $\ydcat{\H}$ of Yetter-Drinfeld modules is a braided monoidal category, i.e. there exists a tensor product operation and a natural isomorphism $c_{M,N}: M\ot N\rar N\ot M$ for all $M,N\in\ydcat{\H},$ called the braiding. It is given by
$$ c_{M,N}(m\ot n):=m\_{-1}.n\ot m\_0,\qquad m\in M, n\in N.$$
The tensor product of two Yetter-Drinfeld modules is just the vector space tensor product with the usual tensor product module and comodule structure. For the compatibility condition we check for $m\in M$ and $n\in N,$
\begin{align*}
\rho(h.(m\ot n))&=\rho(h\_1.m\ot h\_2.n)\\
&=\rho(h\_1.m)\_{-1}\rho(h\_2.n)\_{-1}\ot(\rho(h\_1.m)\_0\ot\rho(h\_2.n)\_0)\\
&=h\_1m\_{-1}\an(h\_3)h\_4n\_{-1}\an(h\_6)\ot(h\_2.m\_0\ot h\_5.n\_0)\\
&=h\_1(m\_{-1}n\_{-1})\an(h\_4)\ot(h\_2.m\_0\ot h\_3.n\_0)\\
&=h\_1(m\_{-1}n\_{-1})\an(h\_3)\ot h\_2.(m\_0\ot n\_0)\\
&=\rho(h\_1(m\ot n)\_{-1})\an(h\_3)\ot h\_2.(m\ot n)\_0.
\end{align*}
The first step is the tensor module formula and the second is the tensor comodule formula. In the third step we used the compatibility condition for $M$ and $N$ separately. The fourth step is the definition of the antipode together with the counit axiom, and the last two steps are again the tensor formulas. In this example we get a glimpse of the usefulness of the Sweedler notation.\\
For further references about braided categories we suggest \cite{JS} and \cite[Chapters XI and XIII]{Kas}.
%As a short remark we want to point out that the category of Yetter-Drinfeld modules is equivalent to the one of bicovariant bimodules or Hopf modules.

The notion of a Hopf algebra still makes sense in a braided category. The tensor product allows one to define algebras and coalgebras in the same way as we have done above. The compatibility condition (\ref{delalg}) of multiplication, here denoted by $\m$, and comultiplication usually involves a flip operation $\tau,$
$$ \del\c\m(a\ot b)=(\m\ot\m)\c(\id\ot\tau\ot\id)\c(\del\ot\del).$$
In braided categories we just have to replace the flip operator $\tau$ by the braiding $c$ and the rest remains as before.

As an example for such a Hopf algebra we give the algebra of coinvariants of a Hopf algebra surjection. Let $\H$ and $\H_0$ be Hopf algebras and \mbox{$\pi: \H\rar\H_0$} and \mbox{$\iota: \H_0\rar\H$} Hopf algebra homomorphisms such that $\pi\iota=\id.$ Then the algebra of coinvariants is defined as
$$\R:=\H^{\co\pi}:=\{h\in\H:(\id\ot\pi)\del(h)=h\ot 1\}.$$
$\R$ is a braided Hopf algebra in $\ydcat{\H_0}$. The module structure is given by $g.r:=\iota(g\_1)r\iota(\an(g\_2)),$ for $g\in\H_0$ and $r\in\R.$ The coaction is $(\pi\ot\id)\del.$ As an algebra, $\R$ is just a subalgebra of $\H.$ The comultiplication is defined by $\del_{\R}(r):=r\_1\iota\pi\an(r\_2)\ot r\_3,$ for all $r\in\R.$

Starting with a braided Hopf algebra $\R$ in $\ydcat{\H_0}$ one can reconstruct an ordinary Hopf algebra $\H$ like the one above. For this we consider the biproduct, or \emph{bosonization} $\R\#\H_0$ \cite{Rad-bos,Maj} of $\R$ and $\H_0.$ This is a Hopf algebra with underlying vector space $\R\ot\H_0,$ with multiplication and comultiplication given by
\begin{align*}
(r\# h)(s\# g)&=r(h\_1.s)\# h\_2g,\qquad r,s\in\R,\; g,h\in\H_0,\\
\del(r\# h)&=(r\^1\# (r\^2)\_{-1}h\_1)\ot((r\^2)\_0\# h\_2).
\end{align*}
Here we used the Sweedler notation with upper indices to indicate the comultiplication in $\R$ in order to distinguish it from the coaction of $\H_0$ on $\R$ which is denoted by lower indices.

In the case where the Hopf algebra $\H_0$ is the group algebra of a finite abelian group $\G,$ the structure theory for Yetter-Drinfeld modules becomes very easy. For a finite dimensional $V\in\ydcat{\H_0}$ there exist a basis  $x_1,\dots,x_{\th}$ of $V,$ elements $g_1,\dots,g_{\th}\in\G$ and characters $\chi_1,\dots,\chi_{\th}\in\Gh,$ such that the action and coaction on $V$ take the form
$$h.x_i=\chi_i(h)x_i\quad\text{and}\quad \rho(x_i)=g_i\ot x_i,\qquad\forall h\in\G, 1\le i\le\th.$$
The braiding can then be expressed by
$$ c(x_i\ot x_j)=b_{ij}x_j\ot x_i,\quad\text{where}\quad \bm{b}=(b_{ij})_{1\le i,j\le\th}:=(\chi_j(g_i))_{1\le i,j\le\th}.$$

For more detailed expositions we refer to \cite{AG}.
%%%%%%%%%%%%%%%%%%%%%%%%%%%%%%%%%%%%%%%%%%%%%%%%%%%%%%%%%%%%%%%%%%%%%%%%%%%%%%%%%%%%%%%
\section{Lie algebras}

Lie algebras are already a classic subject. Having been introduced more than a hundred years ago, they now appear in many branches of mathematics and physics. Their connection with Hopf algebras became really apparent in the mid 1980s with the invention of Quantum groups. In contrast to the usual enveloping algebras of Lie algebras, Quantum groups, being generally non-commutative and non-cocommutative, provided Hopf algebras with a vast class of non-trivial examples.

A very good textbook on Lie algebras is \cite{Hum} and for affine Lie algebras we refer to \cite{Kac}.

\beer{defi}{A \emph{Lie algebra} is a $\k$-vector space $\g$ with a bilinear operation $[\cdot,\cdot]: \g\times\g\rar\g,$ called a \emph{Lie-bracket}, which satisfies $[a,a]=0$ and the \emph{Jacobi identity} $[a,[b,c]]+[b,[c,a]]+[c,[a,b]]=0,$ for all $a,b,c\in\g.$ Note that the Lie bracket is usually non-associative.}

The classification of finite dimensional semisimple Lie algebras by Killing and Cartan is a gem of mathematics. Semisimple Lie algebras are direct sums of simple ones, which in turn are non-abelian algebras that do not have any non-trivial ideals.

Every simple finite dimensional Lie algebra corresponds to one of the diagrams in Figure \ref{fig-lie}. These diagrams are named after Dynkin.
\floatbild[0.6]{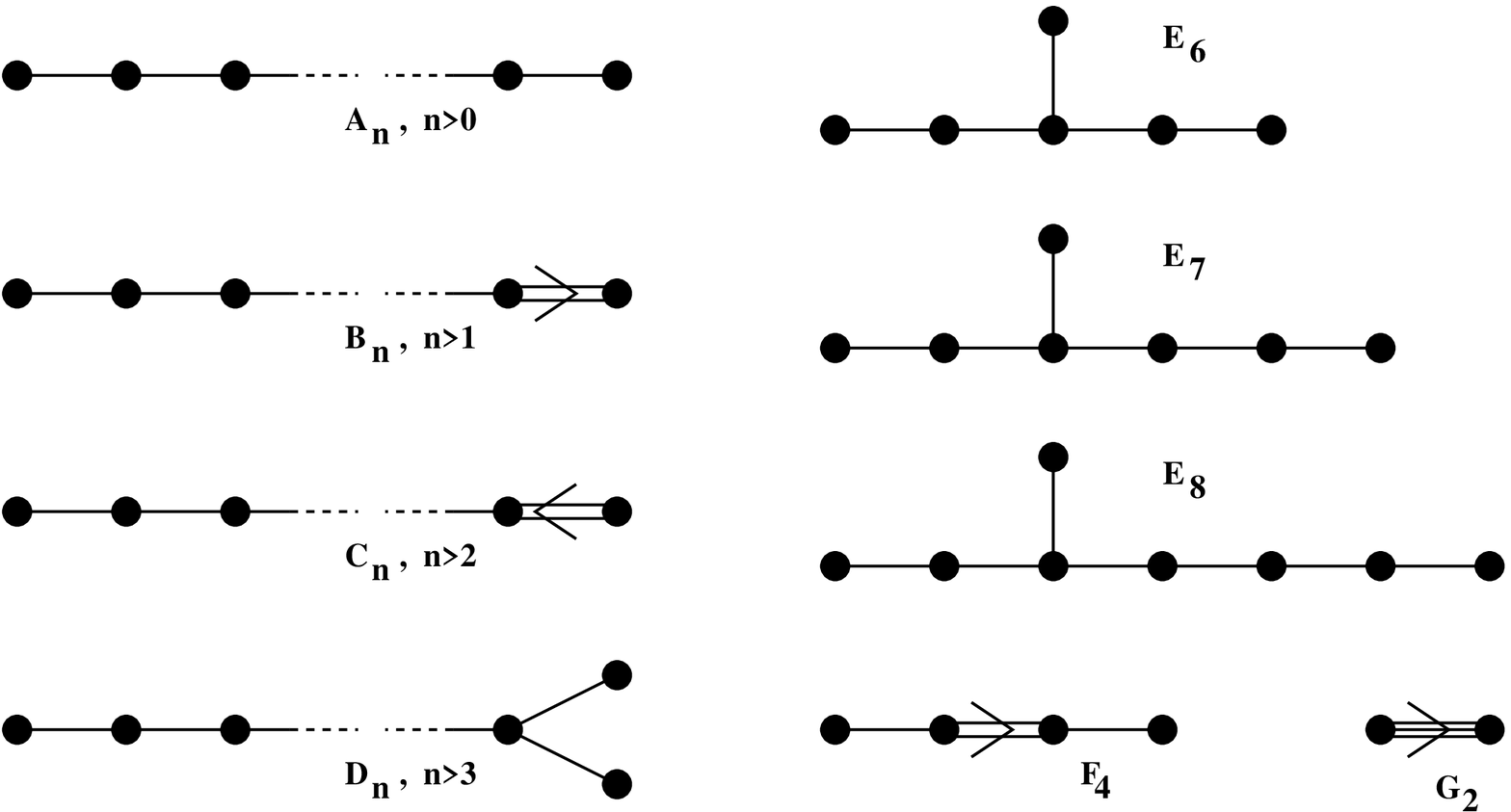}{Dynkin diagrams of finite dimensional simple Lie algebras}{fig-lie}
The four infinite series correspond to the classical examples of $\sln, \son$ and $\spn,$ and then there are 5 exceptional Lie algebras. Instead of this diagrammatic way, there is also the description by Cartan matrices. For a Dynkin diagram $\D$ with $\th$ vertices one takes the $(\th\times\th)$-matrix $\bm{a}=(\a)$ with diagonal entries equal to 2.\\
When two vertices $i$ and $j$ in $\D$ are
\begin{itemize} 
\item not directly connected by a line, we set $\a=a_{ji}=0,$
\item connected by a single line, we set $\a=a_{ji}=-1,$
\item connected by a double line with the arrow pointing at $j$, we set \mbox{$\a=-1,\; a_{ji}=-2,$}
\item connected by a triple line with the arrow pointing at $j$, we set \mbox{$\a=-1,\; a_{ji}=-3.$}
\end{itemize}

There is a so-called root system associated with every semisimple Lie algebra and hence Dynkin diagram. The elements of the root system are called roots and come in two classes: positive and negative roots. For every vertex of the Dynkin diagram we have one simple positive root. Every positive root is a sum of simple positive roots. The number of summands in this presentation is called the \emph{height} of the root. The number of positive and negative roots is the same, and the dimension of the Lie algebra is exactly the number of roots plus the number of vertices of the corresponding Dynkin diagram. We denote the set of positive roots by $\Phi^+$.
  
For affine Lie algebras we have a similar characterization and refer to \cite{Kac} for more details.

%%%%%%%%%%%%%%%%%%%%%%%%%%%%%%%%%%%%%%%%%%%%%%%%%%%%%%%%%%%%%%%%%%%%%%%%%%%%%%%%
\section{Deformation by cocycles}\label{basics.cocyc}
There is a nice deformation operation on Hopf algebras that can provide a Hopf algebra with a new multiplication. For this we need a cocycle.
\beer{defi}{A \emph{2-cocycle $\si$} for the Hopf algebra $\H$ is a linear, convolution-invertible map $\si:\H\ot\H\rar\k$ fulfilling
\begin{align}
\label{cocyc}\si(x\_1,y\_1)\si(x\_2y\_2,z)&=\si(y\_1,z\_1)\si(x,y\_2z\_2)\\
\text{and}\qquad \si(1,1)&=1, \qquad \text{for all }x,y,z\in\H.
\end{align}
Convolution-invertible means that there is another linear map $\si^{-1}:\H\ot\H\rar\k$ such that for all $x,y\in\H,$
\begin{equation}\label{invers}\si\star\si^{-1}(x,y):=\si(x\_1,y\_1)\si^{-1}(x\_2,y\_2)=\ep(x)\ep(y)=%\si^{-1}(x\_1,y\_1)\si(x\_2,y\_2)=:
\si^{-1}\star\si(x,y).
\end{equation}
}

Now, given a 2-cocycle for the Hopf algebra $\H$ we can form a new Hopf algebra $\H_\si$ which, as a coalgebra, is the same as $\H$ but has a new multiplication denoted by $\cdot_\si\:,$
\begin{equation}
\label{newmult}x\cdot_\si y:=\si(x\_1,y\_1)x\_2y\_2\si^{-1}(x\_3,y\_3).
\end{equation}

Given two 2-cocycles $\si,\tau$ for $\H,$ the convolution product $\tau\star\si^{-1}$ is again a 2-cocycle, but for the Hopf algebra $\H_\si.$ The proof of this is a straightforward calculation. By using (\ref{cocyc}) and (\ref{invers}) appropriately we first obtain the cocycle condition for $\si^{-1},$
$$\si^{-1}(x\_1y\_1,z)\si^{-1}(x\_2,y\_2)=\si^{-1}(x,y\_1z\_1)\si^{-1}(y\_2,z\_2).$$
And then we use (\ref{newmult}) to write down (\ref{cocyc}) for $\tau\star\si^{-1}$ more explicitly and get the result by using the cocycle conditions for $\tau$ and $\si^{-1}.$

\section[$q$-Calculus]{$\bm{q}$-Calculus}
We want to collect some basic definitions and results. 
\beer{defi}{For every $q\in\k$ we define for $n,i\in\N$ the 
\begin{itemize}
\item \emph{$q$-numbers} $(n)_q:=\frac{1-q^n}{1-q}=1+q+q^2+\dots +q^{n-1}\:,$
\item \emph{$q$-factorials} $(n)!_q:=(n)_q\cdots (2)_q(1)_q\:,$
\item \emph{$q$-binomial coefficients} $\binom ni_q:=\frac{(n)!_q}{(n-i)!_q(i)!_q}\;.$
\end{itemize}
}
Note that for $q=1$ these are the usual notions.

For $q$-commuting elements $x$ and $y$ in an algebra with $xy=qyx$ we have the quantum binomial formula 
\begin{equation}\label{qbinom}(x+y)^n=\sum_{i=0}^n\binom ni_qy^ix^{n-i}\qquad\text{for all }n\in\N.\end{equation}

%%% Local Variables: 
%%% mode: latex
%%% TeX-master: "main.tex"
%%% End: 

\chapter{The Lifting Method}\label{chap.lift}
We want to give an overview of the so-called lifting method developed by H.-J. Schneider and N. Andruskiewitsch. The method is very general, but most of the results so far concern finite dimensional algebras. For this reason we will limit ourselves mainly to this case. The idea is to break up the classification of finite dimensional pointed Hopf algebras into manageable stages. 

%%%%%%%%%%%%%%%%%%%%%%%%%%%%%%%%%%%%%%%%%%%%%%%%%%%%%%%%%%%%%%%%%%%%%%%%%%%%%%%%
\section{General overview}
We start with a finite dimensional pointed Hopf algebra $\H$ and consider its coradical filtration. The associated graded coalgebra $\gr\H$ is again a pointed Hopf algebra with the \emph{same} coradical, because the coradical is the group algebra of a group $\G$ and hence a Hopf subalgebra.

The algebra $\R$ of coinvariants of the projection $\pi:\gr\H\rar\H_0=\k\G$ is a braided Hopf algebra in the category of Yetter-Drinfeld modules over $\k\G,$ which we will simply denote by $\ydcat{\G}$. By the process of bosonization we can reconstruct $\gr\H$ as $\R\#\H_0.$ We will call $\R$ the \emph{diagram} of $\H$.

The classification of finite dimensional pointed Hopf algebras 
with certain properties
can now be addressed in the following three steps.
\begin{enumerate}
\item Decide what group can be chosen as the coradical, so that it is compatible with the property in question.
\item Find all finite dimensional braided Hopf algebras with the desired property in the Yetter-Drinfeld category of the coradical. %fulfilling the desired property.
\item Find all finite dimensional pointed Hopf algebras whose associated graded version is a bosonization of the ingredients found in the first two steps.
\end{enumerate}
The last step is the actual \emph{lifting} where we have to find ``complicated'' objects over a fairly easy one.

In most applications of this method, the extra property is chosen in such a way as to simplify either step 1 or 2 of the above procedure. For instance, when we fix the coradical of the Hopf algebras in question, step 1 is obvious. If we fix the dimension of the Hopf algebra, then the order of the coradical must divide it and thus, the choice of possible groups is again very limited. Another approach is to fix the diagram $\R$ and hence step 2. Then we need to decide for which groups the diagram $\R$ is actually a Yetter-Drinfeld module. Generally, the lifting method always contains a part that is completely group theoretic.

Step 2 is normally very difficult. The biggest progress was made for the case where the coradical is abelian.
%, because in this case, the structure theory of Yetter-Drinfeld modules is simple.
We will detail this in the next section.

For step 3 it is a priori unclear how to tackle it. But it turns out that in the cases where step 2 can be dealt with satisfactorily, this step becomes manageable too. We present the results for this in Section \ref{lift}.

%%%%%%%%%%%%%%%%%%%%%%%%%%%%%%%%%%%%%%%%%%%%%%%%%%%%%%%%%%%%%%%%%%%%%%%%%%%%%%%%
\section{Nichols algebras}

We introduce a key concept.
\beer{defi}{Let $\S$ be any graded braided Hopf algebra in $\ydcat{\H_0}$ with $\S(0)=\k$ and $\S(1)=P(\S),$ the space of $(1,1)$-primitive elements. The Hopf subalgebra of $\S$ generated as an algebra by $V:=\S(1)$ will be denoted $\B(V)$ and called the \emph{Nichols algebra of $V$}. The dimension of $V$ will be called the \emph{rank} of $\S$.}
An important consequence of the requirement that the generators of the Nichols algebra are \emph{all} the primitive elements is an alternative description of $\B(V)$.
We just take the free Hopf algebra generated by $V$ so that the generators are primitive and divide out all other primitive elements. This allows us to define Nichols algebras for any Yetter-Drinfeld module. 

These algebras appeared first in the work of Nichols \cite{Nic} as the invariant part of his ``bialgebras of part one'', which in turn are the bosonization of a Nichols algebra with the group algebra. Woronowicz used the term ``quantum symmetric algebra'' in \cite{Wor}, and in Lusztig \cite{Lus-book} the algebras $\f$ are examples of Nichols algebras.

The diagram $\R\in\ydcat{\G}$ of a pointed Hopf algebra $\H$ inherits the gradation from $\gr\H,$ where $\R(n)=\gr\H(n)\cap\R.$ Because $\R$ comes from the coradical filtration of $\H$, we can deduce from a Theorem of Taft and Wilson \cite[Theorem 5.4.1]{Mon} that $\R(1)=P(\R)$ and hence we can define the Nichols algebra $\B(V)$ with $V=P(\R).$ We note that $V$ is a Yetter-Drinfeld submodule of $\R.$

A Hopf algebra generated as an algebra by primitive and group-like elements is pointed. This is an easy statement. The converse however, at least for finite dimensional Hopf algebras in characteristic 0, is the main conjecture of Andruskiewitsch and Schneider, cf. \cite[Conjecture 1.4]{AS-cartan}. For the cases where the conjecture is true, we have $\R=\B(V).$

So we see that for step 2 of the lifting method we can limit ourselves to the investigation of the question as to when the Nichols algebra of $V$ is finite dimensional and is compatible with the desired property. For this, $V$ must necessarily be finite dimensional. 

From now on the coradical of $\H$ is the group algebra $\k\G$ with $\G$ finite and abelian. Then we know that $V$ has a basis $x_1,\dots,x_\th$ and there exist $g_1,\dots,g_{\th}\in\G$, $\chi_1,\dots,\chi_{\th}\in\Gh$ such that the action and coaction of $\H_0=\k\G$ take the form
\begin{equation}\label{yddiag}h.x_i=\chi_i(h)x_i\quad\text{and}\quad \rho(x_i)=g_i\ot x_i,\qquad\forall h\in\G, 1\le i\le\th.\end{equation}
On the other hand, given elements $g_i$ and $\chi_i$ as above, we can define a Yetter-Drinfeld module $V\in\ydcat{\G}$ with basis $x_i$ by (\ref{yddiag}) and form the Nichols algebra $\B(V).$

An important role is played by the \emph{braiding matrix}
$$ {\bm b}=(b_{ij})_{1\le i,j\le\th}\qquad\text{with}\quad b_{ij}:=\chi_j(g_i).$$
\beer{defi}{\label{braidmat}A braiding matrix ${\bm b}$ is of 
\begin{itemize}
\item {\bf Cartan type} if
\begin{align}
\label{diagnot1orig}b_{ii}&\neq 1\text{ is a root of unity and}\\
\label{cartcond}b_{ij}b_{ji}&=b_{ii}^{\a}\quad\text{with }\a\in\Z\text{ for all }1\le i,j\le\th.
\end{align}
The integers $\a$ are uniquely determined by requiring $a_{ii}=2$ and $0\ge\a> -\ord b_{ii}$ for $i\neq j.$ Then $(\a)$ is a generalized Cartan matrix, cf. \cite{Kac}.\\
\item {\bf Finite Cartan type} if it is of Cartan type where the Cartan matrix corresponds to a finite dimensional semisimple Lie algebra.
\item {\bf FL-type} if it is of Cartan type with Cartan matrix $(\a)$ and there exist a $q\in\k$ and positive integers $d_1,\dots,d_{\th}$  such that for all $i,j$ 
$$ b_{ij}=q^{d_i\a}\qquad\text{and}\qquad d_i\a=d_ja_{ji}.$$
\item {\bf Local FL-type} if any principal $2\times 2$ submatrix of $(b_{ij})$ is of FL-type.
\end{itemize}
}

The characterization of finite dimensional Nichols algebras over abelian groups is given by the main result of \cite{AS-cartan}.
\beer{thm}{\label{cartthm}\cite[Theorem 1.1.]{AS-cartan} Let $\bm b$ be a braiding of Cartan type and assume that $b_{ii}$ has odd order for all $i$.
\begin{enumerate}
\item If $\bm b$ is of finite Cartan type, then $\B(V)$ is finite dimensional.
\item Assume that $\bm b$ is of local FL-type and that for all $i$, the order of $b_{ii}$ is relatively prime to 3 whenever $\a=-3$ for some $j,$ and is different from 3, 5, 7, 11, 13, 17.\\
If $\B(V)$ is finite dimensional, then $\bm b$ is of finite Cartan type.
\end{enumerate}
}
With this result, the determination of all finite dimensional diagrams $\R\in\ydcat{\G}$ reduces, in many cases, to finding elements $g_i\in\G$ and $\chi_i\in\Gh$ such that the corresponding braiding matrix is of finite Cartan type. This is again partly a group theoretic question. We will deal with a specific problem of this sort in Chapter \ref{chap.group}.

We also have a complete description of $\B(V)$ when the braiding is of finite Cartan type. Let $V$ be a Yetter-Drinfeld module over $\G$ defined by (\ref{yddiag}) and $\bm b,$ the corresponding braiding matrix, is of finite Cartan type. This means that $\bm b$ is associated to a Cartan matrix $(\a)$ and hence to a Dynkin diagram $D$ of a semisimple Lie algebra. We assume that $N_i,$ the order of $b_{ii},$ is odd and not divisible by 3 if $i$ belongs to a connected component of type $G_2$. Here, $i$ is used simultaneously as a vertex in the Dynkin diagram and the corresponding index in the braiding matrix. Let $\X$ be the set of connected components of the Dynkin diagram $D.$ If vertices $i$ and $j$ are in the same component $I\in\X,$ then the orders $N_i$ and $N_j$ of the corresponding braiding matrix entries are equal, due to the symmetry of (\ref{cartcond}). Hence $N_I:=N_i$ is well defined. We define an adjoint action and a braided commutator on the free algebra of $V$ by 
\begin{equation}
  \label{eq:adc}
  (\ad x_i) x_j:=[x_i,x_j]:=x_ix_j-b_{ij}x_jx_i.
\end{equation}
In \cite{Lus-1,Lus-2} Lusztig defined \emph{root vectors}.\label{start.rootvec}\\
For every simple positive root $\alp$ corresponding to the vertex $i$ of the Dynkin diagram, we define the root vector $x_{\alp}:=x_i.$ The root vectors corresponding to all the other positive roots are now defined as iterated braided commutators. The number of commutators equals the height of the root minus one. The entries of the braided commutators are just the $x_i$ corresponding to the simple roots, which form the summands of the positive root. The order of the commutators is the same as in Lusztig's work, where this construction is done for a special braiding. In the second half of the introduction to \cite{Rin} an explicit construction method is given.

As an example we give the root vectors of $G_2$ explicitly. We assume the arrow points at vertex 1, so $a_{12}=-3$ and $a_{21}=-1$. We have the simple roots $\alp_1, \alp_2$ and the positive roots $\alp_1+\alp_2,$ $2\alp_1+\alp_2,$ $3\alp_1+\alp_2$ and $3\alp_1+2\alp_2.$\footnote{There is a slight discrepancy in the notation compared to some literature like \cite{Hum}. For the roots to be the same, we would have to work with the transposed Cartan matrix. This is caused by the Serre relations (\ref{braidserre}), which we want to have in the same form as in the works of Andruskiewitsch and Schneider or \cite{Kac}.} The corresponding root vectors are now%, in agreement with \cite[Section 1]{Rin},
\begin{align}
 x_{\alp_1}&:=x_1, \qquad x_{\alp_2}:=x_2, \\
\label{g2rootz}x_{\alp_1+\alp_2}&:=[x_2,x_1],\\
\label{g2rootu}x_{2\alp_1+\alp_2}&:=[x_{\alp_1+\alp_2},x_1]=[[x_2,x_1],x_1],\\
\label{g2rootv}x_{3\alp_1+\alp_2}&:=[x_{2\alp_1+\alp_2},x_1]=[[[x_2,x_1],x_1],x_1],\\
\label{g2rootw}x_{3\alp_1+2\alp_2}&:=[x_{\alp_1+\alp_2},x_{2\alp_1+\alp_2}]=[[x_2,x_1],[[x_2,x_1],x_1]].
\end{align}
We denote the set of positive roots corresponding to the component $I\in\X$ by $\Phi_I^+.$

\beer{thm}{\label{nichgen}\cite[Theorem 4.5.]{AS-p17} The Nichols algebra $\B(V)$ is presented by generators $x_i, 1\le i\le\th,$ and relations
\begin{align}
\label{braidserre}(ad x_i)^{1-\a}x_j&=0,\qquad i\neq j,\\
x_{\alp}^{N_I}&=0,\qquad \alp\in\Phi_I^+, I\in\X.
\end{align}
The elements $x_{\alp_1}^{n_1}x_{\alp_2}^{n_2}\dots x_{\alp_P}^{n_P}$ with $0\le n_i<N_I,$ if $\alp_i\in\Phi_I^+,$ form a basis of $\B(V).$ Here $P$ is the total number of positive roots and the product involves all root vectors. Hence the dimension of $\B(V)$ is 
$$\dim\B(V)=\prod_{I\in\X}|\Phi_I^+|^{N_I}.$$
}
 
%%%%%%%%%%%%%%%%%%%%%%%%%%%%%%%%%%%%%%%%%%%%%%%%%%%%%%%%%%%%%%%%%%%%%%%%%%%%%%%%
\section{Lifting}\label{lift}

For a more detailed exposition of the material in this section we refer to \cite[Sections 6.2.-6.4.]{AS-Poin}.

We again limit ourselves to pointed Hopf algebras $\H$ with finite abelian coradical $\k\G$. According to the considerations of the last section, the associated graded Hopf algebra $\gr\H$ is, in many cases, just the bosonization of a group algebra and a Nichols algebra of the form given in Theorem \ref{nichgen}. We will give an explicit description.

We fix a presentation $\G=<h_1>\oplus\dots\oplus<h_t>$ and denote by $M_k$ the order of $h_k,$ $1\le k\le t$. Then $\gr\H$ can be presented by generators $y_k,$ $1\le k\le t,$ and $x_i,$ $1\le i\le\th,$ with defining relations
\begin{align}
\label{grprel0}y_k^{M_k}&=1,\quad y_ky_l=y_ly_k,&& 1\le k,l\le t;\\
\label{mixrel0}y_kx_i&=\chi_i(h_k)x_iy_k,&& 1\le k\le t, 1\le i\le\th;\\
\label{serre0}(ad x_i)^{1-\a}(x_j)&=0,&& 1\le i\neq j\le\th;\\
\label{roots0}x_{\alp}^{N_I}&=0,&& \alp\in\Phi_I^+, I\in\X,\\
\intertext{where the Hopf algebra structure is determined by}
\del(y_k)&=y_k\ot y_k,&& 1\le k\le t;\\
\del(x_i)&=x_i\ot 1+g_i\ot x_i,&& 1\le i\le\th.
\end{align}

From the results presented in Section 6 of the survey article \cite{AS-Poin},
% of \cite[Lemma 6.1., Lemma 6.2. and Theorem 6.6.]{AS-Poin} 
we know that apart from a few exceptional cases all finite dimensional pointed Hopf algebras $\H$, with $\gr\H$ as above, can be described in a similar way, i.e. with the same generators and similar relations. The only changes possible are in the ``quantum Serre'' relations (\ref{serre0}) and in the root vector relations (\ref{roots0}). Before we give an explicit formulation, we need some fundamental terminology.

\beer{defi}{\label{link.diag}
Let $(\a)$ be a generalized $(\th\times\th)$-Cartan matrix (cf. \cite{Kac}). The corresponding Dynkin diagram with a number of additional edges, drawn as dotted edges that do not share vertices, will be denoted $D$ and called a \emph{linkable Dynkin diagram}. Two vertices $i$ and $j\neq i$  connected by such dotted edges are called \emph{linkable}. This is written $i\lkb j.$}

\beer{defi}{\label{link.braidmat}
A \emph{linkable braiding matrix of $D$-Cartan type} for a linkable Dynkin diagram $D$ is a $(\th\times\th)$ matrix $(b_{ij})$ with the following properties
\begin{align}
b_{ii}&\neq 1,\\
\label{cartan} b_{ij}b_{ji}&=b_{ii}^{\a},\\
\label{link}  b_{ki}^{1-\a}b_{kj}&=1, \quad k=1,\ldots,\th,\,\text{ if }i \text{ is linkable to } j .
\end{align}
}

\beer{defi}{\label{link.realize}
A linkable braiding matrix $(b_{ij})$ is called \emph{realizable over the abelian group $\G$} if there are elements $g_1,\ldots,g_\th\in\G$ and characters $\chi_1,\ldots\chi_\th\in\Gh$ such that 
\begin{align} 
b_{ij}=\chi_j(g_i),\qquad\text{ for all }i,j\quad\text{ and}\\
\chi_i^{1-\a}\chi_j=1,\qquad\text{ whenever }i\lkb j.
\end{align}
}

\beer{defi}{\label{link.datum}
A \emph{linking datum of Cartan type} is a collection of the following ingredients:
\begin{itemize}
\item an abelian group $\G$,
\item elements $h_1,\ldots h_t\in\G,$ such that $\G=<h_1>\oplus<h_2>\oplus\ldots\oplus<h_t>$, $M_k:=\ord(h_k)$, $k=1,\dots,t$,
\item a linkable Dynkin diagram $D$ with $\th$ vertices and Cartan matrix $(\a)$,
\item a linkable braiding matrix of $D$-Cartan type, which is realizable over $\G$ with elements and characters $g_i\in\G,$ $\chi_i\in\Gh,$ $i=1,\dots,\th,$
\item parameters $\lam\in\k,$ $1\le i\neq j\le\th,$ such that $\lam=0$ if $i$ is not linkable to $j$ and $\lam=-\chi_j(g_i)\l_{ji}$ when $\a=0.$
\end{itemize}
Vertices $i$ and $j$ with $\l_{ij}\neq 0$ are called \emph{linked}. If $i$ and $j$ are linked and lie in the same connected component of the Dynkin diagram, we talk of a \emph{self-linking}.
A \emph{linking datum of finite Cartan type} is a linking datum where the diagonal elements of the braiding matrix have finite order and the Cartan matrix is of finite type, i.e. it corresponds to a finite dimensional semisimple Lie algebra.
}

\rem{The notion of \emph{linkable vertices} is simpler and hence more general than the one given in \cite[Definition 5.1.]{AS-p17}. However, the notions \emph{do} mainly coincide when we require the existence of a linking datum and demand that linkable vertices lie in different connected components of the Dynkin diagram. We note that a linkable braiding matrix of a linking datum of finite Cartan type is a braiding matrix of Cartan type as defined in Definition \ref{braidmat}, and the two notions coincide when there are no linkable vertices.
}

\beer{defi}{\label{link.alg}
The algebra $\Uf(\D)$ for a linking datum $\D$ of Cartan type is given by generators $y_k,$ $1\le k\le t,$ $x_i,$ $1\le i\le\th,$ and relations
\begin{align}
\label{grprel}y_k^{M_k}&=1,\quad y_ky_l=y_ly_k,&& 1\le k,l\le t;\\
\label{mixrel}y_kx_i&=\chi_i(h_k)x_iy_k,&& 1\le k\le t, 1\le i\le\th;\\
\label{serre}(ad x_i)^{1-\a}(x_j)&=\lam(1-g_i^{1-\a}g_j),&& 1\le i\neq j\le\th.
\end{align}
So far, $\ad$ is only a symbol and has the following explicit form:
\begin{equation}
(\ad x_i)^{1-\a}(x_j):=\sum_{k=0}^{1-\a}(-1)^k\qibin{1-\a}kq_i^{\binom k2}b_{ij}^kx_i^{1-\a-k}x_jx_i^k,
\end{equation}
where $q_i:=b_{ii}=\chi_i(g_i).$ The group elements $g_i$ are interpreted as words in the generators $y_k.$

Now let $\D$ be a linking datum of finite Cartan type. We define the root vectors $x_\alp$ for all positive roots $\alp\in\Phi^+$ in the same way as sketched on page \pageref{start.rootvec}. Given a family $\uu=(u_\alp),$  where $u_\alp$ is an expression in the generators $x_i, y_k$ for every positive root $\alp,$ we define the algebra $\uf(\D,\uu)$ in the same way as $\Uf(\D)$, but with the extra relations
\begin{equation} \label{roots} x_{\alp}^{N_I}=u_\alp,\qquad \alp\in\Phi_I^+, I\in\X. \end{equation}
$\X$ is again the set of connected components of the Dynkin diagram, and $N_I$ denotes the common order of those diagonal elements of the braiding matrix that correspond to the component $I$. We call $\uu$ \emph{root vector parameters}.
}

\rem{When $\a=0$, (\ref{serre}) simplifies to
\begin{equation}\label{links} x_ix_j-\chi_j(g_i)x_jx_i=\lam(1-g_ig_j).\end{equation} When there are no self-linkings, then \emph{all} relations (\ref{serre}) with non-vanishing right hand side are of the form (\ref{links}).
}

\beer{prop}{There exists a unique Hopf algebra structure on $\Uf(\D)$  determined by 
\begin{align}
\del(y_k)&=y_k\ot y_k,&& 1\le k\le t;\\
\del(x_i)&=g_i\ot x_i+x_i\ot 1,&& 1\le i\le\th.
\end{align}
If vertices are only linked when they lie in different connected components, then $\uf(\D,\bm 0)$ is also a Hopf algebra with the same comultiplication as above.
}
For $\Uf(\D)$ the proof that the relations define a Hopf ideal is mostly an exercise. Only the ``quantum Serre'' relations (\ref{serre}) need special attention. One has to show that both sides of (\ref{serre}) are $(g_i^{1-\a}g_j,1)$-primitive. For the left hand side one can use, for instance, \cite[Lemma A.1.]{AS-cartan}.\\ The statement about $\uf(\D,\bm 0)$ when all the $u_\alp$ are zero is exactly \cite[Theorem 5.17.]{AS-p17}.

We want to show how these new Hopf algebras are connected with the usual quantized Kac-Moody Hopf algebras $\uqg.$ \\
We start with the direct sum of two copies of the given symmetrizable Cartan matrix. In the associated Dynkin diagram we connect corresponding vertices by dotted lines. This is our linkable Dynkin diagram. We number the vertices of one copy of the original diagram from 1 to $N$ and the remaining ones from $N+1$ to $2N$ in the same order. The group $\G$ is simply $\Z^N.$ For the $g_i,\; 1\leq i\leq N,$ we take the canonical basis of $\G$, set $g_{N+i}:=g_i$ and define characters $\chi_j(g_i):=q^{d_i\a}, \chi_{N+i}:=\chi_i^{-1},$ where $d_i\a=d_ja_{ji}.$ As a linkable braiding matrix of the given Cartan type we can now take $b_{ij}=\chi_j(g_i).$ Then $(b_{ij})$ is even of FL-type. We set $\l_{i(N+i)}:=1,\; 1\leq i\leq N,$ and all other $\l_{ij}:=0$ when $i<j$, so there is no self-linking. The Hopf algebra $\Uf(\D)$ obtained from this linking datum $\D$ is the quantized Kac-Moody algebra. To see this, one sets $K_i:=g_i, K^{-1}_i:=g_i^{-1},E_i:=x_i,F_i:=(q^{-d_i}-q^{d_i})^{-1}x_{N+i}g_i^{-1},\,1\leq i\leq N.$

A linking datum $\D$ where all the $\lam$ are zero is denoted by $\D_0.$ So we see that the Hopf algebras $\gr\H$ given at the beginning of this section by (\ref{grprel0})-(\ref{roots0}) are simply of the form $\uf(\D_0,\bm 0)$. The lifting method seems to indicate now that apart from a few exceptions, all pointed finite dimensional Hopf algebras $\H$ with $\gr\H\simeq\uf(\D_0,\bm 0)$ are of the form $\uf(\D,\uu)$. However, the complete list of possibilities for the root vector parameters $u_\alp$ has been found only in a few cases.

%%%%%%%%%%%%%%%%%%%%%%%%%%%%%%%%%%%%%%%%%%%%%%%%%%%%%%%%%%%%%%%%%%%%%%%%%%%%%%%%
\section{Examples}

We want to present some examples of the successful application of the lifting method.
%%%%%%%%%%%%%%%%%%%%%%%%%%%%%%%%%%%%%%%%%%%%%%%%%%%%%%%%%%%%%%%%%%%%%%%%%%%%%%%%
\subsection[Classification of pointed Hopf algebras of dimension $p^3$]{Classification of pointed Hopf algebras of dimension $\bm{p^3}$}
In \cite{AS-p3} the authors classified all pointed non-cosemisimple Hopf algebras $\H$ of dimension $p^3,$ $p$ an odd prime, with the help of their lifting method. This was done independently in \cite{CD} and \cite{SvO}. According to the Nichols-Zoeller theorem \cite[Theorem 3.1.5]{Mon}, the dimension of the coradical, being a Hopf subalgebra, has to divide $p^3$. Hence for the algebra to be non-cosemisimple, the coradical must have order $p$ or $p^2$. So the diagram $\R$ must have dimension $p^2$ or $p$. The authors proved that in these cases $\R$ is a Nichols algebra of a Yetter-Drinfeld module $V$, and the braiding matrix is of finite Cartan type with Dynkin diagram $A_1$ or $A_1\xdig A_1.$ They determined all possible linking and root vector parameters and gave a complete list of such Hopf algebras of dimension $p^3$. As a bonus they considered the coradical $\Z/(p^2)$ and a two dimensional module $V$ with braiding of type $A_1\xdig A_1$. The lifting of the corresponding graded Hopf algebra produces an infinite family of non-isomorphic pointed Hopf algebras of dimension $p^4$. This was one of the first counterexamples to a conjecture of Kaplansky.

%%%%%%%%%%%%%%%%%%%%%%%%%%%%%%%%%%%%%%%%%%%%%%%%%%%%%%%%%%%%%%%%%%%%%%%%%%%%%%%%
\subsection[Classification of pointed Hopf algebras of dimension $p^n$]{Classification of pointed Hopf algebras of dimension $\bm{p^n}$}
Pointed Hopf algebras of dimension $p$ or $p^2$ are just group algebras or Taft algebras. The case $n=3$ was explained in the previous subsection. For $n=4$ in \cite{AS-A2} and $n=5$ in \cite{Gn-p5}, similar strategies were used to obtain the classification. Again, the coradical can have only special orders. The possible Dynkin diagrams appearing in these cases are $A_2,$ $B_2$ or copies of $A_1.$ One could go on like that for $n>5,$ but the explicit list of the algebras would soon become unmanageable. We refer to Subsection \ref{subp17} for an important class in such a classification.

%%%%%%%%%%%%%%%%%%%%%%%%%%%%%%%%%%%%%%%%%%%%%%%%%%%%%%%%%%%%%%%%%%%%%%%%%%%%%%%%
\subsection[Lifting of Nichols algebras of type $A_n$ and $B_2$]{Lifting of Nichols algebras of type $\bm{A_n}$ and $\bm{B_2}$}
Here the strategy is to start with a Yetter-Drinfeld module $V\in\ydcat{\G}$ whose braiding is of finite Cartan type with Dynkin diagram $A_2$ (see \cite{AS-A2}) or $B_2$ (in \cite{BDR-B2}) or $A_n$ \cite[Section 7]{AS-Poin}. Without specifying the group $\G$, all possible liftings in such a situation are then determined. These are the few cases where the generalized root vector relations (\ref{roots}) are known explicitly. The question of which groups actually admit such Yetter-Drinfeld modules and in how many ways has still to be addressed.

Because all the diagrams considered here have only one connected component, the lifted Hopf algebras have no linking parameters. However, there are a few exceptional cases where the lifting method is not as straightforward as described in the general picture above. In \cite[Section 3]{AS-A2} for instance, the authors could not deal with a case called $p=3$ for the diagram $A_2$. This was then done in \cite{BDR-B2}, but at the same time the authors were not able to treat $p=5$ for $B_2$. We will give an answer to this in Section \ref{rank2comps}. Anticipating further developments we will also provide a partial answer for the exceptional case $p=7$ of the diagram $G_2.$ Here $p$ denotes the order of the diagonal elements of the braiding matrix. 

%%%%%%%%%%%%%%%%%%%%%%%%%%%%%%%%%%%%%%%%%%%%%%%%%%%%%%%%%%%%%%%%%%%%%%%%%%%%%%%%
\subsection[Classification of pointed Hopf algebras with coradical $(\Z/(p))^s$]{Classification of pointed Hopf algebras with coradical $\bm{(\Z/(p))^s}$}\label{subp17}
In \cite{AS-p17} the authors are able to give a complete classification of \emph{all} (and this time there really are \emph{no} exceptions) pointed finite dimensional Hopf algebras whose coradical consists of an arbitrary number of copies of the group with $p$ elements, where $p$ is a prime bigger than 17.

By having the group consisting of cyclic groups of prime order, the root vector parameters can only be zero. And $p>17$ ensures that the exceptional cases for the lifting procedure and the ones mentioned in Theorem \ref{cartthm} do not interfere.

\beer{thm}{\label{p17thm}\cite[Theorem 1.1.]{AS-p17}\\
(a). Let $p>17$ be a prime and $\H$ a pointed finite dimensional Hopf algebra such that $G(\H)\simeq\G:=(\Z/(p))^s$. Then there exists
%\begin{itemize}
%\item a finite Cartan matrix $(\a)\in\k^{\th\times\th}$;
%\item elements $g_1,\dots,g_\th\in\G, \chi_1,\dots,\chi_\th\in\Gh$ such that
%\begin{align}
%\label{diagnot1}\chi_i(g_i)&\neq 1,&& 1\le i\le\th,\\
%\label{cartan}\chi_i(g_j)\chi_j(g_i)&=\chi_i(g_i)^{\a},&& 1\le i,j\le\th;
%\end{align}
%\item and a collection $(\lam)_{1\le i<j\le\th, i\nsim j}$\footnote{The symbol $i\nsim j$ means that vertices $i$ and $j$ are not in the same connected component of the Dynkin diagram.} of elements in $\k$ where $\lam=0$ unless $\chi_i\chi_j=1$ and $g_ig_j\neq 1;$
%\end{itemize}
a linking datum $\D$ of finite Cartan type with group $\G$ and no self-linkings
such that $\H\simeq\uf(\D,\bm 0)$.
% can be presented as an algebra by generators $x_1,\dots x_\th, y_1,\dots,y_s$ and relations
%\begin{align}
%\label{grouprel}y_k^p&=1,\quad y_ky_l=y_ly_k,&& 1\le k,l\le s,\\
%\label{mixrel}y_kx_i&=\chi_i(y_k)x_iy_k,&& 1\le k\le s, 1\le j\le\th,\\
%\label{serrein}(ad x_i)^{1-\a}(x_j)&=0,&& 1\le i\neq j\le\th,\, i\sim j,\\
%\label{serreout}x_ix_j-\chi_j(g_i)x_jx_i&=\lam(1-g_ig_j),&& 1\le i<j\le\th,\, i\nsim j,\\
%\label{rootrel}x_{\alp}^p&=0,&& \alp\in\Phi^+;
%\end{align}
%and where the Hopf algebra structure is determined by
%\begin{align}
%\label{costruct}\del y_k&=y_k\ot y_k,&& 1\le k\le s,\\
%\qquad \del x_i&=x_i\ot 1+g_i\ot x_i,&& 1\le i\le\th.
%\end{align}

(b). Conversely, given a linking datum $\D$ of finite Cartan type with group $\G,$
% let $(\a)\in\Z^{\th\times\th}$ be a finite Cartan matrix, $g_1,\dots,g_{\th}\in\G,$ $\chi_1,\dots,\chi_{\th}\in\Gh$ such that (\ref{diagnot1}), (\ref{cartan}) hold and $\lam$ as above. Assume that $p>3$ if the Cartan matrix has a connected component of type $G_2$. Then
the algebra $\H:=\uf(\D,\bm 0)$ 
%presented by generators $x_1,\dots,x_{\th},$ $y_1,\dots,y_s$ and relations (\ref{grouprel})-(\ref{rootrel}) has a unique Hopf algebra structure determined by (\ref{costruct}). $\H$ 
is pointed, $G(\H)\simeq\G$ and $\dim\H=p^{s+|\Phi^+|}.$
}

Although this result provides a good answer to the classification problem, there are still a few difficulties when we want to know all Hopf algebras of this kind explicitly. This is the starting point of this Ph.D. thesis.

One aspect needing clarification is the linking parameters. Having fixed the Dynkin diagram, the group and character
elements, what possible $\l$ can appear? It is not at all obvious which vertices can be linked. ``Exotic'' linkings like \cite[Example 5.13.]{AS-p17}, where 4 copies of $A_3$ are linked into a circle are possible. The general picture was presented in \cite{D-link} and the next chapter is devoted to this problem. 

Given a fixed prime $p$ and an $s,$ what Dynkin diagrams are realizable? In other words, for which diagrams can we find group elements and characters such that (\ref{cartan}) can be fulfilled? This question has been addressed so far only for $s=1$ in \cite{AS-cartan}.
% There the authors found as well some general results in the form of bounds. 
We will present the answers for $s=2$ in Chapter \ref{chap.group}.

In the last chapter we will be concerned with analyzing how different all these new Hopf algebras actually are.

%%% Local Variables: 
%%% mode: latex
%%% TeX-master: "main.tex"
%%% End: 

\chapter{The structure of linkable Dynkin diagrams}\label{chap.link}

In this chapter we want to address the problem of determining all possible linkings. 
We will be mainly concerned with a detailed investigation of when a linkable braiding matrix of a given Cartan type does exist. This will lead to a characterization of the corresponding linkable Dynkin diagrams. We show how these ideas are related to the usual quantized enveloping algebras and to the finite dimensional Hopf algebras constructed in \cite{AS-p17}, which are themselves variations of the finite dimensional Quantum groups called Frobenius-Lusztig kernels \cite{Lus-book}.

To get a nice result, we slightly specialize some of  our earlier definitions. We will discuss generalisations later.
\begin{quote}From now on all linkable Dynkin diagrams are assumed to be link-connected, i.e. when viewed as a graph they are connected. Furthermore, we will restrict our considerations to diagrams where two vertices are linkable only if they lie in different connected components of the original diagram, i.e. there are no self-linkings. Finally, all diagonal elements of the braiding matrix have finite order and the base field $\k$ is required to contain a $\pth p$ root of unity for a prime $p>3.$
\end{quote}
For two vertices $i,j$ of the Dynkin diagram with $a_{ij}\neq0,$ the symmetry of (\ref{cartan}) implies \begin{equation}\label{normdiag}b_{ii}^{a_{ij}}=b_{jj}^{a_{ji}}.\end{equation}For $i\lkb j$ we have $a_{ij}=0$, as we required the vertices to lie in different connection components. Using (\ref{link}) and (\ref{cartan}) alternately, we arrive at \begin{equation}\label{diaglink}b_{ii}=b_{ij}^{-1}=b_{ji}=b_{jj}^{-1}.\end{equation}

%%%%%%%%%%%%%%%%%%%%%%%%%%%%%%%%%%%%%%%%%%%%%%%%%%%%%%%%%%%%%%%%%%%%%%%%%%%%%%%%

\section{The finite case}
First we will only consider Dynkin diagrams of finite type, i.e. the corresponding Lie algebras are finite dimensional.
In order to get interesting applications in regard of \cite{AS-p17} we further require that a linkable braiding matrix has the following property:\\

\refstepcounter{equation}
\hfill\begin{parbox}{.8\textwidth}{\emph{The order of the diagonal elements $b_{ii}$ is greater than 2 and not divisible by 3 if the linkable Dynkin diagram contains a component of type $G_2$.}}
\end{parbox}\hfill  (\theequation)\\ \label{order}

The first properties are presented in a lemma, which is essentially Lemma 5.6. in \cite{AS-p17}. However, we formulate it on the level of the braiding matrix.

\beer{lemma}{\label{aslemma}We are given a linkable Dynkin diagram $D$ and a corresponding linkable braiding matrix $\bm b.$ Suppose that the vertices $i$ and $j$ are linkable to $k$ and $l$, respectively. Then $a_{ij}=a_{kl}.$
}
\prf{
If $a_{il}\neq 0$ or $a_{jk}\neq{0}$ then we immediately get $a_{ij}=a_{kl}=0$, because linkable vertices must lie in different connected components of $D.$ So we now take $a_{il}=a_{jk}=0.$
Without loss of generality we assume $a_{ij}\le a_{kl}.$
Using (\ref{cartan}) and (\ref{link}) alternately, we get
$$b_{ii}^{a_{ij}}=b_{ij}b_{ji}=b_{il}^{-1}b_{jk}^{-1}=b_{li}b_{kj}=b_{lk}^{-1}b_{kl}^{-1}=b_{kk}^{-a_{kl}}=b_{ii}^{a_{kl}}.$$
In the last step we used (\ref{diaglink}). Hence $a_{ij}=a_{kl}$ modulo the order of $b_{ii}.$  As $b_{ii}\neq\pm 1,$ we either get $a_{ij}=a_{kl}$ or that the order of $b_{ii}$ is 3 and $a_{ij}=-3, a_{kl}=0.$ But in the last case $i$ and $j$ form a $G_2$ component.
So $b_{ii}=3$ is a contradiction to the assumption on the order of the diagonal elements.
}

Before we can state our result on the structure of linkable Dynkin diagrams that admit a corresponding braiding matrix with the above properties, we have to introduce some terminology.
\beer{defi}{
For every cycle\footnote{A cycle is a closed, non self-intersecting path in the diagram.} $c$ in $D$ we choose an orientation and denote by the \emph{weight} $w_c$ the absolute value of the difference of the numbers of double edges in that cycle with the arrow pointing with the orientation and against it. The \emph{length} $l_c$ of the cycle is defined to be the number of dotted edges in that cycle.\\
The \emph{genus} $g_c$ of the cycle is now defined by the following formula:
\begin{equation}
g_c:=2^{w_c}-(-1)^{l_c}.
\end{equation}
}
In preparation for some technicalities in the second part of the proof of our result we also need the following concept.
\beer{defi}{
For two vertices $i$ and $j$ of $D$ we define for every directed path $P$ from $j$ to $i$ a number $h^i_j(P)\geq 0$, called the \emph{height of $i$ over $j$ along $P$}, by the following algorithm.

First we set $h=0$. Then we follow the path $P$  starting at $j$. At every vertex we get to, we 
\begin{quote}\begin{center}{\footnotesize increase the value of $h$ by 1\\ decrease it by 1\\ or leave it unchanged,}\end{center}\end{quote}
 depending on if the edge we just passed was a double edge pointing  
\begin{quote}\begin{center}{\footnotesize with the orientation of $P$\\ against it\\ or was not a double edge.}\end{center}\end{quote} 
The only exception is that the value of $h$ is not decreased when it is 0. $h^i_j(P)$ is then set to be the value of $h$ after we followed through the whole path $P$ arriving at $i.$

For a cycle $c$ we define the \emph{natural orientation} to be the one where the number of double edges in $c$ pointing with this orientation is not bigger than the number of double edges pointing against it\footnote{If the weight $w_c=0$ then the natural orientation is ambiguous. In that case we choose one of the possible two orientations. This will not lead to any problems.}.

For every vertex $i$ of $c$ we define the \emph{absolute height} $h_i(c)\geq 0$ to be the height of $i$ over itself along $c$ following its natural orientation. A vertex of absolute height $0$ in a cycle of genus $g_c>0$ is called a \emph{Level 0 vertex.}
}
This seems to be the right point to illustrate all the notions in an example. We consider the Dynkin diagram in Figure \ref{dynkin-example}, where the vertices are supposed to be linkable in the indicated way.
\floatbild{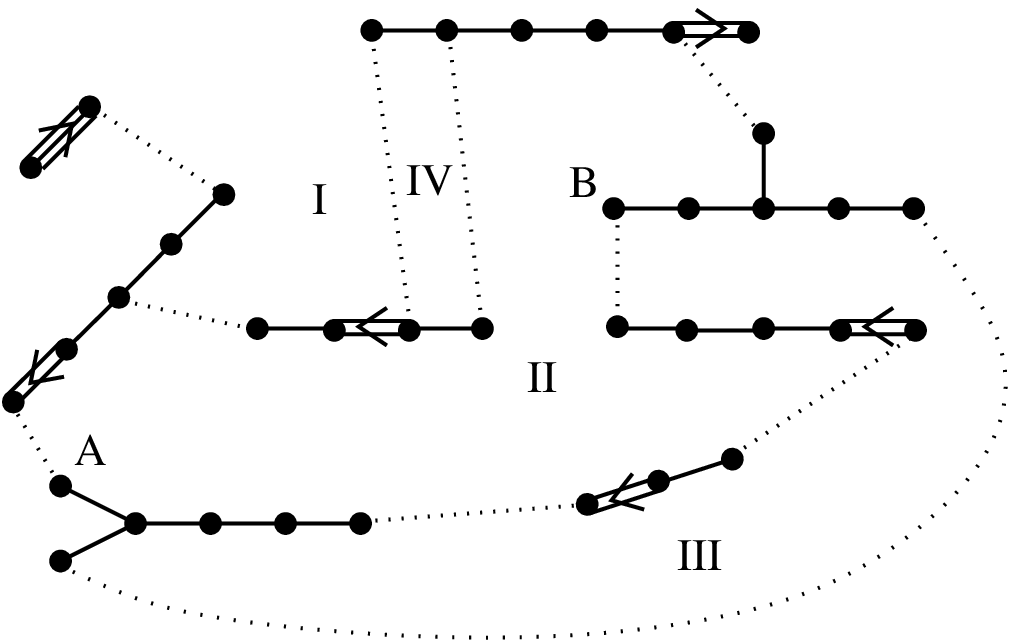}{Example of a linkable Dynkin diagram}{dynkin-example}
For the four cycles denoted by I, II, III and IV (where I is the outside cycle) the values for $w$, $l$ and $g$ are given in this table:
%Table \ref{dynkin-table}.
%\begin{table}
\begin{center}
\begin{tabular}[b]{c | c c c}
&$w_c$&$l_c$&$g_c$\\
\hline
I&2&5&5\\
II&2&7&5\\
III&0&4&0\\
IV&0&2&0\\
\end{tabular}.\end{center}
%\caption{Weights, lengths and genera for cycles of Figure \ref{dynkin-example}}\label{dynkin-table}\end{table}
The natural orientation of cycles I and II is clockwise, whereas the natural orientation in cycles III and IV is ambiguous. The vertex indicated by the letter ``A'' is a vertex of absolute height 1 in cycle II, but a Level 0 vertex for cycle I. And Vertex ``B'' is a Level 0 vertex for cycle II but a vertex of absolute height 1 for cycle III, independent of the natural orientation chosen for that cycle.

We are now able to come to our first main result.
\beer{thm}{
\label{mainthm}
We are given a link-connected linkable Dynkin diagram $D$ and explicitly exclude the case $G_2\xdig G_2$. It will be treated later.

A linkable braiding matrix of $D$-Cartan type exists, iff
\begin{enumerate}
\item\label{cd1} In components of type $G_2$ not both vertices are linkable to other vertices.
\item\label{cd2} $D$ does not contain any induced subgraphs\footnote{An induced subgraph consists of a subset of the original vertices and all the corresponding edges.} of the form:
\abild{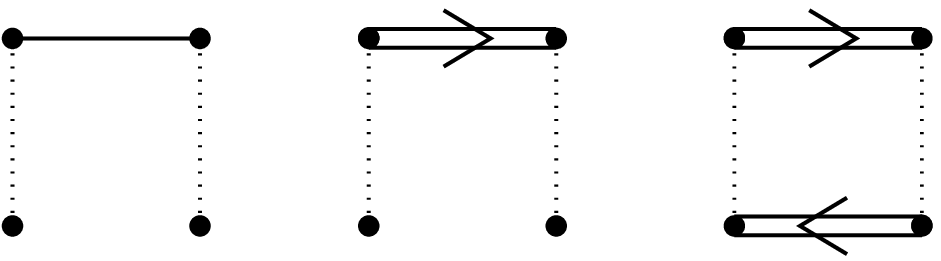}
\item\label{cd3} One of the following conditions is true: 
\begin{itemize}
\item[    ] $D$ contains no cycles or the genera of all cycles are zero.
\item[or  ] $D$ does not contain a component of type $G_2$ and there is a common divisor $d>2$ of all cycle genera and the field $\k$ contains a primitive $\pth d$ root of unity.
\item[or  ] $D$ does contain a component of type $G_2$ and there is a common divisor $d>2$ of all cycle genera, $d$ is not divisible by 3 and the field $\k$ contains a primitive $\pth d$ root of unity.
\end{itemize}
\end{enumerate}
}
\prf{
We first prove the ``if'' part, i.e. we assume conditions \ref{cd1}-\ref{cd3}.

We will construct the braiding matrix explicitly and show that it  fulfills the required identities.

The main observation is that once an element of the diagonal has been chosen, the other diagonal elements are determined (up to possible signs) by (\ref{normdiag}) and (\ref{diaglink}).

We take $d>2$ as given by condition \ref{cd3}. In the case that $D$ contains no cycles or the genera of all cycles are zero we set $d>2$ to be a prime, such that $\k$ contains a primitive $\pth d$ root of unity. This is possible by the general assumption on the field $\k$.  We note in particular that $d$ is always odd.

Now we choose a vertex $i$ 
and set $b_{ii}:=q$, where $q$ is a primitive $\pth d$ root of unity. As the Dynkin diagram $D$ is link-connected, we can choose for every vertex $j\neq i$ a path\footnote{Again we demand that a path does not include a vertex more than once.} $P_{ij}$ connecting $i$ and $j$, which we denote by the sequence of its vertices $(i=p_0,p_1,\dots,p_t=j).$ For every such path $P_{ij}$ we now define the $b_{p_kp_k},\; k=1,\dots,t,$ recursively: 
\begin{equation}\label{recurs}b_{p_{k+1}p_{k+1}}=\begin{cases} b_{p_kp_k}^{-1},&\text{if }a_{p_kp_{k+1}}=0,\\
                                   b_{p_kp_k}^{a_{p_kp_{k+1}}/a_{p_{k+1}p_k}},&\text{otherwise.}\end{cases}\end{equation}
When taking the square root we choose it to be again a $\pth d$ root of unity. As $d$ is odd, this determines exactly one of the two possible roots. When taking third roots we choose them to be the unique $\pth d$ root of unity as well. This is always possible, as we only have to take third roots when we are dealing with vertices of components of type $G_2.$ Then however, $d$ is not divisible by 3.

We now show that this process is well defined.

The only time a problem could arise, would be when a vertex can be reached from the starting vertex by two different paths, i.e. when $D$ contains a cycle. Suppose we are given two different paths $(i=p_0,p_1,\dots,p_t=j)$ and $(i=q_0,q_1,\dots,q_u=j).$ Let $n$ be the smallest integer with $p_{n+1}\neq q_{n+1}$ and $m_1>n$ the smallest integer, such that there is a $m_2>n$ with $p_{m_1}=q_{m_2}$. Then $c=(p_n,p_{n+1},\dots,p_{m_1}=q_{m_2},q_{m_2-1},\dots,q_n=p_n)$ is a cycle. It is now sufficient to show that the recursive procedure (\ref{recurs}) for the paths $P_1:=P_{p_np_{m_1}}$ and $P_2:=P_{q_nq_{m_2}}$ leads to the same value 
\begin{equation}\label{equal}b_{p_{m_1}p_{m_1}}=b_{q_{m_2}q_{m_2}}.\end{equation}
As triple edges are not part of cycles, we easily obtain a closed formula for the desired values:
\begin{equation}b_{p_{m_1}p_{m_1}}=Q^{(-1)^{l_1}2^{w_1}}\qquad \text{and}\qquad b_{q_{m_2}q_{m_2}}=Q^{(-1)^{l_2}2^{w_2}}.\end{equation}
Here $Q:=b_{p_np_n}$, $l_i\ge 0$ denotes the number of dotted edges in the path $P_i$ and $w_i\in\Z$ is the difference of the numbers of double edges in $P_i$ that have the arrow pointing against the path's orientation and with it. Without loss of generality we assume $w_1\ge w_2$ and have $w_c=w_1-w_2$ and $l_{c}=l_1+l_2$.

As $Q$ is a $\pth d$ root of unity and $d$ divides all cycle genera we get $Q^{g_c}=1$ or $$Q^{2^{(w_1-w_2)}}=Q^{(-1)^{(l_1+l_2)}}.$$
Taking both sides to the power of $(-1)^{l_1}2^{w_2},$ we arrive at (\ref{equal}). Here we would like to remind the reader that all values are $\pth d$ roots of unity and hence there is no ambiguity regarding signs.

The so specified diagonal entries of the braiding matrix fulfill the requirements for their orders because $d$ is odd and not divisible by 3 when there are components of type $G_2$ in $D.$\\

We now give the remaining entries for the linkable braiding matrix, i.e. we specify $b_{ij}$ for $i\neq j.$ For this we divide the set $\{(i,j):i\neq j\}$ of pairs of vertices into 4 classes:
\begin{description}
\item[None of the two vertices
is linkable to any
 other one.]
We set $$b_{ji}:=z,\qquad b_{ij}:=b_{ii}^{a_{ij}}z^{-1}.$$
\item[The two vertices
are linkable to each other.]
We set %\begin{align*}
$$b_{ij}:=b_{ii}^{-1},\qquad b_{ji}:=b_{jj}^{-1}.$$
%\end{align*}   
\item[Only one of the two vertices
is linkable to another vertex.] We assume $i$ is linkable to $k$. We set
\begin{align*} b_{ji}&:=z,  &\qquad b_{ij}&:=b_{ii}^{a_{ij}}z^{-1},\\
                  b_{jk}&:=z^{-1}, &\qquad b_{kj}&:=b_{kk}^{a_{kj}}z. \end{align*}
\item[Both vertices
are linkable to other vertices.]
We assume $i$ is linkable to $k$ and $j$ is linkable to $l.$ For $i$ and $k$ to be linkable we can not have $a_{ij}\neq 0$ and $a_{jk}\neq 0.$ So after a possible renaming of the indices $i$ and $k,$ we can assume that $a_{jk}=0.$ By the same reasoning we take $a_{il}=0.$ Now we set
\begin{align*} b_{ji}&:=b_{kj}:=z, &\qquad b_{ij}&:=b_{li}:=b_{ii}^{a_{ij}}z^{-1},\\
                  b_{jk}&:=b_{kl}:=z^{-1}, &\qquad b_{il}&:=b_{lk}:=b_{ii}^{-a_{ij}}z.\end{align*}
\end{description}
In all the cases, $z\neq 0$ can be chosen freely from the field $\k$ and can be different for every class and pair of vertices.

We would like to point out that all pairs of indices fall into one of those classes and that there are no overlapping cases, i.e. each off-diagonal element is only set in one of these.

In this way we have explicitly constructed the matrix $\bm b=(b_{ij})$. We are left to show that (\ref{cartan}) and (\ref{link}) are fulfilled. For the diagonal entries this has already been done. For the entries being set in the first three classes it is immediately clear from the definition.

In the last class only the relation $b_{kl}b_{lk}=b_{kk}^{a_{kl}}$ must still be checked.

We note that neither $i$ and $j$ nor $k$ and $l$ can form a component of type $G_2$, as this would contradict part \ref{cd1} of the assumption. From the construction we get $b_{kl}b_{lk}=b_{ii}^{-a_{ij}}.$ As vertex $i$ is linkable to $k$ we know $b_{kk}=b_{ii}^{-1}.$ We will show that $a_{kl}=a_{ij}.$
 
If $a_{ij}=0$ then we immediately get $a_{kl}=0$, because $a_{kl}a_{lk}=1$ or 2 is not permitted by part \ref{cd2} of the assumption. Analogously we get the result if we assume $a_{kl}=0.$

The case where the 4 indices form a sub-diagram of the kind
\abild{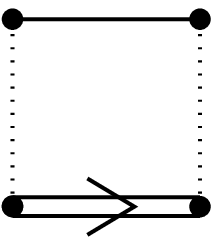}
is excluded, as the cycle genus for this diagram is $1=2^1-(-1)^2.$ So the only other possible diagrams these four vertices can form are
\abild{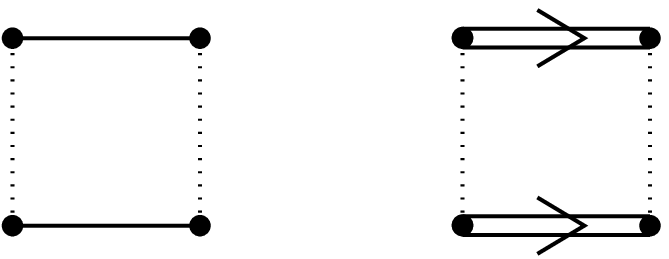}
which all have $a_{ij}=a_{kl}.$ We would like to note that the cycle genus for the last diagram is $0=2^0-(-1)^2.$

This concludes the ``if'' part of the proof.\\

Before we come to the ``only if'' part we prove a lemma to enable us to deal with some arising technicalities.
\beer{lemma}{
\label{tech}
We are given a linkable Dynkin diagram $D$.
\begin{enumerate}
\item In every cycle $c$ of $D$ with $g_c>0$ there exists a Level 0 vertex.
\item \label{simplecyc}Given a linkable braiding matrix $\bm b$ of $D$-Cartan type we have for every Level 0 vertex $i$ of $c$: \quad $b_{ii}^{g_c}=1.$
\item \label{allcyc}Let $G$ be the greatest common divisor of all cycle genera. If there are no cycles with cycle genus 1 or 2 we have for every Level 0 vertex $i$:\quad $b_{ii}^G=1.$
\end{enumerate}
}
\prf[of Lemma]{
\begin{enumerate}
\item Pick a vertex $i$ in $c$
and calculate $h_i(c)$. If $h_i(c)>0$ then take the vertex $j$, where (in the recursive definition) $h$ was 0 for the last time.

Then $j$ is a Level 0 vertex, because:
\begin{itemize}
\item In the recursive calculation of $h_j(c)$ the value of $h$ is positive at least until we pass vertex $i$. (Choice of $j$ and $i$ is not a Level 0 vertex.)
\item Assume $h$ stays positive until it reaches $j$ again, i.e. $j$ is not Level 0. Then however, the number of double edges along the cycle $c$ pointing with the natural orientation is greater than the number of double edges pointing against it, which is a contradiction to the definition of natural orientation.
\item So there is a vertex $k$ between $i$ and $j$ where $h$ becomes 0. But this means that the value of $h$ at $k$ in the calculation of $h_i(c)$ was 0 as well. (This value can not be bigger, having started out smaller at vertex $i.$)
\item Now, the calculation of $h$ from $k$ until $j$ is the same as for $h_i(c).$ So $h=0$ when it reaches $j.$ 
\end{itemize}

\item Set $Q:=b_{ii}.$ Following the cycle in its natural orientation starting at $i$ and using (\ref{normdiag}) and (\ref{diaglink}) we arrive at
$$Q=Q^{(-1)^{l_c}2^{w_c}}.$$
That there are no extra signs from possible square roots in the above formula is ensured by the assumption that $i$ is of Level 0. Raising both sides to the $\pth{{(-1)^{l_c}}}$ power and dividing by the new left hand side we get $1=Q^{g_c}.$
\item If there are no cycles $c$ with $g_c>0$ there is nothing to show. When there is only one such cycle $c$ we have $G=g_c$ and part \ref{simplecyc} of this Lemma establishes the claim.

We now take two cycles $c_1$ and $c_2$ with $g_{c_1}\ge g_{c_2}>2$ and set $g:=\gcd(g_{c_1},g_{c_2}).$ 
We pick in $c_1$ a vertex $i$ of Level 0 and a vertex $j$ of Level 0 in $c_2.$ Following a path from $i$ to $j$ and applying (\ref{normdiag}) and (\ref{diaglink}) appropriately we get $b_{jj}^{2^{w_1}}=b_{ii}^{(-1)^l2^{w_2}}$ for some values $w_1,w_2$ and $l.$ $b_{ii}$ is a $\pth{g_{c_1}}$ root of unity (see previous part of this Lemma) and so $b_{jj}$ is also a $\pth{g_{c_1}}$ root of unity. However, using the previous part again, $b_{jj}$ must be a $\pth{g_{c_2}}$ root of unity. As the cycle genera are not divisible by 2 we conclude that $b_{jj}$ is a $\pth g$ root of unity.

Repeating the argument for all the other cycles with cycle genus bigger than 2, we conclude that for every Level 0 vertex $i,$ the corresponding diagonal entry $b_{ii}$ of the braiding matrix is a $\pth G$ root of unity. 
\end{enumerate}
}
Now we finish the proof of the theorem. We are given a linkable braiding matrix $\mathbb b$ and a corresponding linkable Dynkin diagram $D$ and set $G$ to be the greatest common divisor of all cycle genera. $G:=0$ if all cycle genera are 0.

Assume now:
\begin{itemize}
\item there is a $G_2$ component with both vertices $i$ and $j$ linkable to other vertices $k$ and $l$, respectively.

From Lemma \ref{aslemma} we immediately get $a_{kl}=a_{ij}$ and $a_{lk}=a_{ji}.$ So $k$ and $l$ form another $G_2$ component. Thus, the given diagram must be $G_2\xdig G_2,$ which we do not want to consider here.
\item one induced subgraph is of the kind as in condition \ref{cd2}.

This contradicts Lemma \ref{aslemma}.
\item there is a cycle $c$ with cycle genus $g_c=1$ or 2.

Using Lemma \ref{tech} \ref{simplecyc} we get that a diagonal entry of the braiding matrix must be 1 or $-1$.
\item $G=1$ and there is no cycle with genus 1 or 2.

Lemma \ref{tech} \ref{allcyc} then shows that $b_{ii}=1$ for all Level 0 vertices $i$, a contradiction to our assumption (\ref{order}) about the order.
\item $G>2$ and the base field $\k$ does not contain a primitive $\pth d$ root of unity for any $d>2$ dividing $G$.

Lemma \ref{tech} \ref{allcyc} shows that the field $\k$ must contain a $\pth G$ root of unity. This means there must be a $d$ dividing $G$, such that $\k$ contains a primitive $\pth d$ root of unity. By the order assumption (\ref{order}) 
there is even a $d$ bigger than 2.
\item $G>2$, there is a component of type $G_2$ and the only $d>2$ that divide $G$, such that $\k$ contains a primitive $\pth d$ root of unity, are divisible by 3.

There is a cycle $c$ with $g_c>2.$ According to Lemma \ref{tech} there is a Level 0 vertex $i$ in $c$ with $b_{ii}$ a primitive $\pth d$ root of unity, where $d>2$ and $d$ divides $G$. From the assumption we get that $d$ must be divisible by 3, so the order of $b_{ii}$ is divisible by 3. This contradicts (\ref{order}).
 \end{itemize}
}

This result allows us to decide quickly which vertices in Dynkin diagrams can be linked to produce new Hopf algebras. This also gives a much more explicit description of the algebras presented in Section \ref{subp17}. For every Dynkin diagram we can choose several pairs of vertices which we want to link. We then apply Theorem \ref{mainthm} to determine if such a setup is at all possible. If it is possible, then we set $s$ to the number of vertices in the diagram and choose a prime $p$ which divides all cycle genera. Then the braiding matrix is realizable over the group $(\Z/(p))^s.$ It might be possible to choose a smaller $s,$ but this kind of question is closely related to the general problem of which diagrams can feature at all for a given group. We will discuss this in Chapter \ref{chap.group} for the special case $s=2.$

%%%%%%%%%%%%%%%%%%%%%%%%%%%%%%%%%%%%%%%%%%%%%%%%%%%%%%%%%%%%%%%%%%%%%%%%%%%%%%%%

\section{The affine case}\label{affine-linkable}
Now we turn to the affine case, i.e. we consider Dynkin diagrams that are unions of diagrams of finite and affine type. For notations regarding affine algebras we refer to \cite[Chapter 4]{Kac}. In order to get a similar result to the previous one we have to consider an even more specialized style of braiding matrix.\\

\refstepcounter{equation}
\hfill\begin{parbox}{.8\textwidth}{\emph{We require that the order of all diagonal entries is the same and equal to a prime bigger than 3. This kind of braiding matrix we will call \emph{homogeneous}.}}
\end{parbox}\hfill  (\theequation)\\ \label{afforder}

We have the analogue of Lemma \ref{aslemma} for this situation.
\beer{lemma}{
\label{asaff}We are given a linkable Dynkin diagram $D$ and a corresponding homogeneous linkable braiding matrix $\bm b.$ Suppose that the vertices $i$ and $j$ are linkable to $k$ and $l$, respectively. Then $a_{ij}=a_{kl}.$
}

The proof is the same as before and the last conclusion is straightforward, as $b_{ii}$ has to have at least order 5, according to the assumptions.

As before, we define the notions of weight and natural orientation for every cycle $c$. This time however, we do this as well for triple edges. The former notions will now be denoted by \emph{natural 2-orientation} and $w^2_c$ and the corresponding ones for the triple edges by \emph{natural 3-orientation} and $w^3_c.$ The length of the cycle is exactly as before.

If the two natural orientations coincide we define the genus of this cycle as $$g_c:=3^{w^3_c}2^{w^2_c}-(-1)^{l_c}.$$
In the other case we take $$g_c:=|3^{w^3_c}-2^{w^2_c}(-1)^{l_c}|.$$
Now the theorem can be formulated in the same spirit.
\beer{thm}{\label{affinethm}
We are given a link-connected linkable affine Dynkin diagram $D$ and explicitly exclude the cases $A_1^{(1)}\xdig A_1^{(1)}$ and $A_2^{(2)}\xdig A_2^{(2)}$. These will be treated later.\\
A homogeneous linkable braiding matrix of $D$-Cartan type exists, iff
\begin{enumerate}
\item\label{dc1} In components of type $A_1^{(1)}$ and $A_2^{(2)}$ not both vertices are linkable to other vertices.
\item\label{dc2} $D$ does not contain any induced subgraphs of the form:
\bbild[.8]{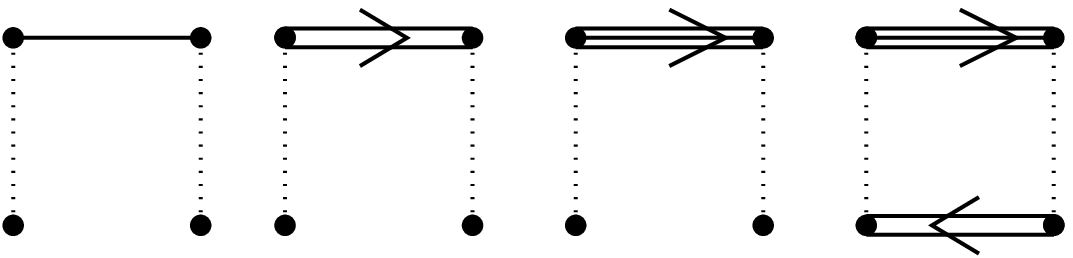}
\item\label{dc3} There is a prime $p>3$ that does divide all cycle genera and the field $\k$ contains a $\pth p$ root of unity.
\end{enumerate}
}
\prf[(Sketch)]{
%(sketch)
We start again with the ``if'' part. Suppose all conditions are fulfilled.

We pick any vertex $i$ and set $b_{ii}$ to be a $\pth p$ root of unity, where $p$ is given by condition \ref{dc3}. As before, we can now set all other diagonal elements of the braiding matrix inductively. Dotted, single, double and triple edges are treated as in the finite case, the edge of type $A_1^{(1)}$ is treated as a single edge, and when we pass quadruple edges we take the $\pth 4$ root or power according to if we go along \emph{with} the arrow or in the opposite direction. This time we only require that all diagonal values are $\pth p$ roots of unity and as $p>3$ is a prime this determines the values uniquely. The independence from the paths chosen in this procedure is ensured again by the condition on the genera, and can be shown by elaborating on the technique used in the finite case.

The off-diagonal elements are set completely in the same way as in the finite case, so we will not repeat the arguments here.\\

The ``only if'' part again follows the same strategy as before. We assume to be given a linkable Dynkin diagram $D$ and a homogeneous linkable braiding matrix of $D$-Cartan type.
Suppose condition \ref{dc1} is not fulfilled. Then Lemma \ref{asaff} immediately gives that the diagram is $A_1^{(1)}\xdig A_1^{(1)}$ or $A_2^{(2)}\xdig A_2^{(2)}.$

The negation of condition \ref{dc2} contradicts Lemma \ref{asaff}.

Suppose now that the biggest prime $p$ dividing all cycle genera is smaller than 5, or that the greatest common divisor $G$ of all cycle genera is 1. 
Here the analogue of Lemma \ref{tech} \ref{allcyc} gives the contradiction. One does not need to use the concept of Level 0 vertices however, as all diagonal entries have to have the same order. Still, the result that every $b_{ii}$ is a $\pth G$ root of unity can be deduced by the same reasoning.
}

%%%%%%%%%%%%%%%%%%%%%%%%%%%%%%%%%%%%%%%%%%%%%%%%%%%%%%%%%%%%%%%%%%%%%%%%%%%%%%%%

\section{The excluded cases}
Suppose the given diagram is one of the above excluded ones, i.e. $G_2\xdig G_2,$ $A_1^{(1)}\xdig A_1^{(1)}$ or $A_2^{(2)}\xdig A_2^{(2)}.$
%When the linkable Dynkin diagram has only one pair of linkable vertices and hence no cycles, the
%When the linkable Dynkin diagram has a cycle, a linkable braiding matrix does exist, according to Lemma \ref{asaff}, only if corresponding vertices are required to be linkable. So
We label the vertices at which the arrows in each copy point as vertices 1 and 3 and the others as 2 and 4, so that vertices 1 and 2 form one of the copies of $G_2,$ $A_1^{(1)}$ or $A_2^{(2)}$ and vertices 3 and 4 form the other.

When the linkable Dynkin diagram has a cycle, a linkable braiding matrix exists, according to Lemma \ref{asaff}, only if vertices 1 and 3 as well as 2 and 4 are linkable. In this case the following linkable braiding matrices exist
%. Then vertices 1 and 3 as well as 2 and 4 are assumed linkable and vertices 1 and 2 form one of the copies of $G_2,$ $A_1^{(1)}$ or $A_2^{(2)}$ and vertices 3 and 4 form the other. Now  the following linkable braiding matrices do exist
%If the linkable diagram is not excluded by Lemma \ref{asaff} then the following linkable braiding matrices do exist for any of the possible linkings
\begin{equation*} \begin{pmatrix} q&z&q&z^{-1}q^{-m}\\z^{-1}q^{-m}&q^n&z&q^n\\q^{-1}&z^{-1}&q^{-1}&zq^m\\zq^m&q^{-n}&z^{-1}&q^{-n}\end{pmatrix}. \end{equation*}
Here $n=m=3$; $n=1,m=2$ or $n=m=4$ for the first, second and last diagram respectively, and $z\neq 0.$

When there is no cycle, the situation is simpler, and braiding matrices similar to the one above can be written down.

%%%%%%%%%%%%%%%%%%%%%%%%%%%%%%%%%%%%%%%%%%%%%%%%%%%%%%%%%%%%%%%%%%%%%%%%%%%%%%%%

\section{Examples}
If we start with two copies of a Dynkin diagram and link corresponding vertices, all cycles have genus 0, so there are no obstructions. This shows once more that braiding matrices necessary for the construction of the quantum universal enveloping algebras $\uqg$ and quantized Kac-Moody algebras exist.

If we take $n$ copies of $A_3$ and link them into a circle \cite[Example 5.13.]{AS-p17}, then we have one cycle of weight 0 and length $n$. So the genus is 0 for $n$ even, and there are no restrictions. For $n$ odd, however, the genus is 2, so under the conditions we imposed there are no corresponding linkable braiding matrices.

For the case of $n$ copies of $B_3$ linked into a circle the genus is $2^n-(-1)^n.$

In Figure \ref{dynkin-example} which we used to demonstrate the various notions, we found that all cycle genera were divisible by 5. So we know that all linkable braiding matrices of that type have diagonal elements whose order is divisible by 5. This gives a limitation on the groups over which these matrices can be realized.

\begin{figure}
\centering
%\vspace{-5cm}
\includegraphics[scale=.7]{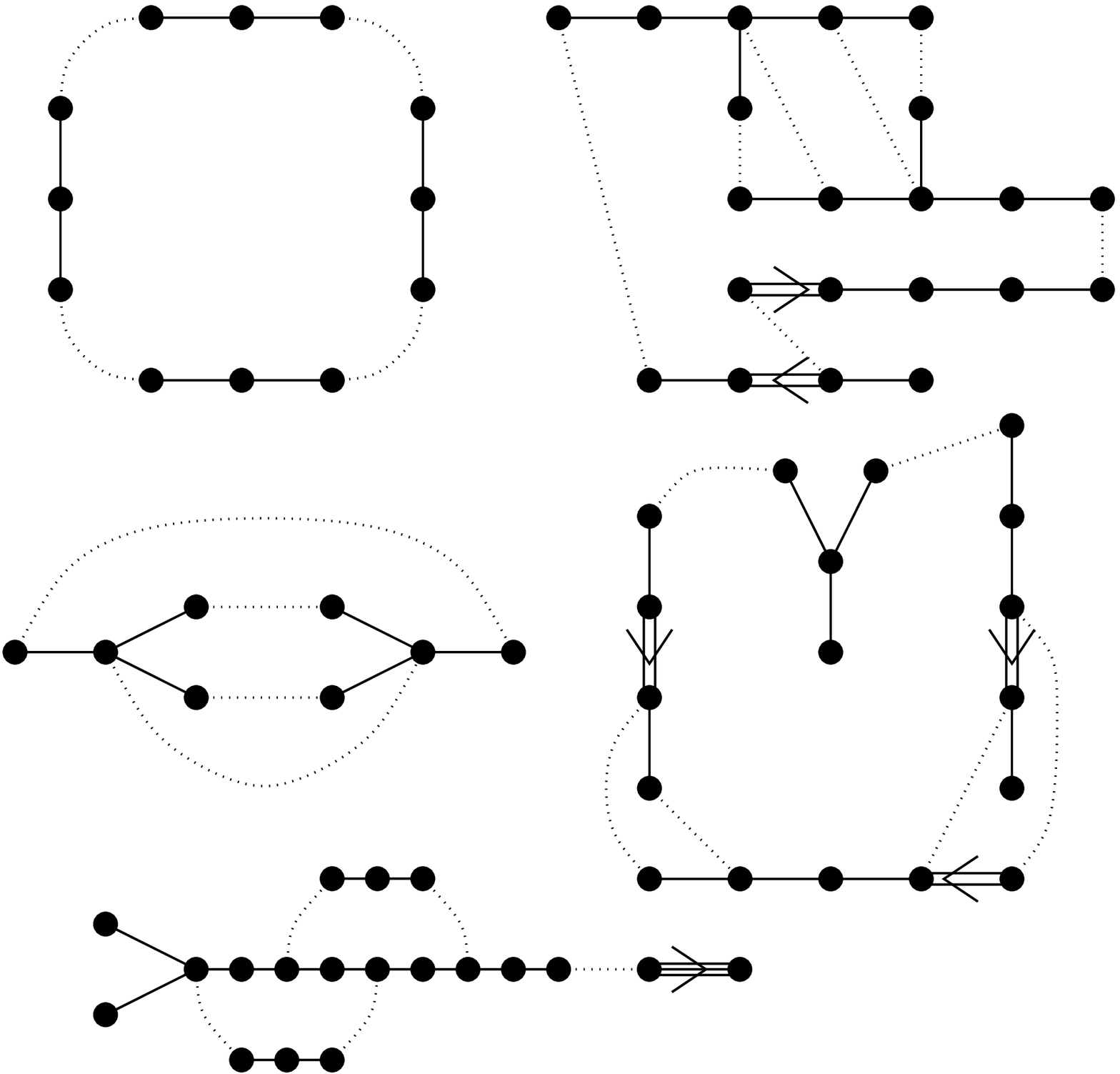}
\caption{Possible exotic linkings}
\label{exotics}
\vspace{4ex}
%\end{figure}
%\begin{figure}
%\centering
\includegraphics[scale=.7]{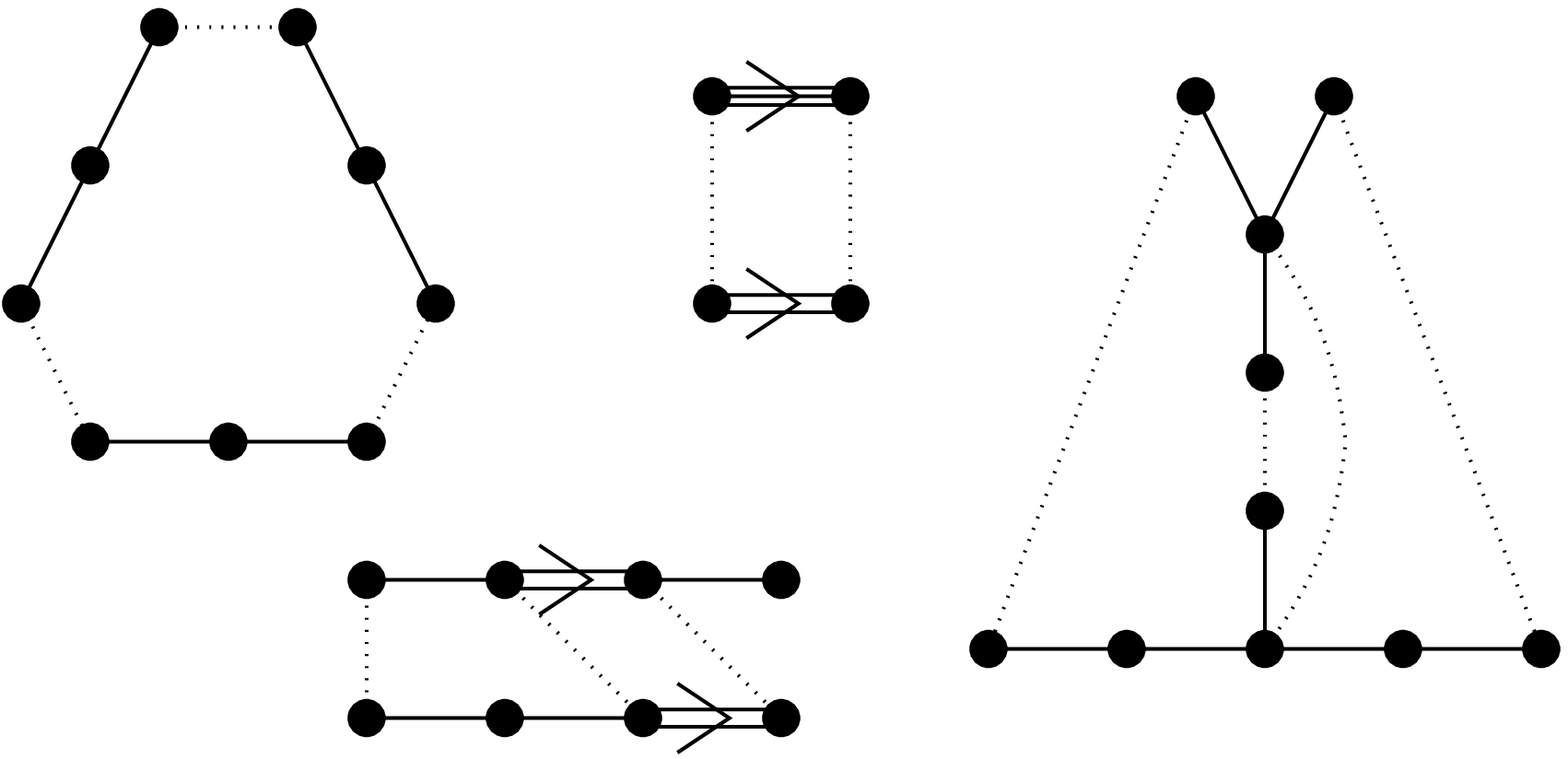}
\caption{Impossible exotic linkings}
\label{noways}
\end{figure}

In Figure \ref{exotics} we give a few more ``exotic'' linkings, that do not impose restrictions on the possible braiding matrices. The original example of 4 copies of $A_3$ linked into a circle is also depicted there. Some ``impossible'' linkable Dynkin diagrams, that is diagrams for which there are no linkable braiding matrices with diagonal elements not equal to one, are displayed in Figure \ref{noways}.
%\begin{minipage}{\textwidth}
%\end{minipage}

%%%%%%%%%%%%%%%%%%%%%%%%%%%%%%%%%%%%%%%%%%%%%%%%%%%%%%%%%%%%%%%%%%%%%%%%%%%

\section{Generalisations}
Here we want to discuss some possible generalisations to restrictions imposed so far. We only considered affine and finite Cartan matrices, as we do not know of any further special classes within generalized Cartan matrices. Apart from that, the combinatorics involved in the classification become more involved with increasing values of the entries of the Cartan matrix. In order to get a nice presentation, one would have to impose even more restrictions on the braiding matrix.

However, the results from Section \ref{affine-linkable} can be extended easily to all Dynkin diagrams with only single, double and triple edges. One possible application of this is for the Cartan matrices given in \cite[(4.9)-(4.17)]{AS-cartan}. As linkable Dynkin diagrams they all have one cycle of length zero. We give the genera of these cycles in Table \ref{except-cartan}, and see that the primes listed with these matrices in the original paper are just the odd divisors of these cycle genera. It is still an open question if finite dimensional Nichols algebras whose braiding matrix is of such Cartan types exist.
\begin{table}\centering
%\begin{minipage}{.5\textwidth}\centering%
\begin{tabular}{c|c|c}
Equation Number in \cite{AS-cartan}&Dynkin diagram&Genus\\
\hline
\extab{9}{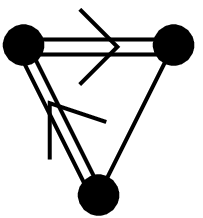}{2^2-1=3}
\extab{10}{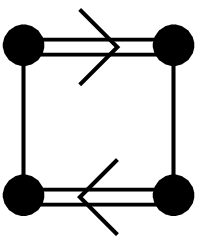}{2^2-1=3}
\extab{11}{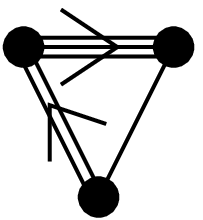}{3\cdot 2-1=5}
\extab{12}{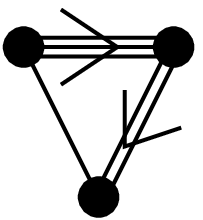}{3\cdot 2-1=5}
\extab{13}{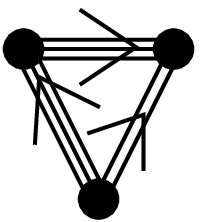}{3^2-2=7}
\extab{14}{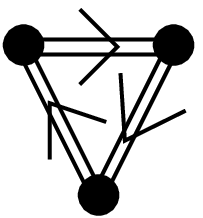}{2^3-1=7}
\extab{15}{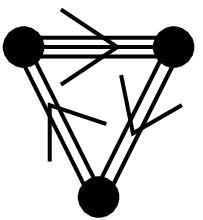}{3\cdot 2^2-1=11}
\extab{16}{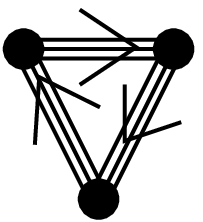}{3^3-1=2\cdot 13}
\extab{17}{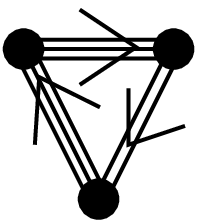}{3^2\cdot 2-1=17}
\hline
\end{tabular}%
%%\end{minipage}\begin{minipage}{.5\textwidth}\centering%
%\begin{tabular}{c|c|c}
%Eq. No&Diagram&Genus\\
%\hline
%%\extab{9}{pic49}{2^2-1=3}
%%\extab{10}{pic410}{2^2-1=3}
%%\extab{11}{pic411}{3\cdot 2-1=5}
%%\extab{12}{pic412}{3\cdot 2-1=5}
%\extab{13}{pic413}{3^2-2=7}
%\extab{14}{pic414}{2^3-1=7}
%\extab{15}{pic415}{3\cdot 2^2-1=11}
%\extab{16}{pic416}{3^3-1=2\cdot 13}
%\extab{17}{pic417}{3^2\cdot 2-1=17}
%\hline
%\end{tabular}
%%\end{minipage}
\caption{Equations (4.9)-(4.17) of \cite{AS-cartan} as Dynkin diagrams}
\label{except-cartan}
\end{table}

%%%%%%%%%%%%%%%%%%%%%%%%%%%%%%%%%%%%%%%%%%%%%%%%%%%%%%%%%%%%%%%%%%%%%%%%%%%%%%%%

\subsection{The order of the diagonal elements}
We would like to make some comments on why we made the various restrictions above on the orders of the diagonal elements.

First, if $b_{ii}=1,$ then $g_i$ commutes with $x_i;$ a case of no interest. Moreover, \cite[Lemma 3.1.]{AS-p3} shows that this  can not emerge when dealing with finite dimensional Hopf algebras.

When the order of $b_{ii}$ can be two, we do not get Lemma \ref{aslemma}. So we would have to deal with a much more difficult structure of possible diagrams.

 Actually, if one disregards the sub-diagram
\abild{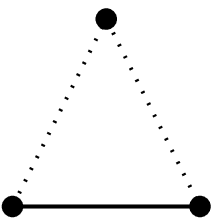}
then Lemma \ref{aslemma} already gives us that a vertex can be linkable to \emph{at most} one other vertex. One just takes $i=j$ and finds $a_{kl}=2$ or $a_{kl}=a_{lk}=-1$ and the order of $b_{ii}$ is 3. If the order of $b_{ii}$ is two, then one could consider a great number of diagrams with vertex $i$ being linkable to more than one vertex, thus making a classification even harder.

The infinite families of 32-dimensional pointed Hopf algebras found in \cite{Gn-32} all arise from the linkable Dynkin diagram where two copies of $A_1$ are linked to each other, and all the entries of the braiding matrix are $-1$.

To simplify the presentation of the theorems we excluded orders divisible by 3 when there are components of type $G_2.$ Without this limitation, a much more thorough examination of the diagrams (with heavy use of the above defined heights) is needed to establish a necessary condition for the existence of a braiding matrix. The problems come from trying to avoid diagonal elements of order 1 and 2.

For the affine diagrams this problem is even more severe, and the easiest way to avoid it is to consider only homogeneous braiding matrices with prime order of diagonal elements greater than 3. This way, there are no extra difficulties stemming from triple arrows either.

If we are dealing with a situation where the diagonal elements of the braiding matrix are required \emph{not} to be  roots of unity, e.g. in \cite{AS-char}, then Lemma \ref{tech} forces us to consider only diagrams where all cycle genera are zero. Hence, in this case Theorem \ref{mainthm} simplifies in the following way: In condition \ref{cd2} we do not need the last diagram as it has genus 3, and we need only the first part of condition \ref{cd3}. Actually, this first part can be reformulated,
% into the condition that in every cycle the number of dotted edges is even and the number of clockwise and anticlockwise pointing double edges is the same. 
thus avoiding the use of cycle genera, lengths and weights altogether. We get the following description.
\beer{thm}{
We are given a link-connected linkable Dynkin diagram $D$ of finite Cartan type, which is not $G_2\xdig G_2.$ A linkable braiding matrix of $D$-Cartan type in which the diagonal entries are not roots of unity exists, iff
\begin{enumerate}
\item In components of type $G_2$ not both vertices are linkable to other vertices.
\item $D$ does not contain any induced subgraphs of the form:
\abild{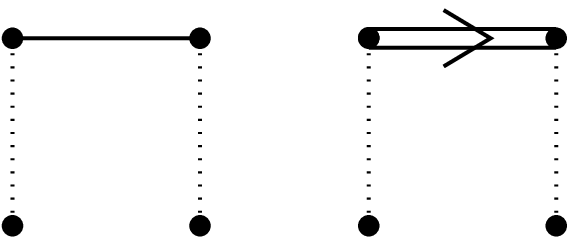}
\item For every cycle in $D$ we have:
        \begin{itemize} 
                \item The number of dotted edges is even.
                \item The number of clockwise and anticlockwise pointing double edges is the same.
        \end{itemize}
\end{enumerate}
}

%%%%%%%%%%%%%%%%%%%%%%%%%%%%%%%%%%%%%%%%%%%%%%%%%%%%%%%%%%%%%%%%%%%%%%%%%%%%%%%%

\subsection{Self-linkings}
So far we have always required that linkable vertices are in different connected components. This is part of the original definition of Andruskiewitsch and Schneider. But the generalization is straightforward and gives us a better understanding of the nature of some limitations required for the order of the diagonal elements. 

We first consider the case where two vertices $i$ and $j$ within the same connected component are linkable, but not neighbouring, i.e. $\a =0.$ In this case we can actually still apply all considerations. In Dynkin diagrams of affine and finite type the possible values of the genera of cycles formed by this special linking can be easily calculated. The cycle contains only one dotted edge and we find as possible values for the genera
\begin{equation}
g_c=\left\{\begin{array}{l}2\\3\\4\\5\end{array}\right.\text{if the diagram between $i$ and $j$ }\begin{array}{l}\text{is of type $A,$ $C^{(1)}$ or $D^{(2)}$}\\\text{has one double edge}\\\text{has one triple edge}\\\text{is of type $A_{2k}^{(2)},\;k\ge 2.$}\end{array}\raisetag{3.5em}
\end{equation}

If $\a\neq0$ and we link $i$ and $j$ we get from (\ref{link}) by setting $k$ first to $i$ then $j$ and exchanging $i$ and $j$ the following identity
\begin{equation}\label{selfcond} b_{ii}^{\a a_{ji}-\a-a_{ji}}=1.\end{equation}
For the possible sub-diagrams
\bbild[.9]{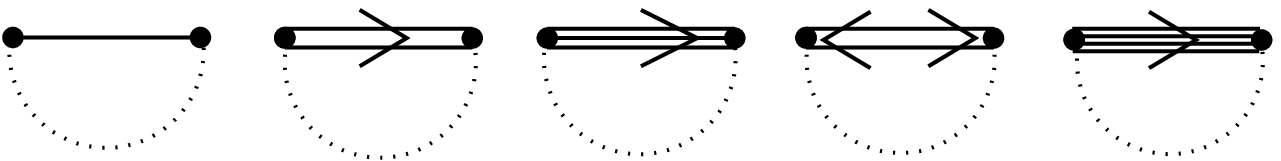}
the condition (\ref{selfcond}) gives us that the order of $b_{ii}$ must divide 3, 5, 7, 8 or 9, respectively.

So we see that admitting self-linkings restricts the possible braiding matrices drastically, especially if we are interested in homogeneous ones.

The values of 3 and 5 that we get when we link the two vertices in diagrams $A_2$ and $B_2$ respectively, are also the problematic ones when one tries to find all liftings of the Nichols algebras of this type, cf. \cite{AS-A2,BDR-B2}. 

%%%%%%%%%%%%%%%%%%%%%%%%%%%%%%%%%%%%%%%%%%%%%%%%%%%%%%%%%%%%%%%%%%%%%%%%%%%%%%%%

\subsection{Link-disconnected diagrams}
The results on link-connected diagrams can easily be extended to arbitrary diagrams. For each link-connected component the considerations can be carried out and a possible braiding matrix constructed. The direct sum of these matrices, with the zeros in the off-diagonal blocks replaced by ones, is then a braiding matrix for the whole diagram. If we are interested only in homogeneous braiding matrices, we have to check that we can choose the orders of the diagonal elements in all matrices (corresponding to the various link-connected components) to be the same.

%%%%%%%%%%%%%%%%%%%%%%%%%%%%%%%%%%%%%%%%%%%%%%%%%%%%%%%%%%%%%%%%%%%%%%%%%%%%%%%%

\section{Self-linkings in the rank 2 diagrams}\label{rank2comps}

We would like to present some results about Hopf algebras arising from self-linking the two vertices in the diagrams $A_2,$ $B_2$ and $G_2$. As the calculations are very involved and we have no insight into any simplifying structure theory, most of the assertions rely on results derived with Computer algebra programs. The explicit examples we construct can be used to generate infinitely many Hopf algebras of the same dimension with many parameters. As we have seen in the subsection discussing self-linkings in general, such linkings are only possible for very special orders $p$ of the diagonal elements of the braiding matrix. 
%%%%%%%%%%%%%%%%%%%%%%%%%%%%%%%%%%%%%%%%%%%%%%%%%%%%%%%%%%%%%%%%%%%%%%%%%%%%%%%%
\subsection[$A_2$ for $p=3$]{$\bm{A_2}$ for $\bm{p=3}$}\label{selfA2sec}

This is the case that was left open in \cite{AS-A2} and fully treated in \cite{BDR-B2}. We will derive the results from this paper in this relatively simple case again, so we can show how our method works. The notations are like in \cite{BDR-B2}.

So we are given the Cartan matrix $(\a)=(\begin{smallmatrix}2&-1\\-1&2\end{smallmatrix})$ and take as the linkable Dynkin diagram $A_2$ with the two vertices linkable to each other. For the corresponding braiding matrix we set $b_{11}:=q$ and get $b_{22}=q$ from (\ref{cartan}). Using (\ref{link}) we find 
$$b_{11}^2b_{12}=b_{21}^2b_{22}=b_{12}^2b_{11}=b_{22}^2b_{21}=1.$$
Thus $b_{12}=b_{21}=q^{-2}$ and $q^3=1.$ Therefore, the only possible matrix is $(b_{ij})=(\begin{smallmatrix}q&q\\q&q\end{smallmatrix}),$ where $q$ is a primitive root of unity of order 3. In \cite{AS-A2} this order was denoted by $p$ and that is the reason for us to call this the $p=3$ case.

Assume now, that there is an abelian group $\G$ with elements $g_1,$ $g_2$ and characters $\chi_1,$ $\chi_2$ such that $b_{ij}=\chi_j(g_i)$ and $\chi_1^2\chi_2=\chi_2^2\chi_1=1$. Take $\lambda_{12}$ and $\lambda_{21}$ arbitrary and consider the Hopf algebra $\Uf(A_2):=\Uf(\D)$ for this linking datum $\D$ as in Definition \ref{link.alg}. As in \cite[Section 3]{BDR-B2}, we introduce the root vector $z:=x_1x_2-b_{12}x_2x_1$
and see that the $\ga_i$ used there are related to \emph{our} $\lam$ by $\ga_1=q\lambda_{12}$ and $\ga_2=-\lambda_{21}.$ Plugging the quantum Serre relations (\ref{serre})

\begin{align*}
x_1z-q^2zx_1&=\lambda_{12}(1-g_1^2g_2),\\ x_2(x_2x_1-qx_1x_2)-q^2(x_2x_1-qx_1x_2)x_2&=\lambda_{21}(1-g_2^2g_1),
\end{align*}
into a computer algebra program we get back a PBW-basis of $\Uf(A_2)$ of the form $z^ix_1^jx_2^kg$ with $i,j,k\ge 0$ and $g\in\G.$

For this calculation we use {\bf felix} \cite{felix}, a program capable of computations with non-commutative structures. The program is freely available at {\tt http://felix.hgb-leipzig.de/} for various platforms. Unfortunately, it is not being developed any further and the documentation is rather sparse. In Appendix \ref{listA2} we give a listing of the code we used.

Using the quantum binomial formula (\ref{qbinom}) with $x=x_i\ot 1$ and $y=g_i\ot x_i,$ we see immediately that $x_i^3$ is $(g_i^3,1)$-primitive. Hence, the relations
$x_i^3=\mu_i(1-g_i^3)$ form a Hopf ideal. In order to turn $\Uf(A_2)$ into a finite dimensional Hopf algebra, we would like to divide out the third power of $z$ as well. We have to ensure that the new relation will lead to a Hopf ideal. By adding appropriate elements to $z^3$ we find that $$v:=z^3+(1-q)^3\mu_1\mu_2(1-g_2^3)+(1-q)q\lambda_{12}\lambda_{21}(1-g_2^2g_1)$$ is $(g_1^3g_2^3,1)$-primitive, and hence requiring the relations $$v=\lambda(1-g_1^3g_2^3),\qquad x_i^3=\mu_i(1-g_i^3) $$ in $\Uf(A_2),$ we get a Hopf algebra $\uf(A_2)$ of dimension $27\cdot\ord(\G)$ with 5 parameters $\lambda, \mu_1, \mu_2, \lambda_{12}$ and $\lambda_{21}.$ However, by rescaling $x_i,$ $\lambda, \lambda_{12}$ and $\lambda_{21}$ appropriately, we can always set $\mu_i$ equal to 1 or 0. So effectively, we are left with 3 free parameters when $\mu_i\neq 0$.

We see that the root vector parameters $u_{\alp_i}:=\mu_i(1-g_i^3)$ and $u_{\alp_1+\alp_2}:=(1-q)^3\mu_1\mu_2(1-g_2^3)+(1-q)q\lambda_{12}\lambda_{21}(1-g_2^2g_1)$ are elements of the group algebra and because the $g_i$ $q$-commute with the generators $x_i,$ all the $u_{\alp}$ are central.

If we take $\G=\Z/(9)$ and a generator of $\G$ as $g_1=g_2$, we get a 3 parameter family of $3^5$-dimensional Hopf algebras. For $\G=(\Z/(3))^2$ we have $g_i^3=1$ and all parameters appear in a product with one factor equal to zero, so the resulting Hopf algebras are independent of the parameters.

%%%%%%%%%%%%%%%%%%%%%%%%%%%%%%%%%%%%%%%%%%%%%%%%%%%%%%%%%%%%%%%%%%%%%%%%%%%%%%%%
\subsection[$B_2$ for $p=5$]{$\bm{B_2}$ for $\bm{p=5}$}

We want to treat the Dynkin diagram $B_2$ with its two vertices being linkable. This is the missing case in \cite[Theorem 2.7]{BDR-B2}. We will keep the notation as close as possible to theirs. The Cartan matrix is $(\a)=(\begin{smallmatrix}2&-1\\-2&2\end{smallmatrix})$ and we set $q:=b_{22}.$ Using (\ref{cartan}) and (\ref{link}) we find that the braiding matrix is $(b_{ij})=(\begin{smallmatrix}q^2&q\\q^2&q\end{smallmatrix})$ and $q$ is a primitive $\pth 5$ root of unity. $\lambda_{12}$ and $\lambda_{21}$ can be chosen arbitrarily, and we consider $\Uf(B_2):=\Uf(\D)$ for this linking datum with an appropriately chosen group $\G$ and elements $g_i$.

\begin{figure}
%\vspace{-4ex}
\parbox{.3\textwidth}{\centering\includegraphics{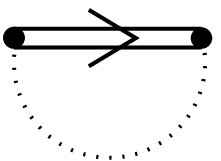}}\hfill $a_{ij}=\begin{pmatrix}2&-1\\-2&2\end{pmatrix}\quad b_{ij}=\chi_j(g_i)=\begin{pmatrix}q^2&q\\q^2&q\end{pmatrix}\quad q^5=1$\\

{\bf Generators: }$x_1, x_2, z, u$\\
{\bf Group elements: }$g_1, g_2,\qquad G_1:=g_1^2g_2, G_2:=g_2^3g_1$\\
{\bf Parameters: }$\lambda_{12}, \lambda_{21}, \mu_1, \mu_2, \ga_1, \ga_2$\\
{\bf Relations: }
\begin{align*}
x_2x_1&=q^2x_1x_2+z\\
x_2z&=q^3zx_2+u\\
x_1z&=qzx_1-q^2\lambda_{12}(1-G_1)\\
x_2u&=q^4ux_2+\lambda_{21}(1-G_2)\\
x_1u&=q^4ux_1+q^2(1-q)z^2+(1-q^3)\lambda_{12}x_2G_1\\
zu&=quz+q(1-q)(1-q^2)\lambda_{12}x_2^2G_1-q(1-q^2)\lambda_{21}x_1\\
x_i^5&=\mu_i(1-g_i^5)\\
v:&=z^5+\mu_1\mu_2(1-q^2)^5(1-g_1^5)\\
&\quad-(q^3-q^2-q+1)\lambda_{12}^2\lambda_{21}G_1^2G_2+(2q^3+2q^2+1)\lambda_{12}^2\lambda_{21}G_1^2\\
&\quad-q(q^2+3q+1)\lambda_{12}^2\lambda_{21}G_1+(q^3+2q^2+3q-1)\lambda_{12}\lambda_{21}zx_1G_1\\
&=\ga_1(1-(g_1g_2)^5)\\
w:&=u^5+2(1-q)^5\mu_2v-(1-q^2)^5(1-q)^5\mu_2^2x_1^5\\
&\quad-(3q^3+q^2-q+2)\lambda_{12}\lambda_{21}^3G_1G_2^3+(3q^3+q^2+4q+2)\lambda_{12}\lambda_{21}^3G_1G_2^2\\
&\quad-2(q^3+2q^2+3q-1)\lambda_{12}\lambda_{21}^3G_1G_2+(2q^3+4q^2+q-2)\lambda_{12}\lambda_{21}^3G_1\\
&\quad-5(q^3+1)\lambda_{12}\lambda_{21}^2ux_2G_1G_2-5(q^3+q^2-1)\lambda_{12}\lambda_{21}^2ux_2G_1\\
&\quad+5(q^3-1)\lambda_{12}\lambda_{21}u^2x_2^2G_1-5(2q^3-q^2+q-2)\lambda_{12}^2\lambda_{21}x_2^5G_1^2G_2\\
&=\ga_2(1-(g_1g_2^2)^5)
\end{align*}

\vspace{-63ex}
\hfill\parbox{10em}{
\begin{align*}
g_ix_2&=qx_2g_i\\
g_ix_1&=q^2x_1g_i\\
g_iz&=q^3zg_i\\
g_iu&=q^4ug_i
\end{align*}
}
\vspace{44ex}

{\bf Comultiplication: }
\begin{align*}
\del x_i&=g_i\ot x_i+x_i\ot 1\\%\qquad\del x_i^5=g_i^5\ot x_i^5+x_i^5\ot 1 \\
\del z&=g_1g_2\ot z+ z\ot 1+(1-q^3)x_2g_1\ot x_1\\
\del u&=g_1g_2^2\ot u+u\ot 1+(1-q^3)(1-q^4)x_2^2g_1\ot x_1+q(1-q^3)x_2g_1g_2\ot z\\
\del v&=(g_1g_2)^5\ot v+v\ot 1\\
\del w&=(g_1g_2^2)^5\ot w+w\ot 1
\end{align*}
%\end{minipage}
\caption{Relations for $B_2$}
\label{B2rel}
\end{figure}

Using the script in Appendix \ref{listB2}, we compute all the commutation relations between the root vectors $x_1, x_2, z:=x_2x_1-q^2x_1x_2$ and $u:=x_2z-q^3zx_2.$ This shows that $\{u^iz^jx_1^kx_2^lg : i,j,k,l\ge 0, g\in\G\}$ is a PBW-basis of $\Uf(B_2)$. $x_i^5$ is clearly $(g_i^5,1)$-primitive. Using felix and a separate module for tensor operations which was kindly provided by Istv\'an Heckenberger, we compute the $\pth 5$ powers of $\del(z)$ and $\del(u)$. The expressions returned by the program  are starting to be slightly unmanageable. After some reordering we get 18 and 87 terms respectively. Using Equation (14) in \cite{BDR-B2} and all the expressions from the felix-output ending on ``$\ot 1$'' we arrive at an educated guess for the right hand sides of the root vector relations for $z^5$ and $u^5,$ see page \pageref{B2rel}. Felix validates these guesses by showing that the resulting expressions $v$ and $w$ are indeed skew-primitive. Dividing out the $\pth 5$ powers of the root vectors with the appropriate right hand side we know that we get a Hopf algebra $\uf(B_2)$ of dimension $5^4\cdot\ord(\G).$ It is given explicitly in Figure \ref{B2rel}.

%\begin{minipage}{\textwidth}\label{B2rel}

As before, we can rescale the $x_i$ and all the 6 parameters appropriately, so that $\mu_i$ is either 1 or 0. Taking, for instance, $\G=\Z/(20)$ and $g_1=g_2$ as a generator of $\G$, we get a family of 12500-dimensional Hopf algebras with 4 parameters. Again, for some special values of the orders of the $g_i,$ the resulting Hopf algebras can become independent of some of the parameters.

As a final remark we note that the root vector parameters, i.e. the right hand sides, for $z^5$ and $u^5$ are not in the group algebra. Moreover, using felix we see that $v$ and $w$ are central and $z^5$ commutes with $x_1$ but not with $x_2$ and $u^5$ commutes with $x_2$ but not with $x_1.$ Hence, the counter terms are not central either if $\lambda_{12}\lambda_{21}\neq 0.$

%%%%%%%%%%%%%%%%%%%%%%%%%%%%%%%%%%%%%%%%%%%%%%%%%%%%%%%%%%%%%%%%%%%%%%%%%%%%%%%%
\subsection[$G_2$ for $p=7$]{$\bm{G_2}$ for $\bm{p=7}$}

\begin{figure}
%\vspace{-4ex}
\parbox{.3\textwidth}{\centering\includegraphics{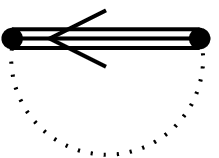}}\hfill $a_{ij}=\begin{pmatrix}2&-3\\-1&2\end{pmatrix}\quad b_{ij}=\chi_j(g_i)=\begin{pmatrix}q&q^3\\q&q^3\end{pmatrix}\quad q^7=1$\\

{\bf Generators: }$x_1, x_2, z, u, v, w$\\
{\bf Group elements: }$g_1, g_2,\qquad G_1:=g_1^4g_2, G_2:=g_2^2g_1$\\
{\bf Parameters: }$\lambda_{12}, \lambda_{21}$\\%, \mu_1, \mu_2, \ga_1, \ga_2$\\
{\bf Relations: }
\begin{align*}
x_2x_1&=qx_1x_2+z\\
x_2z&=q^4zx_2+\lambda_{21}(1-G_2)\\
x_1z&=q^5zx_1-q^5u\\
x_2u&=q^5ux_2+(q^4-q^2)z^2-(q^3-1)\lambda_{21}x_1\\
x_1u&=q^4ux_1-q^4v\\
zu&=q^3uz+w\\
x_2v&=q^6vx_2-(q^4-1)uz+(q^4-q^3-q^2)w-(q^4+q^3-2)\lambda_{21}x_1^2\\
x_1v&=q^3vx_1-q^6\lambda_{12}(1-G_1)\\
zv&=q^5vz+(2q^5+q^4+q^3+q^2+1)u^2-(q^5-q^4+3q^3+q^2+2q+1)\lambda_{21}x_1^3\\&\quad+(2q^5+2q^3+q+2)\lambda_{12}x_2G_1\\
x_2w&=q^2wx_2+(2q^5+q^4+q^3+q^2+2q)z^3+(2q^4-q-1)\lambda_{21}zx_1-(q^4-1)\lambda_{21}u\\
x_1w&=q^2wx_1+(q^3-q^2+q-1)u^2+(4q^5+2q^4+3q^3+2q^2+q+2)\lambda_{21}x_1^3\\&\quad-(2q^5+q^3+2q^2+2)\lambda_{12}x_2G_1\\
uv&=q^4vu+(3q^5+4q^4+3q^3+2q+2)\lambda_{21}x_1^4-(q^5-q^2)\lambda_{12}zG_1\\
zw&=q^4wz+(2q^4-q-1)\lambda_{21}ux_1-(q^4-1)\lambda_{21}v+(2q^5-q^2-q)\lambda_{12}x_2^2G_1\\
uw&=q^3wu+(q^5+2q^4+4q^3+4q^2+2q+1)\lambda_{21}vx_1+(q^4-q^3-q^2+q)\lambda_{12}zx_2G_1\\&\quad-(3q^5+4q^4+3q^3+2q^2+q+1)\lambda_{12}\lambda_{21}G_1G_2\\&\quad+(q^5+2q^4+2q^3-q^2-3q-1)\lambda_{12}\lambda_{21}G_1\\&\quad+(2q^5+2q^4+q^3+3q^2+4q+2)\lambda_{12}\lambda_{21}\\
vw&=q^5wv-(q^5+q^4+3q^2-q+3)u^3+(5q^5+2q^4+5q^3+q+1)\lambda_{21}ux_1^3\\&\quad+(2q^5-5q^3-5q^2-6q)\lambda_{21}vx_1^2-(3q^5+2q^4+2q^3+3q^2+4)\lambda_{12}ux_2G_1\\&\quad+(q^4-2q^3+q^2)\lambda_{12}z^2G_1+(5q^5+2q^4-4q^2-4q+1)\lambda_{12}\lambda_{21}x_1G_1\\&\quad-(q^5+q^4+q^3+2q^2+5q+4)\lambda_{12}\lambda_{21}x_1
\end{align*}

\vspace{-83ex}
\hfill\parbox{10em}{
\begin{align*}
g_ix_1&=qx_1g_i\\
g_ix_2&=q^3x_2g_i\\
g_iz&=q^4zg_i\\
g_iu&=q^5ug_i\\
g_iv&=q^6vg_i\\
g_iw&=q^2wg_i
\end{align*}
}
\vspace{60ex}

%\end{minipage}
\caption{Relations for $G_2$}
\label{G2rel}
\end{figure}

Here we want to treat the remaining rank 2 case of the finite Dynkin diagrams. Starting from the diagram $G_2$ with $(\a)=(\begin{smallmatrix}2&-3\\-1&2\end{smallmatrix})$ and the two vertices linkable, we find, as before, that the braiding matrix must be $(b_{ij})=(\begin{smallmatrix}q&q^3\\q&q^3\end{smallmatrix})$ where $q$ is a primitive $\pth 7$ root of unity. We define the root vectors $z:=x_2x_1-qx_1x_2,$ $u:=zx_1-q^2x_1z,$ $v:=ux_1-q^3x_1u$ and $w:=zu+q^3uz,$ cf. Equations (\ref{g2rootz})-(\ref{g2rootw}). The quantum Serre relations can now be expressed as $x_1v-q^3vx_1=-q^6\lambda_{12}(1-G_1)$ and $x_2z-q^4zx_2=\lambda_{21}(1-G_2)$ with $G_1:=g_1^4g_2$ and $G_2:=g_2^2g_1.$ Plugging this into felix as in Appendix \ref{listG2}, we get all commutation relations and hence a PBW-basis. The results are listed in Figure \ref{G2rel}.% and reduce to [ORIGINAL LUSZTIG/RINGEL] if lam12lam21=0.

%\begin{minipage}{\textwidth}\label{G2rel}

To be able to divide out the $\pth 7$ powers of all the root vectors, we again have to find appropriate counter terms first. Using felix to compute $\del(z^7)$ and trying the analogous treatment as for $B_2$ we indeed immediately get that
\begin{align*}
Z:=z^7&+(1-q^3)^7\mu_1\mu_2(1-g_1^7)\\
&-(2q^5+4q^4-q^3+q^2-4q-2)\l_{12}\l_{21}z^2x_2^2G_1\\
&+(6q^5+8q^4+6q^3-3q-3)\l_{12}\l_{21}^2zx_2G_1G_2\\
&-(q^4+3q^3-q^2+3q+1)\l_{12}\l_{21}^2zx_2G_1\\
&+(2q^5+2q^4+4q^3+5q^2+2q-1)\l_{12}\l_{21}^3G_1G_2^3\\
&+(q^5-2q^4-4q^3-7q^2-6q-3)\l_{12}\l_{21}^3G_1G_2^2\\
&-(4q^5+2q^4+2q^3-2q^2-2q-4)\l_{12}\l_{21}^3G_1G_2\\
&+(q^5+2q^4+2q^3+2q)\l_{12}\l_{21}^3G_1
\end{align*}
is $(g_1^7g_2^7,1)$-primitive and hence we can divide out $Z=\ga_1(1-g_1^7g_2^7).$ Trying to do the same for the root vector $u$ we find that the calculation powers and memory space needed are immense. The expression for $\del(u^7)$ has more than a thousand terms and needs some tricks in order for felix to be able to calculate it at all. Based on the experience gained so far, we can conclude from the felix-output that $u^7$ needs 46 counterterms to become skew-primitive. We do not list these here. For $v^7$ and $w^7$ we did not even attempt an answer. As before $z^7$ is not central and so neither are the counter terms.  

%%% Local Variables: 
%%% mode: latex
%%% TeX-master: "main.tex"
%%% End: 

\chapter{Group realization}\label{chap.group}

In this chapter we want to address the problem of deciding  which Dynkin diagrams of finite dimensional semisimple Lie algebras can be realized over a given group. This means that we have to find group elements and characters so that the braiding matrix associated with these elements is of the given Cartan type.  We are still a long way from dealing with this question in the general case and this chapter is meant more as a first step and example towards a deeper theory.

We will present results only for the groups $(\Z/(p))^2$ where $p$ is a prime bigger than 3. This will help us to describe the Hopf algebras of Theorem \ref{p17thm} even more explicitly in the case $s=2$. The case $s=1$ was done in \cite[Section 5]{AS-cartan}. One of the results established there was that the diagrams $A_2,$ $B_2$ and $G_2$ are only realizable over $\Z/(p)$ with $p\ge 5$ if $-3,$ $-1$ or $-3$ respectively, is a square modulo $p$. In Section 8 of that paper the authors prove some very nice results that place bounds on the number of vertices a Dynkin diagram can have if it is realizable over a group. For instance, for the groups $(\Z/(p))^s$ the maximal number of vertices is bound by $2s$ unless the Dynkin diagram has components of type $A_{p-1}$. In this case the bound can be raised by the number of such components. This estimate uses just simple linear algebra and the explicit values of the determinants of the Cartan matrices.

The results in this chapter will show that the bounds established in \cite[Propositon 8.3]{AS-cartan} are not only necessary, but also sufficient in most cases when $s=2$. This raises the hope that there is an $s_0>2$ such that for all $s\ge s_0,$ all Dynkin diagrams with no more than $2s$ vertices are realizable over $(\Z/(p))^s.$

We consider from now on $\G:=(\Z/(p))^2$ with $p\ge 5$ a prime and denote the group operation multiplicatively. From the bounds mentioned above we know that all diagrams realizable over $\G$ have at most 4 vertices, unless $p=5$ when the diagram could as well be $A_4\xdig A_1$. We prove that apart from a few exceptions the converse is also true.

\beer{thm}{\label{groupthm}All Dynkin diagrams of finite type with at most 4 vertices are realizable over $\G=(\Z/(p))^2$ for all $p\ge 5,$ except 
\begin{itemize}
\item $A_4$ is only realizable when $p\equiv 1\text{ or }9 \mod 10,$
\item $A_2\xdig B_2$ and $B_2\xdig G_2$ are only realizable when $p=12d\pm 1, d\ge 1.$
\end{itemize}
}
\prf{
We will prove this theorem by giving a nearly explicit description of the realization, i.e. we present a way of determining the necessary group elements and characters. The idea is the same as in \cite[Section 5]{AS-cartan} but the calculations are much more involved. First, we will show that all diagrams with 4 vertices except $A_4$, $A_2\xdig B_2$ and $B_2\xdig G_2$ are realizable over $\G$. All diagrams with less than 4 vertices are then realizable over $\G$ as well, because they are subdiagrams of the diagrams for which we established the result already.

To keep the presentation short we will not deal with each case separately but consider many diagrams simultaneously. However, we will not be too general, so that readability is still preserved.

{\bf Part (a)\hspace{2em}}\addcontentsline{toc}{section}{Proof Part (a)}The first and biggest class of diagrams we consider is 
\begin{align*}A_4, B_4, C_4, F_4,& A_3\xdig A_1, B_3\xdig A_1, C_3\xdig A_1,\\& A_2\xdig A_2, A_2\xdig B_2, A_2\xdig G_2, A_2\xdig A_1\xdig A_1.\end{align*}
We fix the enumeration of the vertices so that it follows the way we wrote the diagrams intuitively. With this we mean that the arrows are drawn in the directions given by Figure \ref{fig-lie}, and the order of the connected components is the same as in the written notation. The double and triple edges are always last in each component. As an example we make this precise for the diagram $C_3\xdig A_1$. The first two vertices form an $A_2$ subdiagram. The arrow points from vertex 3 to vertex 2 and vertex 4 corresponds to the $A_1$ component. The characteristic feature of this class of diagrams is that the first two vertices always form an $A_2$ subdiagram and $a_{24}=0$.

Now, assume that the diagrams are realizable over $\G,$ i.e. there are $g_1,\dots g_4\in\G$ and characters $\chi_1,\dots\chi_4\in\Gh$ so that the braiding matrix $(\chi_j(g_i))_{ij}$ is of the corresponding finite Cartan type. We assume that $g_1$ and $g_2$ generate $\G$. Then $g_3=g_1^ng_2^m$ and $g_4=g_1^kg_2^l$ for some $0\le n,m,k,l<p.$ We set $q:=\chi_1(g_1)$ and know that $q^p=1.$ Hence, $\chi_1(g_2)=q^t$ for some $0\le t<p.$ We can now express $\chi_1(g_3)$ and $\chi_1(g_4)$ in terms of $q,t,m,n,k$ and $l$. Using the Cartan condition (\ref{cartcond})
\begin{equation}\label{localcart}\chi_i(g_j)\chi_j(g_i)=\chi_i(g_i)^{\a}=\chi_j(g_j)^{a_{ji}}\end{equation}
we can also write down $\chi_2(g_1),$ $\chi_3(g_1)$ and $\chi_4(g_1)$ explicitly. As the first two vertices form an $A_2$ diagram and the third and fourth vertex are not connected to the first vertex directly in all the diagrams of our class, all the braiding matrix entries determined so far are the same for all diagrams. $\chi_2(g_2)=q$ because of (\ref{localcart}), and this lets us express $\chi_2(g_3)$ and $\chi_2(g_4)$. Using (\ref{localcart}) again, we now find $\chi_3(g_2)$ and $\chi_4(g_2)$. With those values, all remaining entries of the braiding matrix can be calculated. As a final result we get $(b_{ij})=(\chi_j(g_i))=$
%\begin{equation}
%\label{class1braid}
$$\begin{pmatrix}%{llll}
q & q^{-t-1} & q^{-n-tm} & q^{-k-tl}\\
q^t & q & q^{tn+n-m+a_{23}} & q^{tk+k-l}\\
q^{n+tm} & q^{-tn-n+m} & q^{-n^2+nm-m^2+a_{23}m} & q^{-kn-tln+tkm+km-lm}\\
q^{k+tl} & q^{-tk-k-l} & q^{-nk-tmk+tnl+nl-ml+a_{23}l} & q^{-k^2+kl-l^2}
\end{pmatrix}.$$
%\end{equation}
We note that the parameter $t$ is cancelled in the expressions for the diagonal elements. Some of the entries can be calculated using (\ref{localcart}) in a different way which leads to three important compatibility conditions.

First, $\chi_3(g_3)$ can be determined from $\chi_2(g_2)$ directly by the last part of (\ref{localcart}). We find $\chi_3(g_3)=q^x$ where $x:=a_{23}/a_{32}$ unless $a_{32}=0,$ as is the case for the last four diagrams of our class. Then $x$ is not further determined. Comparing this with the expression found so far, we get the condition (\ref{nm}). Next, we can determine $\chi_4(g_4)$ directly from $\chi_3(g_3)$, when vertices 3 and 4 are connected. In this case $\chi_4(g_4)=q^z$ with $z:=xa_{34}/a_{43}.$ When $a_{43}=0$ we again leave $z$ not further specified. This time the compatibility between the different ways of calculating the matrix elements leads us to (\ref{kl}). Finally, we can apply (\ref{localcart}) to find that $\chi_3(g_3)\chi_4(g_4)=q^{xa_{34}}$. Comparing this to the expressions found so far gives us (\ref{nmkl}). These are \emph{all} relations.

We list the possible values for the above parameters for each diagram explicitly in Table \ref{para-tab}. Let us now look at the system of equations we have found.
\begin{align}
\label{nm} -n^2+nm-m^2+a_{23}m&=x\quad\mod p\\
\label{kl} -k^2+kl-l^2&=z\quad\mod p\\
\label{nmkl} k(m-2n)+l(n-2m+a_{23})&=xa_{34}\quad\mod p
\end{align}
From now on all calculations will be made modulo $p$.

A linear transformation helps to bring the above equations into a nicer form. We set
$$a:=\frac{m+n-a_{23}}{2},\qquad b:=\frac{3m-3n-a_{23}}{2}\quad\text{and}\quad y:=3x-a_{23}^2.$$
Then (\ref{nm})-(\ref{nmkl}) can be rewritten as
\begin{align}
\label{ab} 3a^2+b^2+y&=0,\\
\label{newkl} k^2-kl+l^2+z&=0,\\
\label{abkl} k(a-b)+l(a+b)+xa_{34}&=0.
\end{align}
The first thing to notice is that (\ref{ab}) always has a solution\footnote{We thank Richard Pink, Bernd Herzog and Istv\'an Heckenberger for some useful tips concerning this equation.} . This can easily be seen by the following argument. As $a$ runs through the numbers 0 till $p-1,$ the expression $3a^2$ takes $\frac{p+1}2$ different values. The expression $-b^2-y$ has also $\frac{p+1}2$ different values for all the possible values for $b$. Assume that all these values are different. This means we have altogether $2\frac{p+1}2=p+1$ values. But modulo $p$ there are only $p$ numbers. Hence there is at least one solution. Actually, there are even $6d$ solutions for every prime $p$ of the form $p=6d\pm 1$. This can be seen by considering the homogeneous equation
\begin{equation}\label{homo3}3\tilde a^2+\tilde b^2+yc^2=0,\end{equation} which we get by multiplying (\ref{ab}) by $c^2$ and setting $\tilde a=ac$ and $\tilde b=bc.$ For such equations there is a nice theory, and we can get the number of solutions from \cite[Theorem 2, p.25]{BS}. The Gaussian sums appearing there are zero in our case and the number of solutions of (\ref{homo3}) is just $p^2$ for all primes $p$. If $-3$ is not a square modulo $p,$ then there is only one solution (the trivial one) with $c=0$. This solution does not lead to a solution of the original problem (\ref{ab}), and all other solutions produce $p-1$ times the same solution for the original equation, as $c$ acts like a scaling. Hence, when $-3$ is not a square (\ref{ab}) has $\frac{p^2-1}{p-1}=p+1$ solutions. When $-3$ is a square, there are $2(p-1)+1$ solutions of (\ref{homo3}) with $c=0,$ and therefore, $p^2-2p+1$ solutions which give a solution of (\ref{ab}). Again, due to the scaling, $p-1$ of the solutions lead to the same values for $a$ and $b$. Hence, if $-3$ is a square (\ref{ab}) has $\frac{p^2-2p+1}{p-1}=p-1$ solutions. Using the quadratic reciprocity law we find that $-3$ is a square only for primes $p$ of the form $p=6d+1.$

As $p\ge 5$ we know that (\ref{ab}) has at least 6 solutions. Therefore, we can always pick a solution with $a\neq b$, as this corresponds to 2 solutions at most. Then we can solve (\ref{abkl}) for $k$ and substitute the result in (\ref{newkl}). After some reordering we get
\begin{align*} \left[\left(\frac{a+b}{a-b}\right)^2+\frac{a+b}{a-b}+1\right]l^2&+\left[\frac{2(a+b)xa_{34}}{(a-b)^2}+\frac{xa_{34}}{a-b}\right]l+\\
&+\left[\left(\frac{xa_{34}}{a-b}\right)^2+z\right]=0. 
\end{align*}
Multiplying by the common denominator, simplifying and using (\ref{ab}) we arrive at 
$$yl^2-(3a+b)xa_{34}l-(x^2a_{34}^2+z(a-b)^2)=0.$$
For the discriminant of this quadratic equation we calculate 
$$ \left(\frac{(3a+b)xa_{34}}{2y}\right)^2+\frac{x^2a_{34}^2+z(a-b)^2}y=\left(\frac{a-b}{2y}\right)^2\left[4yz-3x^2a_{34}^2\right].$$
Hence, our system of equations has a solution modulo $p$ if and only if $D:=4yz-3x^2a_{34}^2$ is a square modulo $p$. All the values of the parameters for our class of Dynkin diagrams are given in Table \ref{para-tab}.
\begin{table}
\begin{center}
\begin{tabular}{c|c c c c c|c}
\dtab{\text{Diagram}}{a_{23}}{a_{34}}xyz{D=4yz-3x^2a_{34}^2}
%Diagram&$a_{23}$&$a_{34}$&$x$&$y$&$z$&$D$\\
\hline
\dtab{A_4}{-1}{-1}121{8-3=5}
\dtab{B_4}{-1}{-1}12{1/2}{4-3=1^2}
\dtab{C_4}{-1}{-2}122{16-12=2^2}
\dtab{F_4}{-2}{-1}222{16-12=2^2}
\dtab{A_3\xdig A_1}{-1}012z{8z}
\dtab{B_3\xdig A_1}{-1}0{1/2}{1/2}z{2z}
\dtab{C_3\xdig A_1}{-2}022z{8z}
\dtab{A_2\xdig A_2}0{-1}x{3x}x{12x^2-3x^2=(3x)^2}
\dtab{A_2\xdig B_2}0{-1}x{3x}{x/2}{6x^2-3x^2=3x^2}
\dtab{A_2\xdig G_2}0{-1}x{3x}{x/3}{4x^2-3x^2=x^2}
\dtab{A_2\xdig A_1\xdig A_1}00x{3x}z{12xz}
\hline
\end{tabular}
\end{center}
\caption{Parameters for Dynkin diagrams}
\label{para-tab}
\end{table}
In the cases where $a_{23}=0$ or $a_{34}=0$ we can freely choose the values of $x$ or $z$ respectively. So we see that $D$ is or can be made a square for all diagrams except for $A_4$ and $A_2\xdig B_2$. For these we need 5 or 3 respectively to be a square modulo $p$. According to the quadratic reciprocity law this happens if $p=\pm 1 \mod 10$ or $p=12d\pm 1,$ respectively. After choosing a solution of (\ref{ab}) with $a\neq b$ we can determine $k$ and $l$ and thus find an explicit realization of these diagrams over $\G$.

One could argue that the assumptions made at the beginning of this calculation are too special and other choices will make it possible to realize $A_4$ and $A_2\xdig B_2$ as well. So let us consider this more closely.

First we see that we can not have $a=b$ and $a=-b$ simultaneously, as (\ref{ab}) would give $y=0$. However, Table \ref{para-tab} shows that $y\neq0.$ We recall that $x$ and $z$ are not zero, because the diagonal elements of a braiding matrix of Cartan type are different from 1 (\ref{diagnot1orig}). Now we notice the following symmetry. (\ref{newkl}) does not change if $k$ and $l$ are interchanged. If we also send $b$ to $-b,$ we see that the whole system of equations (\ref{ab})-(\ref{abkl}) is unchanged. So the case $a=b\neq 0$ leads to basically the \emph{same} solution as $a\neq b$ if we set $a=-b.$ Therefore, the crucial discriminant is the same too and we get the same conditions.

We now only consider $A_4.$ We will come back to the diagram $A_2\xdig B_2$ when we deal with $B_2\xdig G_2$. So what if $g_1$ and $g_2$ do not generate $\G$? Then $g_2=g_1^a$ for some $a\in\Z/(p).$ We know from \cite[Proposition 5.1.]{AS-cartan} that $g_3$ can not be a power of $g_1$ as well, because then we would have a diagram with 3 vertices realizable over $\Z/(p)$ with $p\ge 5.$ If $g_4$ is not a power of $g_3$ then $g_3$ and $g_4$ generate $\G$. In this case we simply enumerate the vertices of $A_4$ from the other end and get back the original case, because $A_4$ is symmetric. But if $g_4=g_3^b$ for some $b\in\Z/(p)$, we argue in the following way. Set $q:=\chi_1(g_1).$  Then $\chi_1(g_3)=q^t$ and $\chi_1(g_4)=q^{tb}$ for some $t\in\Z/(p).$ Using (\ref{localcart}) we know that $\chi_3(g_1)=q^{-t}$ and $\chi_4(g_1)=q^{-tb}$ and hence $\chi_3(g_2)=q^{-at}$ and $\chi_4(g_2)=q^{-atb}.$ Again we use (\ref{localcart}) to calculate $\chi_2(g_3)=q^{at-1}$ and $\chi_2(g_4)=q^{atb}$. But as $g_4=g_3^b$ we see that $b$ must be zero. On the other hand, we have the condition $b^2+b+1=0$ because vertices 3 and 4 form an $A_2$ diagram, cf. \cite[page 25]{AS-cartan}. This is a contradiction. So we have fully proved that $A_4$ can be realized over $\G$ only if 5 is a square modulo $p$.

{\bf Part (b)\hspace{2em}}\addcontentsline{toc}{section}{Proof Part (b)} We consider most of the remaining 4-vertex diagrams:
$$B_2\xdig B_2, B_2\xdig G_2, G_2\xdig G_2, B_2\xdig A_1\xdig A_1, G_2\xdig A_1\xdig A_1,A_1\xdig A_1\xdig A_1\xdig A_1.$$
These diagrams are characterized by $a_{23}=0,$ the enumeration of the vertices is like before, but all the arrows are assumed to be pointing at vertices 1 or 3.

Again we assume that these diagrams are realizable over $\G$ and that $g_1$ and $g_2$ generate $\G$. We then carry out the completely analogous steps as in part (a). We introduce the parameters $x, y$ and $z$ to denote the exponents of the diagonal braiding matrix elements, so that
$$\chi_1(g_1)=q,\quad\chi_2(g_2)=q^x,\quad\chi_3(g_3)=q^y,\quad\chi_4(g_4)=q^z.$$
With $g_3=g_1^ng_2^m$ and $g_4=g_1^kg_2^l$ the three compatibility conditions come out as
\begin{align}
\label{nmpartb} n^2+a_{12}nm+xm^2+y&=0,\\
\label{klpartb} k^2+a_{12}kl+xl^2+z&=0,\\
\label{nmklpartb} k(2n+a_{12}m)+l(2xm+a_{12}n)+a_{34}y&=0.
\end{align}
As vertex 2 and 3 are not directly connected in all the diagrams of our class, we can choose $y\neq 0$ freely. Hence, we can always arrange that (\ref{nmpartb}) has a solution. Assume that $2n+a_{12}m\neq 0$. Then we can solve (\ref{nmklpartb}) for $k$ and plug the result into (\ref{klpartb}). After some reordering and simplifying we see that we can apply (\ref{nmpartb}) to get 
$$(4x-a_{12}^2)yl^2-ma_{34}(4x-a_{12}^2)yl-(a_{34}^2y^2+z(2n+a_{12}m)^2)=0.$$
For solving this quadratic equation we consider the discriminant
$$ \left(\frac{ma_{34}}2\right)^2+\frac{a_{34}^2y^2+z(2n+a_{12}m)^2}{(4x-a_{12}^2)y}=\left(\frac{2n+a_{12}m}2\right)^2\left[\frac{4z-a_{34}^2y}{(4x-a_{12}^2)y}\right]$$
and see that we have a solution to the system of equations (\ref{nmpartb})-(\ref{nmklpartb}) if $D:=\frac{4z-a_{34}^2y}{(4x-a_{12}^2)y}$ is a square modulo $p$. All the parameters are listed in Table \ref{para-tab2}.
\begin{table}
\begin{center}
\begin{tabular}{c|c c c c|c}
\dtabb{\text{Diagram}}{a_{12}}{a_{34}}xz{D=\frac{4z-a_{34}^2y}{(4x-a_{12}^2)y}}
\hline
\dtabb{B_2\xdig B_2}{-2}{-2}2{2y}{\frac{4y}{4y}=1}
\dtabb{B_2\xdig G_2}{-2}{-3}2{3y}{\frac{3y}{4y}=\frac 3{2^2}}
\dtabb{G_2\xdig G_2}{-3}{-3}3{3y}{\frac{3y}{3y}=1}
\dtabb{B_2\xdig A_1\xdig A_1}{-2}02z{\frac{4z}{4y}=\frac zy}
\dtabb{G_2\xdig A_1\xdig A_1}{-3}03z{\frac{4z}{3y}=2^2\frac z{3y}}
\dtabb{A_1\xdig A_1\xdig A_1\xdig A_1}00xz{\frac{4z}{4xy}=\frac z{xy}}
\hline
\end{tabular}
\end{center}
\caption{More parameters for Dynkin diagrams}
\label{para-tab2}
\end{table}
When $x$ and $z$ can be chosen freely, we can always arrange for $D$ to be a square. And so we see that all these diagrams can be realized over $\G$ except $B_2\xdig G_2$, where we need 3 to be a square modulo $p$.

Again we show that our assumptions are not too special. We only consider the problematic diagram $B_2\xdig G_2.$ If $2n+a_{12}m=0$ and $2xm+a_{12}n=0$, we find that $4x-a_{12}^2=0$. But that is not the case as we see from Table \ref{para-tab2}. So we assume $2n+a_{12}m=0$ and solve (\ref{nmklpartb}) for $l$. From (\ref{nmpartb}) we find $y=-(4x-a_{12}^2)m^2/4$ and so $l=a_{34}m/2.$ Plugging this into (\ref{klpartb}) we solve the quadratic equation and find that it has a solution if 3 is a square modulo $p$. Note that $z=3y,$ $x=2=-a_{12}$ and $a_{34}=-3.$

Now we will show that the case where $g_2$ is a power of $g_1$ does not lead to any simplification of the condition that 3 must be a square for $B_2\xdig G_2$ and $A_2\xdig B_2.$ If $g_3$ and $g_4$ generate $\G$, we just enumerate the vertices, so they correspond to the notation $G_2\xdig B_2$ and $B_2\xdig A_2,$ respectively. Now all considerations still apply and we find for the first diagram $D=4/3$ and the second $D=3/4.$ If however $g_4$ is a power of $g_3,$ we have the case that one copy of $\Z/(p)$ has to realize $B_2$ and the other one has to realize $A_2$ or $G_2$. This is only possible if $-1$ is a square ($B_2$) and $-3$ is square ($A_2$ and $G_2$), but then 3 is also a square.

{\bf Part (c)\hspace{2em}}\addcontentsline{toc}{section}{Proof Part (c)}The only diagram left is $D_4.$ We number the vertices so that $a_{12}=a_{23}=a_{24}=-1$. Then we proceed as in the last two parts. We get the three compatibility conditions
\begin{align*}
n^2-nm+m^2+m+1&=0,\\
k^2-kl+l^2+l+1&=0,\\
k(m-2n)+l(n-2m-1)-m&=0.
\end{align*}
The first two equations are the same and can be transformed into (\ref{ab}) in the same way as before. So we know that we can find $n$ and $m$ fulfilling the first equation, such that $m\neq 2n.$ Solving the last equation for $k$ and plugging this into the middle equation, we find that $l=-\frac{m+2\pm(m-2n)}2.$ No square roots are necessary and we have found a realization.
}

This proof was rather direct and does not explain in any way, \emph{why} the quadratic systems of 3 equations and 4 indeterminates can always be solved independently of the choices of the numerous parameters. It seems that this proof is just a shadow of a much deeper principle. Especially when one wants to generalize this kind of result to groups with more copies of $\Z/(p),$ a more fundamental understanding might be necessary. For instance, we believe that following the same strategy as we have just displayed, we can get a system of $\binom{t+1}2$ compatibility equations with $st$ indeterminates when we try to realize a diagram with $s+t$ vertices over $\G(s):=(\Z/(p))^s.$ As far as we understand this general case, the equations should remain quadratic. This should be fairly easy to prove. So judging from the case $s=2$ presented here, things should get simpler, as the quotient of indeterminates to equations rises quite comfortably. We could expect more freedom in the choice of the variables and hence hope that all diagrams with at most $2s$ vertices will be realizable over $\G(s).$ However, this reasoning, of course, completely fails to explain why we can not realize diagrams with $2s+1$ vertices in general. The only observation we can make is that at this step (from $s+s$ vertices to $s+s+1$ vertices), the number of new equations becomes bigger than the number of new variables. But even if all this can be made rigorous, what guarantees the solvability of the system of equations for all $p$? One last remark is that for Lusztig's finite dimensional Hopf algebras, two copies of the same diagram with $s$ vertices each is always realizable over $\G(s).$ Moreover, the two copies can even be linked.

We want to discuss possible linkings in the diagrams for $\G(2)$. We only consider linking vertices in different connected components. So $i$ and $j$ can be linked only if $a_{ij}=0$ and $\chi_i\chi_j=1.$ Using (\ref{localcart}) we find that $\chi_i(g_i)=\chi_j(g_j)^{-1}.$ Hence, when we link vertices, one of the free parameters (the exponent of a diagonal element) becomes fixed. From Tables \ref{para-tab} and \ref{para-tab2} we see that the only problematic linkings are the cases where the diagram is $A_2/B_2/G_2\xdig A_1\xdig A_1$ and we link the $A_1$ components or $A_3/B_3/C_3\xdig A_1$. In the first cases the crucial discriminant $D$ is only a square if $-3,$ $-1$ or $-3$ are squares modulo $p,$ respectively. These are \emph{exactly} the conditions necessary to realize $A_2,$ $B_2$ and $G_2$ over $\Z/(p)$. So assuming $g_1$ and $g_2$ to generate $\G$ is not limiting and we can link $A_1$ to $A_1$ in $A_2/G_2\xdig A_1\xdig A_1$ or $B_2\xdig A_1\xdig A_1$ only if $-3$ or $-1,$ respectively, is a square. This is only a necessary condition and might not be sufficient, as we still have to check $\chi_i\chi_j=1$ on the group elements.

In the diagrams $A_3/B_3/C_3\xdig A_1$ we have $z=-1$ if we link vertex 4 to vertex 1 or 2. If we link it to vertex 3, then $z=-1,$ $z=-1/2$ or $z=-2,$ respectively. From Table \ref{para-tab} we see that $2z$ needs to be a square. This means, we need $-2$ to be a square if $z=-1,$ or $-1$ to be a square for the other two cases. This does not change if we assume that $g_2$ is a power of $g_1.$ We would only get the added condition that $-3$ has to be a square, too. Again, this might not be sufficient, and the two extra conditions coming from evaluating $\chi_i\chi_j=1$ on the group generators have to be taken into consideration.

Finally, let us mention that for $p=5$ the diagram $A_4\xdig A_1$ is realizable over $(\Z/(p))^2$. 5 is of course a square modulo $p=5$. Therefore, part (a) of the proof gives us a realization of $A_4.$ We choose a special solution with $n=4,$ $m=2$ and $k=l=3.$ Setting $g_5:=g_1g_2^4$ and $q$ a primitive $\pth 5$ root of unity, we have this special realization of $A_4\xdig A_1:$
$$ (\chi_j(g_i))=
\begin{pmatrix}
q&q^4&q&q^2&q^4\\
1&q&q&1&q^2\\
q^4&q^3&q&q^3&1\\
q^3&1&q&q&q^3\\
q&q^3&1&q^2&q^2\\
\end{pmatrix}.
\qquad
\begin{bmatrix}
g_1\\g_2\\g_1^4g_2^2=g_3\\g_1^3g_2^3=g_4\\g_1g_2^4=g_5\\
\end{bmatrix}
$$

%%% Local Variables: 
%%% mode: latex
%%% TeX-master: "main"
%%% End: 

\chapter{Quasi-isomorphisms}

We will show that a large class of finite dimensional pointed Hopf algebras is quasi-isomorphic to their associated graded version coming from the coradical filtration, i.e.~they are $2$-cocycle deformations of the latter. This supports a slightly specialized form of a conjecture in \cite{Mas}. Most of the material in this chapter will be published as \cite{D2}.

Recently there has been a lot of progress in determining the structure of pointed Hopf algebras. This has led to a discovery of whole new classes of such Hopf algebras and to some important classification results. A lot of these classes contain infinitely many non-isomorphic Hopf algebras of the same dimension, thus disproving an old conjecture of Kaplansky. Masuoka showed in \cite{Mas} and in a private note \cite{Mas-priv}, that for certain of these new families the Hopf algebras are all $2$-cocycle deformations of each other. This led him to weaken Kaplansky's conjecture, stating that up to quasi-isomorphisms there should only be a finite number of Hopf algebras of a given dimension. This was disproved in \cite{EG} for families of Hopf algebras of dimension 32. However, our results support Masuoka's conjecture in a really big class of examples. So we suggest that the conjecture could still be saved by specializing it slightly. We have the feeling that there is a fundamental difference between the even and odd dimensional case. We propose for a base field of characteristic zero:
\begin{quote}\label{conjec} In a given odd dimension there are only finitely many non quasi-isomorphic pointed Hopf algebras whose coradical is abelian.\end{quote}
Two Hopf algebras are called \emph{quasi-isomorphic} or \emph{monoidally co-Morita equivalent} if their categories of comodules are monoidally equivalent. In \cite{S} Schauenburg showed that for finite dimensional Hopf algebras this is equivalent to the Hopf algebras being $2$-cocycle deformations of each other, cf. Section \ref{basics.cocyc}.

We fix a linking datum $\D$ of finite Cartan type with abelian group $\G$, cf. Definition \ref{link.datum} and the comments thereafter. We require that there are no self-linkings. Consider now a family $\uu$ of root vector parameters, such that $\uf(\D,\uu)$ is a pointed finite dimensional Hopf algebra. The ultimate goal is to show that all Hopf algebras which differ only in their choice of linking and root vector parameters are cocycle deformations of each other. Hence, there is only one quasi-isomorphism class. We present three major steps towards this goal.

First, we prove the statement in the case with only linking parameters and no root vector parameters, i.e. all $u_\alp$ are zero. In the second part we show the result for the case where the Cartan matrix is of type $A_n.$ Here, all the root vector parameters are known explicitly and linking parameters do not appear. In the last part we combine all results to treat the mixed case, where the Dynkin diagram is a union of $A_n$'s.

Once the root vector parameters for the other Cartan matrices of finite type have been determined explicitly, an analogous treatment should provide the same result as for the $A_n$ case.

\section{The linking case}

We assume that we have no root vector parameters, i.e. all $u_\alp$ are zero. This is the case, for instance, for the Hopf algebras of Theorem \ref{p17thm}. Instead of $\uf(\D,\bm 0)$ we simply write $\uf(\D)$.

\beer{thm}{\label{linkthm} For two linking data $\D$ and $\D'$ of finite Cartan type which only differ in their choice of the $\lam$ we have that $\uf(\D)$ and $\uf(\D')$ are quasi-isomorphic.
}

\prf{The strategy of the proof is as follows. We will prove the statement for an arbitrary linking datum and the datum where all the $\lam$ are 0.
This will be achieved inductively by proving the statement for an arbitrary linking datum with at least one pair of linked vertices and one datum with the same data, except that the number of connected components that are not linked to any other vertex is increased by one. This means, if in the original datum we have $\lam\neq 0$ for some $i,j,$ then we take for the other datum $\lambda_{kl}=0$ for all $1\le l\le\theta$ and for all $k$ that are in the same connected component $I$ as $i$. The transitivity of the quasi-isomorphism relation will then yield the result.

So let $\D=\{\G,(\a)_{1\le i,j\le\theta},(g_i)_{1\le i\le\theta},(\chi_j)_{1\le j\le\theta},(\lam)_{1\le i<j\le\theta}\} $ be a linking datum, where there is a connected component $I$ and an $i\in I$ such that there is a $j$ with $\lam\neq 0$. Let $\tilde{\D}$ be the same linking datum except that $\lambda_{kl}=0$ for all $k\in I$. It suffices to show that $\uf(\D)$ and $\uf(\tilde{\D})$ are quasi-isomorphic.

Now, from the proof of Theorem 5.17.~in \cite{AS-p17} we see that $\uf(\D)$ is isomorphic to $(\U\ot\B)_{\si}/\K^+,$ where $\U,\B$ and $\K$ are Hopf algebras to be described shortly. $\si$ is a 2-cocycle constructed from the linking datum. The key observation is now that $\uf(\tilde{\D})$ is isomorphic to $(\U\ot\B)_{\tilde{\si}}/\K^+$ with the \emph{same} Hopf algebras involved. Only the cocycle $\tilde{\si}$ is different. 

We reorder the vertices so that $I=\{1,\dots,\tilde{\theta}\}$ and set $\Y:={<\!Z_1\!>}\oplus\cdots\oplus {<\!Z_{\tth}\!>},$ where the order of $Z_i$ is the least common multiple of $\ord g_i$ and $\ord \chi_i.$ Let $\eta_j$ be the unique character of $\Y$ such that $\eta_j(Z_i)=\chi_j(g_i),$ $1\le i,j\le\tilde{\theta}.$ This is well defined because $\ord g_i$ divides $\ord Z_i$ for all $i.$

Now take as linking datum $\D_1$ with the original group $\G$ $$\D_1=\{\G,(\a)_{\tilde{\theta}<i,j\le\theta},(g_i)_{\tilde{\theta}<i\le\theta},(\chi_j)_{\tilde{\theta}<j\le\theta},(\lam)_{\tilde{\theta}<i<j\le\theta}\},$$ and as $\D_2$ with the group $\Y$
$$\D_2=\{\Y,(\a)_{1\le i,j\le\tilde{\theta}},(Z_i)_{1\le i\le\tilde{\theta}},(\eta_j)_{1\le j\le\tilde{\theta}},(\lam)_{1\le i<j\le\tilde{\theta}}\}.$$
Then $\B:=\uf(\D_1)$ and $\U:=\uf(\D_2).$ We denote the generators of $\B$ by $b_{\tth+1},\dots,b_{\theta}$ and $y_1,\dots,y_s,$ and the generators of $\U$ by $u_1,\dots,u_{\tth}$ and $z_1,\dots,z_{\tth}.$

The central Hopf subalgebra $\K$ of $\U\ot\B$ is $\k[z_i\ot g_i^{-1}:1\le i\le\tth].$

The cocycle $\si$ for $\U\ot\B$ is defined by \begin{equation*} \si(u\ot a,v\ot b):=\ep_{\U}(u)\tau(v,a)\ep_{\B}(b).\end{equation*}
$\tau: \U\ot\B\rightarrow\k$ is a linear map with the following list of properties;
\begin{enumerate}
\item $\tau(uv,a)=\tau(u,a\_1)\tau(v,a\_2),$
\item $\tau(u,ab)=\tau(u\_1,b)\tau(u\_2,a),$
\item $\tau(1,a)=\ep_{\B}(a),$
\item $\tau(u,1)=\ep_{\U}(u).$
\end{enumerate}
It is given by $$\tau(u\ot b):=\varphi(u)(b),$$ where $\varphi:\U\rightarrow(\B^*)^{cop}$ is a Hopf algebra homomorphism defined on the generators of $\U$ by
$$\varphi(z_i):=\ga_i\qquad\text{and}\qquad\varphi(u_j):=\de_j.$$
Here $\ga_i:\B\rightarrow\k$ is a character defined on the generators of $\B$ by
$$ \ga_i(y_k):=\chi_i(y_k)\qquad\text{and}\qquad\ga_i(b_l):=0,$$
and $\de_j:\B\rightarrow\k$ is a $(\ep_{\B},\ga_j)$-derivation defined by
$$ \de_j(y_k):=0\qquad\text{and}\qquad\de_j(b_l):=-\chi_j(g_l)\lambda_{jl}.$$
In all the above formulas $1\le i,j\le\tth$ and $\tth<k,l\le\theta.$ The cocycle $\tilde{\si}$ is now defined in exactly the same way as $\si,$ only in the last part $\tilde{\de}_j(b_l):=0$ as all the $\lambda_{jl}=0.$ This is because for $\tilde{\D},$ we wanted the component $I=\{1,\dots,\tth\}$ not to be linked to any other component.

The inverse $\si^{-1}$ is given in the same way by $\tau^{-1},$ where $$\tau^{-1}(u,b):=\varphi(\an u)(b)=\varphi(u)(\an^{-1}b).$$
If $\si,\tilde{\si}$ are two 2-cocycles for a Hopf algebra $\A,$ then $\rho:=\tilde{\si}\si^{-1}$ is again a 2-cocycle, but for the Hopf algebra $\A_{\si}.$ In our case this means $\rho$ is a 2-cocycle for the Hopf algebra $\A:=(\U\ot\B)_{\si}.$ Then $(\U\ot\B)_{\tilde{\si}}=\A_{\rho}$. If we can show that $\rho$ passes down naturally to a 2-cocycle $\rho'$ on $\A/\K^+$ such that $(\A/\K^+)_{\rho'}$ is isomorphic to $\A_{\rho}/\K^+,$ then the statement is clear:
\begin{equation} \uf(\tilde{\D})\simeq \A_{\rho}/\K^+\simeq (\A/\K^+)_{\rho'}\simeq \uf(\D)_{\rho'}. \end{equation}
Therefore, we want the following situation.
$$
\xy \xymatrix@C+1cm@R+1cm{\A\ot\A  \ar[r]^{\rho} \ar[d]^{\pi} & \k \\
     (\A\ot\A)/(\A\ot\K^++\K^+\ot\A) \simeq \A/\K^+\ot\A/\K^+ \ar@{{}-->}[ru]^{\rho'}\\} \endxy
$$
For this it suffices to show that $\rho$ is 0 on the kernel of the natural projection $\pi$ and we get the factorization and $(\A/\K^+)_{\rho'}\simeq \A_{\rho}/\K^+$ by definition of $\rho'.$
So we see that it is enough to show $\rho(\A,\K^+)=0=\rho(\K^+,\A).$ This means that for all $u\in\U,$ $b\in\B$ and $1\le i\le\tth$ we need
\begin{align*} \rho(u\ot b,z_i\ot g_i^{-1})&=\rho(u\ot b,1\ot 1),\\
                \rho(z_i\ot g_i^{-1},u\ot b)&=\rho(1\ot 1,u\ot b).
\end{align*}
We calculate using the definition of $\rho$, the convolution product and property 3 of $\tau$ and $\tilde{\tau}:$
\begin{align*}\rho(u\ot b,z_i\ot g_i^{-1})&=\tilde{\si}(u\_1\ot b\_1,z_i\ot g_i^{-1})\si^{-1}(u\_2\ot b\_2,z_i\ot g_i^{-1})\\
                &=\ep_{\U}(u\_1)\tilde{\tau}(z_i,b\_1)\ep_{\B}(g_i^{-1})\ep_{\U}(u\_2)\tau(\an z_i,b\_2)\ep_{\B}(g_i^{-1})\\
                &=\ep_{\U}(u)\tilde{\tau}(z_i,b\_1)\tau(z_i^{-1},b\_2),\\
        \rho(u\ot b,1\ot 1)&=\tilde{\si}(u\_1\ot b\_1,1\ot 1)\si^{-1}(u\_2\ot b\_2,1\ot 1)\\
                &=\ep_{\U}(u\_1)\tilde{\tau}(1,b\_1)\ep_{\B}(1)\ep_{\U}(u\_2)\tau(\an 1,b\_2)\ep_{\B}(1)\\
                &=\ep_{\U}(u)\ep_{\B}(b).
\end{align*}
As $z_i$ is group-like, we see by using property 2 of $\tau$ and $\tilde{\tau}$ that it is enough to verify 
\begin{equation} \label{chk1}\tilde{\tau}(z_i,b\_1)\tau(z_i^{-1},b\_2)=\ep_{\B}(b) \end{equation} on the generators of $\B.$
An analogue calculation, using property 1 and 4 this time, shows that for the second condition we have to verify for the generators of $\U$
\begin{equation} \label{chk2}\tilde{\tau}(u\_1,g_i^{-1})\tau(u\_2,g_i)=\ep_{\U}(u). \end{equation}
The verification goes as follows ($1\le j\le\tth, \tth<k\le\theta$):
\begin{align*}
\intertext{$b=y_k$ in (\ref{chk1})}
  \tilde{\tau}(z_i,y_k)\tau(z_i^{-1},y_k)&=\tilde{\ga}_i(y_k)\ga_i^{-1}(y_k)\\
                                &=\chi_i(y_k)\chi_i^{-1}(y_k)=1=\ep_{\B}(y_k),\\
\intertext{$b=b_k$ in (\ref{chk1})}
  \tilde{\tau}(z_i,b_k)\tau(z_i^{-1},1)+\tilde{\tau}(z_i,g_k)\tau(z_i^{-1},b_k)&=\tilde{\ga}_i(b_k)\ga_i^{-1}(1)+\tilde{\ga}_i(g_k)\ga_i^{-1}(b_k)\\
                                        &=0+0=\ep_{\B}(b_k),\\
\intertext{$u=z_j$ in (\ref{chk2})}
 \tilde{\tau}(z_j,g_i^{-1})\tau(z_j,g_i)&=\chi_j(g_i^{-1})\chi_j(g_i)=\chi_j(1)\\
                                &=1=\ep_{\U}(z_j),\\
\intertext{$u=u_j$ in (\ref{chk2})}
 \tilde{\tau}(u_j,g_i^{-1})\tau(1,g_i)+\tilde{\tau}(z_j,g_i^{-1})\tau(u_j,g_i)&=\tilde{\de}_j(g_i^{-1})\ep_{\B}(g_i)+\tilde{\ga}_j(g_i^{-1})\de_j(g_i)\\
                                &=0+0=\ep_{\U}(u_j).
\end{align*}
}
%%%%%%%%%%%%%%%%%%%%%%%%%%%%%%%%%%%%%%%%%%%%%%%%%%%%%
\section{A special root vector case}

Now we deal with Hopf algebras where the Dynkin diagram is just one copy of $A_n,$ $n\ge 1,$ and hence there can be no linking. To fix the order of the indices we require for the $n\times n$ Cartan matrix $a_{ij}=0$ whenever $|i-j|\ge 2.$
Given a linking datum $\D$ of $A_n$ type with group $\G$, we consider the algebra $\uf(\D,\uu)$ with certain central elements $u_{\alp}$ from the group algebra.

In \cite{AS-Poin}, the explicit form of the root vectors $x_{\alp}$ and all the possible families $u_{\alp}$ such that $\uf(\D,\uu)$ is a Hopf algebra have been determined. We repeat these here for later calculations.

Set $q_{i,j}:=\chi_j(g_i),\;1\le i,j\le n,$ $N:=\ord \chi_i(g_i)$ and $$B^{i,j}_{p,r}:=\prod_{i\le l<j, p\le h<r}q_{l,h}.$$ $N$ is defined independently of the choice of $i$ as $\chi_i(g_i)=\chi_j(g_j)$ due to (\ref{cartan}). Then $$C^{j}_{i,p}:=(1-q^{-1})^N(B^{p,j}_{i,p})^{\binom N2},\qquad 1\le i<p<j\le n+1.$$
The root vectors are defined inductively on the height of the corresponding root:
\begin{align}
e_{i,i+1}&:=x_i,&& 1\le i\le n;\\
\label{defroot}e_{i,j}&:=e_{i,j-1}e_{j-1,j}-B^{i,j-1}_{j-1,j}e_{j-1,j}e_{i,j-1} && j-i\ge 2.
\end{align}
We set $$ \chi_{i,j}:=\prod_{i\le l<j}\chi_l,\quad g_{i,j}:=\prod_{i\le l<j}g_l\quad\text{and}\quad h_{i,j}:=g_{i,j}^N.$$
Then $B^{i,j}_{p,r}=\chi_{p,r}(g_{i,j}).$
Next we introduce a parameter family $\ga$ containing for every root vector $e_{i,j}$ a scalar $\ga_{i,j}\in\k.$ The elements $u_{i,j}\in\k\G$ are also defined inductively:
\begin{equation}\label{defu} u_{i,j}(\ga):=\ga_{i,j}(1-h_{i,j})+\sum_{i<p<j}C^j_{i,p}\ga_{i,p}u_{p,j}(\ga),\quad 1\le i<j\le n+1. \end{equation}
%Here $h_{i,j}:=\prod_{i\le l<j}g_l^N$ are group elements. 
We require that for all $i<j,$
\begin{equation}\label{charn}\ga_{i,j}=0\quad\text{if}\quad\chi_{i,j}^N\neq\ep\;\text{or}\; h_{i,j}=1.\end{equation}
We will call such families \emph{admissible}. In \cite[Lemma 7.20.]{AS-Poin} it was shown that condition (\ref{charn}) implies the same statement for the $u_{i,j},$ and proves all the $u_{i,j}$ to be central in $H_0.$ 

For an admissible family $\ga$ we simply write $\A(\ga)$ instead of $\uf(\D,\uu=(u_{i,j}(\ga))).$ Theorem 7.25.(i) of \cite{AS-Poin} now states that $\A(\ga)$ is a finite dimensional pointed Hopf algebra. Our result is:
\beer{thm}{\label{main}For a linking datum of $A_n$ type and two admissible families $\ga$ and $\ga'$ of scalars the Hopf algebras $\A(\ga)$ and $\A(\ga')$ are quasi-isomorphic.
}
For the proof of this theorem we will first slightly generalize Theorem 7.24.~in \cite{AS-Poin}.

Let $\G$ be an abelian group and  $H_0$ a $\k$-algebra containing $\k\G$ as a subalgebra. Assume that there is an integer $P$, a subset $Z$ of $\{1,2,\dots,P\},$ integers $N_j>1$ where $j\in Z$ and elements $x_i\in H_0,$ $h_i\in\G,$ $\eta_i\in\Gh,$ $1\le i\le P,$ such that
\begin{align*}
g x_i&=\eta_i(g)x_i g\quad\text{for all }g\in\G,\; 1\le i\le P.\\
x_ix_j^{N_j}&=\eta_j^{N_j}(h_i)x_j^{N_j}x_i\quad\text{for all }1\le i\le P,\;j\in Z.\\
\text{The ele}&\text{ments }x_1^{a_1}\dots x_P^{a_P}g,\;a_1,\dots a_P\ge 0,g\in\G \text{ form a }\k\text{-basis of }H_0.
\end{align*}
\beer{thm}{\label{basisthm}
Let $u_j,\;j\in Z,$ be a family of elements in $\k\G,$ and $I$ the ideal in $H_0$ generated by all $x_j^{N_j}-u_j,$ $j\in Z.$ Let $A=H_0/I$ be the quotient algebra.

If $u_j$ is central in $H_0$ for all $j\in Z$ and $u_j=0$ if $\eta_j^{N_j}\neq\ep$,\\
then the residue classes of $x_1^{a_1}\dots x_P^{a_P}g,$ $a_i\ge 0,$ $a_j<N_j$ if $j\in Z,$
$g\in\G,$ form a $\k$-basis of $A.$
}
The proof of this theorem is exactly the same as in the original paper, where $Z$ included all indices from $1$ to $P.$\\

To see how this can be applied in our situation, we recall two more results (Theorem 7.21.~and Lemma 7.22.) from \cite{AS-Poin}.
First, the elements $e_{1,2}^{b_{12}}e_{1,3}^{b_{13}}\dots e_{n,n+1}^{b_{nn+1}}g,$ $g\in\G,$ $b_{ij}\ge 0,$ where the root vectors are arranged in the lexicographic order, form a $\k$-basis of $H_0.$ Furthermore, we have the following crucial commutation rule for all $1\le i<j\le n+1,$ $1\le s<t\le N+1,$
\begin{equation}\label{crucial}
e_{i,j}e_{s,t}^N=\chi_{s,t}^N(g_{i,j})e_{s,t}^Ne_{i,j}.
\end{equation}
We set now $H_0:=\Uf(\D)$ and know that the coradical is a subalgebra. So the last theorem gives us, for instance, a basis of $A(\ga),$ as we get $\uf$ from $\Uf$ exactly by dividing out the root vector relations.

We are now ready to give the proof of our Theorem.
\prf[of Theorem \ref{main}]{As in the linking case, it is enough to prove that for any admissible family $\ga$, $A(\ga)$ is quasi-isomorphic to $A(\ga_0)$ where $\ga_0$ denotes the family in which all the $\ga_{i,j}$ are zero. Again, this  will be achieved by following a stepwise procedure.

$\bu$ Fix $i_0$ with $1\le i_0\le n$ such that for the given admissible family $\ga$ we have $$\ga_{i,j}=0 \text{ for all }1\le i<i_0\text{ and }i<j\le n+1.$$ Set $\tilde{\ga}_{i,j}:=\ga_{i,j}$ if $i\neq i_0$ and $\tilde{\ga}_{i_0,j}:=0$ for all $i_0<j\le n+1.$ Then $\tilde{\ga}$ is again admissible. It is now sufficient to prove that $\A(\ga)$ and $\A(\tilde{\ga})$ are quasi-isomorphic and then to repeat this step with increasing $i_0,$ replacing $\ga$ by $\tilde{\ga}.$
Set
$$H=H_0/I,\; I:=(e_{i,j}^N-u_{i,j}\,:\,1\le i<j\le n+1,\;i\neq i_0\text{ or }\chi_{i,j}^N\neq\ep).$$
Note that $u_{i,j}(\ga)=u_{i,j}(\tilde{\ga})$ for all the $u_{i,j}$ appearing in the ideal $I.$

$\bu$ $H$ is a Hopf algebra.\\
To see this, we recall the comultiplication on the root vectors
\begin{equation}\begin{split}\label{delroot}&\del(e_{i,j}^N-u_{i,j})=(e_{i,j}^N-u_{i,j})\ot 1+h_{i,j}\ot(e_{i,j}^N-u_{i,j})+\\
&+\sum_{i<p<j}C_{i,p}^je_{i,p}^Nh_{p,j}\ot(e_{p,j}^N-u_{p,j})+\sum_{i<p<j}C_{i,p}^j(e_{i,p}^N-u_{i,p})h_{p,j}\ot u_{p,j}.
\raisetag{4em}
\end{split}\end{equation}
If $i<i_0$ we have $\ga_{i,j}=0$ for all $i<j\le n+1$ and hence $u_{i,j}=0.$ So all the summands in (\ref{delroot}) are in $H_0\ot I$ or $I\ot H_0,$ because $e_{i,p}^N\in I$ for all $i<p<j.$\\
The case $i>i_0$ is obvious. When $i=i_0$ then, according to the definition of $I$, $e_{i_0,j}^N-u_{i_0,j}$ is in $I$ only if $\chi_{i_0,j}^N\neq\ep.$ From (\ref{charn}) follows then $\ga_{i_0,j}=0$ and $u_{i_0,j}=0.$ We have 
$$\chi_{i_0,j}^N=\chi_{i_0,p}^N\chi_{p,j}^N\qquad\text{for all }i_0<p<j.$$
And hence $\chi_{i_0,p}^N\neq\ep$ and $u_{i_0,p}=0$ or $\chi_{p,j}^N\neq\ep$ and $u_{p,j}=0.$ This proves that 
$$\del(I)\subset H_0\ot I+I\ot H_0.$$
A simple calculation shows $\ep(e_{i,j})=0=\ep(u_{i,j}).$ So we can get two recursion formulas for the antipode of $e_{i,j}^N-u_{i,j}$ from the comultiplication formula, as $\an(a\_1)a\_2=\ep(a)=a\_1\an(a\_2).$ Using the second of these formulas for the case $i<i_0$ we immediately have $\an(e_{i,j}^N)\in I.$ When $i>i_0,$ an inductive argument using the first formula and $$\an(e_{i,i+1}^N-u_{i,i+1})=-g_i^{-N}(e_{i,i+1}^N-u_{i,i+1})$$ gives again $\an(e_{i,j}^N)\in I.$ For $i=i_0,$ a combination of the reasoning from the discussion of the comultiplication and the inductive argument from the last case give the desired result, establishing $\an(I)\subset I.$ Hence $I$ is a Hopf ideal.

$\bu$ Next we take as $K$ the subalgebra of $H$ generated by the group $\G$ and the remaining $e_{i_0,j}^N,$ i.e.~$i_0<j\le n+1,$ such that $\chi_{i_0,j}^N=\ep.$ A similar calculation to the one above reveals that $K$ is actually a Hopf subalgebra of $H.$ This time, one has to use the fact that for any $p$ between $i_0$ and $j,$ either $\chi_{i_0,p}^N=\ep$ or that $\chi_{p,j}^N\neq\ep$ and hence $u_{p,j}=0=e_{p,j}^N$ in $H.$

As all the $u_{i,j}$ fulfill the conditions of Theorem \ref{basisthm} we immediately get a basis of $H.$ The commutation relations (\ref{crucial}) for the $\pth N$ powers of the root vectors show that the generators of $K$ all commute with each other, because the factor is $1,$ as $\chi_{i_0,j}^N=\ep$ for all generators. Hence, any monomial in $K$ can be reordered and is then a basis element of $H.$ So $K$ is just the polynomial algebra on its generators.

$\bu$ We define an algebra map $f:K\rightarrow \k$ by setting 
$$f(e_{i_0,j}^N):=\ga_{i_0,j},\qquad f(g):=1\qquad\text{on all the generators of }K,\,g\in\G.$$
Algebra maps, from a Hopf algebra $K$ to the base field, form a group under the convolution product where the inverse is given by the composition with the antipode. This group acts on the Hopf algebra $K$ from the left and the right by
$$f.x=x\_1f(x\_2),\qquad x.f=f(x\_1)x\_2.$$
$\bu$ To be able to apply Theorem 2.~of \cite{Mas}, which will give us the desired quasi-isomorphism, we have to calculate $f.e_{i_0,j}^N.f^{-1}.$ As a preparation for this we first calculate $f(u_{i,i+1})=0$ for all $1\le i\le n$ and see then from the inductive definition (\ref{defu}) that $f(u_{i,j})=0$ for all $i<j.$ So for the generators of $K$ we have
\begin{equation*}\begin{split} f.e_{i_0,j}^N&=e_{i_0,j}^Nf(1)+h_{i_0,j}f(e_{i_0,j}^N)+\sum_{i_0<p<j}C_{i_0,p}^je_{i_0,p}^Nh_{p,j}f(e_{p,j}^N)\\&=e_{i_0,j}^N+\ga_{i_0,j}h_{i_0,j},
\end{split}\end{equation*}
as $e_{p,j}^N=u_{p,j}$ in $H.$ The recursive formula for the antipode of the $\pth N$ powers of the root vectors is 
\begin{equation*} \an(e_{i_0,j}^N)=-h_{i_0,j}^{-1}e_{i_0,j}^N-\sum_{i_0<p<j}C_{i_0,p}^jh_{p,j}^{-1}\an(e_{i_0,p}^N)e_{p,j}^N.
\end{equation*}
So we get
\begin{equation}\begin{split}\label{fif} (f.&e_{i_0,j}^N).f^{-1}=\ga_{i_0,j}f(h_{i_0,j}^{-1})h_{i_0,j}+f(\an(e_{i_0,j}^N))+f(h_{i_0,j}^{-1})e_{i_0,j}^N+\\ 
&\qquad\qquad\qquad+\sum_{i<p<j}C_{i,p}^jf(h_{i_0,j}^{-1})f(\an(e_{i_0,p}^N))e_{p,j}^N\\
&=e_{i_0,j}^N+\ga_{i_0,j}h_{i_0,j}+\\
&\quad+f\left(-h_{i_0,j}^{-1}e_{i_0,j}^N-\!\sum_{i_0<p<j}C_{i_0,p}^jh_{p,j}^{-1}\an(e_{i_0,p}^N)e_{p,j}^N\right)+\\
&\quad+\sum_{i<p<j}C_{i,p}^jf\left(-h_{i_0,p}^{-1}e_{i_0,p}^N-\!\sum_{i_0<q<p}C_{i_0,q}^ph_{q,p}^{-1}\an(e_{i_0,q}^N)e_{q,p}^N\right)u_{p,j}\\
&=e_{i_0,j}^N+\ga_{i_0,j}h_{i_0,j}-f(e_{i_0,j}^N)-\!\sum_{i<p<j}C_{i,p}^jf(e_{i_0,p}^N)u_{p,j}\\
&=e_{i_0,j}^N-\ga_{i_0,j}(1-h_{i_0,j})-\!\sum_{i<p<j}C_{i,p}^j\ga_{i_0,p}u_{p,j}\\
&=e_{i_0,j}^N-u_{i_0,j}\qquad\text{by (\ref{defu})}.
\end{split}\end{equation}
Two parts of the sum vanish as they are multiples of $f(e_{p,j}^N=u_{p,j})=0=f(e_{q,p}^N=u_{q,p}).$ Note that in the second to last step, if $j>i_0$ is such that $\chi_{i_0,j}^N\neq\ep,$ then $e_{i_0,j}^N$ is not a generator of $K$ and we can not simply apply the definition of $f.$ But in this case, $u_{i_0,j}=0=e_{i_0,j}^N$ in $H$ and $\ga_{i_0,j}$ is zero as well, as this is required for an admissible family. So we still have $f(e_{i_0,j}^N)=\ga_{i_0,j}.$

$\bu$ Let $J$ be the Hopf ideal of $K$ generated by all the generators $e_{i_0,j}^N$ of $K.$ Then,
according to \cite[Theorem 2.]{Mas}, $H/(f.J)$ is an $(H/(J),H/(f.J.f^{-1}))$-bi-Galois object and hence $H/(J)$ and $H/(f.J.f^{-1})$ are quasi-isomorphic if the bi-Galois object is not zero.
We see that $u_{i_0,j}(\tilde{\ga})=0$ and so $A(\tilde{\ga})=H/(J).$ Calculation (\ref{fif}) showed that $A(\ga)=H/(f.J.f^{-1}).$ We are left to show that $B:=H/(f.J)$ is not zero.

$\bu$ $B=H_0/(I,f.J)$ by construction. We have a basis of $H_0$ and see that we could apply Theorem \ref{basisthm} to get a basis of $B$. It just remains to check that the elements $\ga_{i_0,j}h_{i_0,j}$ appearing in $f.J$ satisfy the conditions of the theorem.

If $\chi_{i_0,j}^N\neq\ep$ then $\ga_{i_0,j}=0$, because $\ga$ is admissible. $h_{i_0,j}$ is in $\G$ and so commutes with all group elements. We will show now that $h_{i_0,j}$ commutes also with all the generators of $H_0.$ For this we calculate
\begin{equation*}\begin{split}h_{i_0,j}x_k&=\chi_k(h_{i_0,j})x_kh_{i_0,j},\\
\chi_k(h_{i_0,j})&=\chi_k(\prod_{i_0\le p<j}g_p^N)
=\prod_{i_0\le p<j}\chi_k^N(g_p)\\&=\prod_{i_0\le p<j}\chi_p^N(g_k^{-1})
=\chi_{i_0,j}^N(g_k^{-1})=1.
\end{split}\end{equation*}
Here we used (\ref{cartan}) and the fact that $N$ is exactly the order of any diagonal element $\chi_p(g_p).$ So $B$ is not zero and the statement is proven.}

%%%%%%%%%%%%%%%%%%%%%%%%%%%%%%%%%%%%%%%%%%%%%%%%%%%%%%%%%%%%%%%%%%%%%%%%%%%%%%%%
\section{The mixed case}

Now finally, we consider the case where linking parameters \emph{and} root vector parameters appear. As we do not know yet how to generalize our considerations to arbitrary Dynkin diagrams of finite type, we will consider only copies of $A_n.$

So we take a linking datum $\D$ of type $A_{n_1}\xdig\dots\xdig A_{n_t},$ $n_k,t>0,$ with a fixed finite abelian group $\G.$
The algebra $\Uf(\D)$ 
%is defined in the same way as $\uf(\D)$ in the introduction except that we leave out relation (\ref{rootrel}). This 
is a Hopf algebra.%, cf. \cite{D,New-AS}.
We order the vertices in the Dynkin diagram so that the root vectors of the $\pth k$ component are $e_{S_k+i,S_k+j},$ $1\le i<j\le n_k+1,$ 
where $S_k=n_1+\dots+n_{k-1}.$ The root vectors within one component are defined in the same way as in the previous section. Then the monomials
$$ e_{1,2}^{b_{1,2}}\cdots e_{1,n_1+1}^{b_{1,n_1+1}}e_{2,3}^{b_{2,3}}\cdots e_{S_t+n_t,S_t+n_t+1}^{b_{S_t+n_t,S_t+n_t+1}}g,\quad 0\le b_{S_k+i,S_k+j}, g\in\G,$$
form a PBW-basis of $\Uf(\D).$  A proof can be found in \cite[Theorem 4.2.]{AS-char}. The idea is that a PBW-basis is known for the linking datum $\D_0$ where all the $\lam$ are zero. Using, for instance, the considerations in the first section of this chapter one knows explicitly the 2-cocycle relating $\Uf(\D)$ and $\Uf(\D_0).$ Expressing now the above monomials in $\Uf(\D_0),$ one sees that they are basis elements of the common underlying vector space.

As in the previous section, we introduce for every component $k$ of the diagram an admissible parameter family $\ga_k$ and define the corresponding elements $u_{S_k+i,S_k+j}(\ga_k)\in\k\G,$ $1\le i<j\le n_k+1.$ The collection of all the parameters $\ga_k$ will be denoted by $\ga.$

$\A(\D,\ga)$ is now defined as the quotient of $\Uf(\D)$ by the ideal generated by the root vector relations
$$ e_{S_k+i,S_k+j}^{N_k}=u_{S_k+i,S_k+j}(\ga_k),\qquad 1\le k\le t,\,1\le i<j\le n_k+1. $$
Here $N_k$ is again the common order of the diagonal elements $\chi_i(g_i)$ with $i$ in the $\pth k$ component of the diagram. Hence $\A(\D,\ga)=\uf(\D,(u_{i,j}(\ga))).$
% and because there are no self-linkings, $\a=0$ when $i$ is linked to $j$ and (\ref{serre}) simplifies into
%\begin{equation}\label{links}x_ix_j-\chi_j(g_i)x_jx_i
\beer{thm}{\label{dimthm}
The so defined algebra $\A(\D,\ga)$ is a Hopf algebra of dimension $N_1^{\binom{n_1+1}2}\cdots N_t^{\binom{n_t+1}2},$ whose basis consists of the monomials
\begin{equation}\label{basismonoms}
e_{1,2}^{b_{1,2}}\cdots e_{1,n_1+1}^{b_{1,n_1+1}}e_{2,3}^{b_{2,3}}\cdots e_{S_t+n_t,S_t+n_t+1}^{b_{S_t+n_t,S_t+n_t+1}}g,\quad 0\le b_{S_k+i,S_k+j}<N_k, g\in\G.
\end{equation}
}
\prf{The statement about the Hopf algebra is clear because we already know from earlier considerations in every component, that the ideal is a Hopf ideal. For the dimension and basis we will use Theorem \ref{basisthm}. We have to check all the conditions.

We first need to prove the commutation relation between root vectors of one component and $\pth N$ powers of root vectors of the other components. Secondly, we have to show that the $u_{i,j}$ are central with regard to all generators $x_a$ of $\Uf(\D)$:
\begin{align}
\label{commute}e_{S_k+i,S_k+j}^{\phantom{N_l}}e_{S_l+r,S_l+s}^{N_l}&=\chi_{S_l+r,S_l+s}^{N_l}(g_{S_k+i,S_k+j})e_{S_l+r,S_l+s}^{N_l}e_{S_k+i,S_k+j}^{\phantom{N_l}},\\
\label{central}u_{S_k+i,S_k+j}x_a&=x_au_{S_k+i,S_k+j},
\end{align}
where $1\le k,l\le t,\,1\le i<j\le n_k+1,\,1\le r<s\le n_l+1,\,1\le a\le S_t+n_t.$

When we are within one component, i.e. $k=l$ or $S_k<a\le S_{k+1},$ the above equations follow already from the original paper \cite{AS-Poin}.\\
The root vectors are linear combinations of homogeneous monomials in the generators $x_a$ of $\Uf(\D)$. Hence we see that Lemma \ref{techcom} will establish (\ref{commute}) for $k\neq l$.
 
(\ref{central}) is shown by induction on $j-i.$ Because of the recursive definition (\ref{defu}) of the $u_{i,j}$, the crucial part is:
\begin{align*}\ga_{k_{i,j}}(1-h_{S_k+i,S_k+j})x_a&=\ga_{k_{i,j}}x_a(1-\chi_a(h_{S_k+i,S_k+j})h_{S_k+i,S_k+j})\\&=\ga_{k_{i,j}}x_a(1-\chi_{S_k+i,S_k+j}^{N_k}(g_a^{-1})h_{S_k+i,S_k+j})\\&=x_a\ga_{k_{i,j}}(1-h_{S_k+i,S_k+j}).
\end{align*}
In the second step we used (\ref{cartan}) and as $a$ is not in the $\pth k$ component, all the corresponding entries of the Cartan matrix are zero. The third step uses the premise that $\ga$ is admissible.

Now we can apply Theorem \ref{basisthm} and the proof is finished.
}
\beer{lemma}{\label{techcom} For all indices $i$ not in the $\pth k$ component of the diagram and $S_k<j\le l\le S_k+n_k$ we have with $q=\chi_j(g_j)=\chi_l(g_l)$\\

$i)\quad x_ie_{j,l+1}=$
\begin{subnumcases}{\hspace{-3em}=}
\chi_j(g_i)x_jx_i+\lam(1-g_ig_j),&\hspace{-5em}if $j=l,$\label{case1}\\
\chi_{j,l+1}(g_i)e_{j,l+1}x_i+\lam(1-q^{-1})e_{j+1,l+1},&\hspace{-3em}if $\lambda_{il}=0,\;j<l,$\label{case2}\\
\chi_{j,l+1}(g_i)e_{j,l+1}x_i-\lambda_{il}(1-q^{-1})\chi_{j,l}(g_i)e_{j,l}g_ig_l,&\hspace{-1em}otherwise.\label{case3}
\end{subnumcases}
$ii)\quad x_ie_{j,l+1}^{N_k}=\chi_{j,l+1}^{N_k}(g_i)e_{j,l+1}^{N_k}x_i.$
}

\prf{
$i)$ $\bu$ The case $j=l$ is simply the defining relation (\ref{links}).

From now on $j<l$. We consider all possible linkings.

$\bu$ If $i$ is not linked to any vertex $p$ with $j\le p\le l$, then $\lam=\lambda_{il}=0$ and a repeated use of (\ref{links}) gives (\ref{case2}).

$\bu$ If $i$ is linked to $j$, then it can not be linked to $l$ as well. Hence $\lambda_{il}=0$ and we have to show the second case. We proceed by induction on $l-j$ and use the recursive definition (\ref{defroot}) of the root vectors.

For $l-j=1$ we have
\begin{align*} x_ie_{j,l+1}&=x_i[x_jx_l-\chi_l(g_j)x_lx_j]\\&=\chi_j(g_i)x_jx_ix_l+\lam(1-g_ig_j)x_l-\chi_l(g_j)\chi_l(g_i)x_lx_ix_j\\&=\chi_j(g_i)\chi_l(g_i)x_jx_lx_i+\lam x_l-\lam\chi_l(g_ig_j)x_lg_ig_j\\&\quad-\chi_l(g_j)\chi_l(g_i)\chi_j(g_i)x_lx_jx_i-\chi_l(g_j)\chi_l(g_i)x_l\lam(1-g_ig_j)\\&=\chi_{j,l+1}(g_i)[x_jx_l-\chi_l(g_j)x_lx_j]x_i+\lam(1-\chi_l(g_i)\chi_l(g_j))x_l\\&=\chi_{j,l+1}(g_i)e_{j,l+1}x_i+\lam(1-\underbrace{\chi_i(g_l^{-1})\chi_j(g_l^{-1})}_{=1}\underbrace{\chi_j^{x_{jl}}(g_j)}_{=q^{-1}})x_l.
\end{align*}
We used (\ref{cartan}) and the condition $\chi_i\chi_j=1$ as $\lam\neq 0.$ For the induction step we use an analogue calculation. The last steps are as follows
\begin{align*} x_ie_{j,l+1}&=\chi_{j,l+1}(g_i)e_{j,l+1}x_i+\lam(1-q^{-1})[e_{j+1,l}x_l-\chi_l(g_{j,l})\chi_l(g_i)x_le_{j+1,l}]\\
&=\chi_{j,l+1}(g_i)e_{j,l+1}x_i+\\&\qquad+\lam(1-q^{-1})[e_{j+1,l}x_l-\chi_l(g_{j+1,l})\underbrace{\chi_l(g_j)\chi_l(g_i)}_{\chi_j(g_l^{-1})\chi_i(g_l^{-1})}x_le_{j+1,l}]\\&=\chi_{j,l+1}(g_i)e_{j,l+1}x_i+\lam(1-q^{-1})e_{j+1,l+1}.
\end{align*}
$\bu$ If $i$ is linked to $l$ we have $\lambda_{il}\neq 0$ and hence we need to prove (\ref{case3}). A direct calculation using the definition of the root vectors and (\ref{case2}) gives
\begin{align*}
x_ie_{j,l+1}&=x_i[e_{j,l}x_l-\chi_l(g_{j,l})x_le_{j,l}]\\
&=\chi_{j,l}(g_i)e_{j,l}x_ix_l-\chi_l(g_{j,l})\chi_l(g_i)x_lx_ie_{j,l}-\chi_l(g_{j,l})\lambda_{il}(1-g_ig_l)e_{j,l}\\
&=\chi_{j,l}(g_i)e_{j,l}\chi_l(g_i)x_lx_i+\chi_{j,l}(g_i)e_{j,l}\lambda_{il}(1-g_ig_l)\\
&\quad -\chi_l(g_{j,l})\chi_l(g_i)x_l\chi_{j,l}(g_i)e_{j,l}x_i-\chi_l(g_{j,l})\lambda_{il}e_{j,l}\\
&\quad +\chi_l(g_{j,l})\lambda_{il}\chi_{j,l}(g_ig_l)e_{j,l}g_ig_l\\
&=\chi_{j,l+1}(g_i)[e_{j,l}x_l-\chi_l(g_{j,l})x_le_{j,l}]x_i+[\chi_{j,l}(g_i)-\chi_l(g_{j,l})]\lambda_{il}e_{j,l}\\
&\quad +[\chi_l(g_{j,l})\chi_{j,l}(g_ig_l)-\chi_{j,l}(g_i)]\lambda_{il}e_{j,l}g_ig_l.
\end{align*}
As $i$ is not in the $\pth k$ component we have by (\ref{cartan}) and $\chi_i\chi_l=1$
$$\chi_{j,l}(g_i)=\chi_i^{-1}(g_{j,l})=\chi_l(g_{j,l}).$$
Hence the second term in the last step of the above calculation vanishes, and for the bracket of the third term we calculate
\begin{align*}
\chi_l(g_{j,l})\chi_{j,l}(g_ig_l)-\chi_{j,l}(g_i)&=\chi_{j,l}(g_i)(\chi_l(g_{j,l})\chi_{j,l}(g_l)-1)\\
\chi_l(g_{j,l})\chi_{j,l}(g_l)&=\chi_l(g_j)\chi_l(g_{j+1})\cdots\chi_l(g_{l-1})\\
&\quad\cdot\chi_j(g_l)\chi_{j+1}(g_l)\cdots\chi_{l-1}(g_l)\\
&=1\cdot1\cdots1\cdot\chi_l(g_l)^{-1}.
\end{align*}
$\bu$ For the last case where $i$ is linked to a vertex $p$ with $j<p<l,$ we again proceed by induction on $l-p.$
As $\lam=\lambda_{jl}=0$ we have to show (\ref{case2}).

If $l-p=1$ we use the recursive definition of the root vectors and then (\ref{case3}). We set $F:=\lambda_{i(l-1)}(1-q^{-1})\chi_{j,l-1}(g_i)$ and have
\begin{align*}
x_ie_{j,l+1}&=x_i[e_{j,l}x_l-\chi_l(g_{i,l})x_le_{j,l}]\\
&=[\chi_{j,l}(g_i)e_{j,l}x_i+F\cdot e_{j,l-1}g_ig_{l-1}]x_l\\
&\quad -\chi_l(g_{i,l})\chi_l(g_i)x_l[\chi_{j,l}(g_i)e_{j,l}x_i+F\cdot e_{j,l-1}g_ig_{l-1}]\\
&=\chi_{j,l}(g_i)e_{j,l}\chi_l(g_i)x_lx_i-\chi_l(g_{i,l})\chi_l(g_i)\chi_{j,l}(g_i)x_le_{j,l}x_i\\
&\quad +F\cdot\chi_l(g_ig_{l-1})e_{j,l-1}x_lg_ig_{l-1}-F\cdot\chi_l(g_{i,l})\chi_l(g_i)x_le_{j,l-1}g_ig_{l-1}\\
&=\chi_{j,l+1}(g_i)e_{j,l+1}x_i\\&\quad+F\cdot\chi_l(g_ig_{l-1})[e_{j,l-1}x_l-\chi_l(g_{i,l-1})x_le_{j,l-1}]g_ig_{l-1}.
\end{align*}
However, the last square bracket is zero according to \cite[(7.17)]{AS-Poin}. The induction step is now simple.\\

$ii)$ $\bu$ The case $\lam=\lambda_{il}=0$ is trivial.

$\bu$ So let now $j=l$ and $\lam\neq 0$. Then $\chi_j(g_i)=\chi_i^{-1}(g_j)=\chi_j(g_j)=q$ and using (\ref{case1}) we have
\begin{align*}x_ix_j^{N_k}&=q^{N_k}x_j^{N_k}x_i+\lam(1+q+q^2+\dots+q^{N_k-1})x_j^{N_k-1}(1-q^{N_k-1}g_ig_j)\\
&=x_j^{N_k}x_i.\end{align*}
Here we used the fact that $N_k$ is the order of $q$.\\
From now on again $j<l.$

$\bu$ If $\lam\neq 0,$ we set $x=e_{j,l+1},$ $y=x_i,$ $z=\lam(1-q^{-1})e_{j+1,l+1},$ $\alp=\chi_{j,l+1}(g_i)$ and $\beta=\chi_{j+1,l+1}^{-1}(g_{j,l+1}).$ Then, because of \cite[(7.24)]{AS-Poin} $zx=\beta xz.$ Hence, cf. \cite[Lemma 3.4 (1)]{AS-A2},
$$yx^{N_k}=\alp^{N_k}x^{N_k}y+\left(\sum_{m=0}^{N_k-1}\alp^m\beta^{N_k-1-m}\right)x^{N_k-1}z.$$
Using $\chi_i\chi_j=1$ and (\ref{cartan}) we see that $\alp=\chi_i^{-1}(g_{j,l+1})=\chi_j(g_{j,l+1})$ and so
\begin{align*}
\alp^m\beta^{N_k-1-m}&=\beta^{N_k-1}\chi_j^m(g_{j,l+1})\chi_{j+1,l+1}^m(g_{j,l+1})\\
&=\beta^{N_k-1}\chi_{j,l+1}^m(g_{j,l+1})=\beta^{N_k-1}(B_{j,l+1}^{j,l+1})^m=\beta^{N_k-1}q^m.
\end{align*}
The last equality follows from \cite[(7.5)]{AS-Poin}. The geometric sum gives zero again.

$\bu$ The final case $\lambda_{il}\neq 0$ is treated similarly to the previous one. This time $z=-\lambda_{il}(1-q^{-1})\chi_{j,l}(g_i)e_{j,l}g_ig_l$ and $\beta=\chi_{j,l+1}(g_ig_l)\chi_{j,l+1}(g_{j,l}),$ because of \cite[(7.23)]{AS-Poin}. So we have
\begin{align*}
\alp^m\beta^{N_k-1-m}&=\beta^{N_k-1}\chi_{j,l+1}^m(g_i)\chi_{j,l+1}^{-m}(g_i)\chi_{j,l+1}^{-m}(g_{j,l+1})\\
&=\beta^{N_k-1}(B_{j,l+1}^{j,l+1})^{-m}=\beta^{N_k-1}q^{-m}.
\end{align*}
}
Here is the final result.
\beer{thm}{\label{mixedthm}
Let $\D$ and $\D'$ be two linking data as above and $\ga$ and $\ga'$ two admissible parameter families. Then $\A(\D,\ga)$ and $\A(\D',\ga')$ are quasi-isomorphic.
}
\prf{The proof is just a combination of the results obtained in the previous sections. First we show that the Hopf algebras are quasi-isomorphic to ones where all the parameters $\ga$ are zero. Hence, $\A(\D,\ga)$ is a cocycle deformation of $\uf(\D)$ and $\A(\D',\ga')$ is a cocycle deformation of $\uf(\D')$. Then we can use Theorem \ref{linkthm} to see that $\uf(\D)$ is quasi-isomorphic to $\uf(\D')$. Because of the transitivity of the quasi-isomorphism relation, this is the desired result.%This is done component-wise and each step is identical to the proof of Theorem \ref{main}. The necessary commutation relations are in Lemma \ref{techcom}.
%This way we find that the Hopf algebras $\A(\D,\ga)$ and $\A(\D',\ga')$ are quasi-isomorphic to $\uf(\D)$ and $\uf(\D'),$ respectively. Now Theorem \ref{linkthm} completes the proof.\\

To show that the $\ga$ can all be set zero, we proceed again stepwise. Let $1\le i_0\le S_t+n_t$ be such that $\ga_{i,j}=0$ for all $i<i_0$ and all $j>i$ for which there is a root vector $e_{i,j}$. We set $\tilde{\ga}_{i,j}:=\ga_{i,j}$ when $i\neq i_0$ and zero otherwise. It is enough to prove that $\A(\D,\ga)$ is quasi-isomorphic to $\A(\D,\tilde{\ga})$. Repeating the last step with increased $i_0$ and $\tilde{\ga}$ as the new $\ga,$ we find that $\A(\D,\ga)$ is quasi-isomorphic to $\A(\D,0)=\uf(\D)$ for all admissible $\ga$.

We set $H_0:=\Uf(\D)$ and consider for $1\le k\le t$ its ideals 
$$I_k:=(e_{i,j}^{N_k}-u_{i,j}\,:\,S_k< i<j\le S_{k+1}+1,\;i\neq i_0\text{ or }\chi_{i,j}^{N_k}\neq\ep).$$
When $i_0$ is not in the $\pth k$ component, i.e. $i_0\le S_k$ or $i_0>S_{k+1}$, we know that $I_k$ is a Hopf ideal from the proof of Theorem 7.25.(i) in \cite{AS-Poin}. For $S_{k_0}<i_0\le S_{k_0+1}$ the considerations that prove $I_{k_0}$ to be a Hopf ideal have been carried out explicitly in the proof of our Theorem \ref{main}. Hence $H:=H_0/I$ is a Hopf algebra, where $I$ denotes the sum of all the ideals $I_k.$ As before, define $K$ as the Hopf subalgebra of $H$ which is generated by $\G$ and the remaining $e_{i_0,j}^{N_{k_0}},$ $S_{k_0}<i_0<j\le S_{k_0+1}+1,$ with $\chi_{i,j}^{N_{k_0}}=\ep.$ Using (\ref{crucial}) we see that $K$ is commutative. Because Lemma \ref{techcom} establishes (\ref{commute}), we can apply Theorem \ref{basisthm} to find a basis of $H$ and see that $K$ is just the polynomial algebra on its generators. Hence, the algebra map $f:K\rightarrow \k$ is well defined by setting
$f(e_{i_0,j}^{N_{k_0}}):=\ga_{i_0,j}$ on all the generators of K and $f(g):=1$ for all $g\in\G.$ The analogue of computation (\ref{fif}) establishes $f.e_{i_0,j}^{N_{k_0}}.f^{-1}=e_{i_0,j}^{N_{k_0}}-u_{i_0,j}.$ We define $J$ as the Hopf ideal of $K$ generated by all the $e_{i_0,j}^{N_{k_0}}$ in $K$. Now we can apply \cite[Theorem 2]{Mas} to prove that $\A(\D,\ga)=H/(f.J.f^{-1})$ and $\A(\D,\tilde{\ga})=H/(J)$ are quasi-isomorphic if $B:=H/(f.J)\neq 0.$ For this last step we calculate $f.e_{i_0,j}^{N_{k_0}}=e_{i_0,j}^{N_{k_0}}+\ga_{i_0,j}h_{i_0,j}.$ As $\ga$ is assumed to be admissible, we see from Lemma \ref{techcom} that we can use Theorem \ref{basisthm} again to find a basis of $B=H_0/(I,f.J)$ consisting of the monomials (\ref{basismonoms}). So $B$ is not zero and everything is proved.
}

The theorem extends the original results of \cite{Mas}, which dealt with the case of copies of $A_1$ only, and \cite{BDR-B2} that includes a proof for the diagram $A_2$. For this proof the authors need to express the Hopf algebra, however, as an Ore extension.

It should be possible to extend our proof to arbitrary Dynkin diagrams of finite Cartan type. The only crucial part is the commutation relations between the root vectors and their powers. They should be checked for the other diagrams. An analogue of Lemma \ref{techcom} for the commutation relations between root vectors belonging to different Cartan types would be especially useful. 

\rem{The Hopf algebras in Section \ref{rank2comps} coming from self-linkings provide a nice class of exceptional examples, where most of the considerations of this chapter are not applicable. The stepwise approach used in the above theorems to prove quasi-isomorphism is not possible there, because the linking parameters appear in the root vector parameters. We can not even apply Theorem \ref{basisthm}, simply because most of the root vector parameters are not in the group algebra and are not central. But having got all relations from felix, we know that the diamond lemma gives us a basis already. 

In a private note \cite{Mas-priv}, however, Masuoka showed that all the Hopf algebras arising from self-linking $A_2$ \emph{are} quasi-isomorphic. For terminology and conventions we refer to Section \ref{selfA2sec}. His approach is as follows. In a first step he deals with the linking relations and the root vector relations for the simple roots simultaneously. For this he 
sets $I:=(x_i^3, zx_i-q^ix_iz; i=1,2)$ and $H$ is $\Uf(A_2)$, but \emph{without} the Serre relations (\ref{serre}). Then he shows using his \cite[Theorem 2]{Mas} that $T:=H/(I.f^{-1})$ is a bi-Galois object for $L':=H/(f.I.f^{-1})$ and $L:=H/(I)$. Here $f$ is the algebra map sending the generators of $I$ to $-\mu_i$ and $-\ga_i$ respectively. Defining
\begin{align*}
u&:=z^3-(q-1)^3\mu_1\mu_2-(1-q)\ga_1\ga_2+\l\qquad\text{in }T,\\
v&:=z^3+(q-1)^3\mu_1\mu_2(g_2^3-1)+\\
&\quad +(1-q)\ga_1\ga_2(g_1g_2^2-1)-\l(g_1^3g_2^3-1)\qquad\text{in }L',
\end{align*}
gives that $\uf(A_2)=L'/(v)$ and $\uf(\D_0,\bm 0)=L/(z^3)$ with $\D_0$ the linking datum with no linking. With ingenious insight and explicit calculations he then shows that $Q:=T/(u)$ is actually a bi-Galois object for these two Hopf algebras and hence they are quasi-isomorphic. The calculations in this case are slightly trickier than for \cite[Proposition 3.3]{BDR-B2}, as the structure maps for $T$ map $T$ into $L'\ot T$ and $T\ot L,$ respectively. When calculating the image of $u$ under these maps, one therefore has to apply different commutation rules, depending on which tensor factor one is in. All technical difficulties can be dealt with directly, as the commutation relations are explicit and the diamond lemma can be applied.

The proof for $A_2$ should be easily transferable to the self-linking of $B_2$, as given in Figure \ref{B2rel}. However, the calculations are much more involved. Even if we could guess the right ``$u$'', giving us a bi-Galois object for the algebras that incorporate the root vector relation for $z$, proving it directly seems hopeless. The expressions for $\del(z^5)$ are just too messy. And a computer algebra program like felix can not be applied directly, as the separate tensor factors have different commutation rules. Besides, after dealing with the relations for $z^5$, we would then still have to incorporate the relations for $u^5$. So, again we need a deeper insight into the general structure theory for these self-linked algebras first, before we attempt a quasi-isomorphism theorem.
}

%defines an algebra map and constructs a bi-Galois object $T$ for two algebras $L', L$ in a similiar way as was done above. Through ingenious insight Masuoka then incorporates the root vector relation for $z$ into these objects.
%%% Local Variables: 
%%% mode: latex
%%% TeX-master: "main"
%%% End: 

\appendix
  \renewcommand{\thesection}{A.\arabic{section}}
  \setcounter{section}{0}   
\chapter{Felix programs}

Here we give the programs used for the calculations in Section \ref{rank2comps}. To run these scripts one has to have a working copy of {\bf felix} \cite{felix}, which can be downloaded at {\tt http://felix.hgb-leipzig.de/}. For the results concerning the coproduct of higher powers of the root vectors one needs a tensor module. Istv\'an Heckenberger (heckenbe@mathematik.uni-leipzig.de) was kind enough to provide his extension module \verb1tensor.cmp1, of which we use the \verb1ttimes1 function.
Saving one of the scripts below in a {\it file}, one can execute it by calling \verb1 felix < 1{\it file} in a Unix-like system environment. The tensor module is compiled in the same way.

We want to give some comments on these scripts.
\begin{itemize}
\item If the tensor module is not available, then only the Groebner basis calculation is possible. In this case remove the first line from the scripts and the first \%-sign in front of the first \verb1bye1.
\item Each script starts by defining a domain with its parameters and the variables. The symbol \verb1ixi1 is defined as a variable and treated by the tensor module as the tensor sign $\ot$.
\item Then a matrix is given which gives an appropriate term-ordering.
\item The ideal has all the defining relations of the algebra. Because of the special property $q^p=1,$ we need to treat $q$ formally as a variable and not as a parameter. By assigning a zero to its position in the ordering matrix and giving all the commutation relations for $q$, felix treats $q$ in effect like a parameter.
\item The function \verb1standard1 computes the Groebner basis. As the algorithm for the non-commutative case is non-deterministic, the right term-ordering can be essential. This is especially important for the $G_2$ case. We want to thank Istv\'an Heckenberger for showing us how to determine a useful ordering.
\item The coproduct for the variables is denoted by \verb1del1.
\item After calculating the higher powers of the coproduct of the root vectors, the counter terms are guessed from the output and the result is already incorporated in the scripts. The summands are all on separate lines and the coefficients are also separately defined.
\item The last part of the calculations is always a test, checking if the new expressions are skew-primitive and if the powers of the root vectors are central. To be able to see these tests better in the output, a short message is printed before them. This causes felix to evaluate these \verb1print1 commands and produces \verb1FALSE1 in the output. This is not a problem. Successful tests will give \verb1@ := 01 as output afterwards.
\end{itemize}

\section[Listing for self-linking of $A_2$]{Listing for self-linking of $\bm{A_2}$}\label{listA2}
%\hrulefill
\begin{verbatim}
link("tensor.mdl")$					

select int(lam12,lam21,mu1,mu2)<z,x1,x2,g1,g2,ixi,q;	
{{0,0,0,0,0,1,0},					
 {1,1,1,0,0,0,0},					
 {1,0,0,0,0,0,0},					
 {0,1,0,0,0,0,0},
 {0,0,0,1,1,0,0},
 {0,0,0,1,0,0,0},}
>$

si:=ideal(						
q*x1-x1*q, q*x2-x2*q, q*z-z*q, 				
q*g1-g1*q, q*g2-g2*q, q*ixi-ixi*q, 
q^2+q+1, g1*g2-g2*g1,					

g1*x1-q*x1*g1, g2*x1-q*x1*g2,				
g1*x2-q*x2*g1, g2*x2-q*x2*g2,
g1*z-q^2*z*g1, g2*z-q^2*z*g2, 

x1*x2-q*x2*x1-z,					
x1*z-q^2*z*x1-lam12*(1-g1^2*g2),			
x2*(x2*x1-q*x1*x2)-q^2*(x2*x1-q*x1*x2)*x2-lam21*(1-g2^2*g1),	

x1^3-mu1*(1-g1^3),					
x2^3-mu2*(1-g2^3)
)$
 
si:=standard(si)$					
%bye$%

delx1:=g1*ixi*x1+x1*ixi$				
delx2:=g2*ixi*x2+x2*ixi$				
delz:=remainder(ttimes(delx1,delx2)
           -ttimes(q*ixi,ttimes(delx2,delx1)),si)$	
delz2:=remainder(ttimes(delz,delz),si)$			
delz3:=remainder(ttimes(delz2,delz),si)$		

f0:=(1-q)^3*mu1*mu2$					
f1:=(1-q)*q*lam12*lam21$

v:=z^3
+f0*(1-g2^3)						
+f1*(1-g2^2*g1)$

delv:=delz3						
+f0*(ixi-g2^3*ixi*g2^3)
+f1*(ixi-g2^2*g1*ixi*g2^2*g1)$

print("TEST IF v SKEW-PRIMITIVE")$
remainder(delv-v*ixi-g1^3*g2^3*ixi*v,si)$		
print("TEST IF z^3 CENTRAL")$
remainder(z^3*x1-x1*z^3,si)$				
remainder(z^3*x2-x2*z^3,si)$				
bye$
\end{verbatim}
%\hrulefill

\section[Listing for self-linking of $B_2$]{Listing for self-linking of $\bm{B_2}$}\label{listB2}
%\hrulefill
\begin{verbatim}
link("tensor.mdl")$

select int(lam12,lam21,mu1,mu2)<u,z,x1,x2,g1,g2,ixi,q;
{{0,0,0,0,0,0,1,0},
 {1,1,1,1,0,0,0,0},
 {1,0,0,0,0,0,0,0},
 {0,1,0,0,0,0,0,0},
 {0,0,1,0,0,0,0,0},
 {0,0,0,0,1,1,0,0},
 {0,0,0,0,1,0,0,0}}
>$

si:=ideal(
 
q*u-u*q, q*z-z*q, q*x1-x1*q, q*x2-x2*q, 
q*g1-g1*q, q*g2-g2*q, q*ixi-ixi*q,
q^4+q^3+q^2+q+1, g1*g2-g2*g1,

g1*u-q^4*u*g1, g2*u-q^4*u*g2,
g1*z-q^3*z*g1, g2*z-q^3*z*g2,
g1*x2-q*x2*g1, g2*x2-q*x2*g2,
g1*x1-q^2*x1*g1, g2*x1-q^2*x1*g2,
 
x2*x1-q^2*x1*x2-z,
x2*z-q^3*z*x2-u,

x1*(x1*x2-q*x2*x1)-q^3*(x1*x2-q*x2*x1)*x1-lam12*(1-g1^2*g2),
x2*u-q^4*u*x2-lam21*(1-g2^3*g1),

x1^5-mu1*(1-g1^5),
x2^5-mu2*(1-g2^5)
)$
 
si:=standard(si)$
%bye$%

delx1:=g1*ixi*x1+x1*ixi*1$
delx2:=g2*ixi*x2+x2*ixi*1$
delz:=remainder(ttimes(delx2,delx1)
            -ttimes(q^2*ixi,ttimes(delx1,delx2)),si)$
delz2:=remainder(ttimes(delz,delz),si)$
delz3:=remainder(ttimes(delz2,delz),si)$
delz4:=remainder(ttimes(delz3,delz),si)$
delz5:=remainder(ttimes(delz4,delz),si)$

f0:=-mu1*mu2*(q^2-1)^5$
f1:=(-q^3+q^2+q-1)*lam12^2*lam21$
f2:=(2*q^3+2*q^2+1)*lam12^2*lam21$
f3:=(q^3+2*q^2+3*q-1)*lam12*lam21$
f4:=-q*(q^2+3*q+1)*lam12^2*lam21$

v:=z^5
+f0*(1-g1^5)
+f1*g1^5*g2^5
+f2*g1^4*g2^2
+f3*z*x1*g1^2*g2
+f4*g1^2*g2$

delv:=remainder(delz5
+f0*(ixi-g1^5*ixi*g1^5)
+f1*g1^5*g2^5*ixi*g1^5*g2^5
+f2*g1^4*g2^2*ixi*g1^4*g2^2
+f3*ttimes(ttimes(delz,delx1),g1^2*g2*ixi*g1^2*g2)
+f4*g1^2*g2*ixi*g1^2*g2
,si)$

print("TEST IF v SKEWPRIMITIVE")$
remainder(delv-v*ixi-g1^5*g2^5*ixi*v,si)$
print("TEST IF z^5 CENTRAL")$
remainder(z^5*x1-x1*z^5,si)$
remainder(z^5*x2-x2*z^5,si)$
%bye$%

delu:=remainder(ttimes(delx2,delz)
            -ttimes(q^3*ixi,ttimes(delz,delx2)),si)$
delu2:=remainder(ttimes(delu,delu),si)$
delu3:=remainder(ttimes(delu2,delu),si)$
delu4:=remainder(ttimes(delu3,delu),si)$
delu5:=remainder(ttimes(delu4,delu),si)$

k0:=-2*(q-1)^5*mu2$
k1:=-(q^2-1)^5*(q-1)^5*mu2^2$
k2:=-5*(2*q^3-q^2+q-2)*lam12^2*lam21$
k3:=-(3*q^3+q^2-q+2)*lam12*lam21^3$
k4:=(3*q^3+q^2+4*q+2)*lam12*lam21^3$
k5:=-5*(q^3+1)*lam12*lam21^2$
k6:=5*(q^3-1)*lam12*lam21$
k7:=-2*(q^3+2*q^2+3*q-1)*lam12*lam21^3$
k8:=-5*(q^3+q^2-1)*lam12*lam21^2$
k9:=(2*q^3+4*q^2+q-2)*lam12*lam21^3$

w:=u^5
+k0*v
+k1*x1^5
+k2*x2^5*g1^5*g2^5
+k3*g1^5*g2^10
+k4*g1^4*g2^7
+k5*u*x2*g1^3*g2^4
+k6*u^2*x2^2*g1^2*g2
+k7*g1^3*g2^4
+k8*u*x2*g1^2*g2
+k9*g1^2*g2$

delw:=remainder(delu5
+k0*delv
+k1*(x1^5*ixi+g1^5*ixi*x1^5)
+k2*ttimes(x2^5*ixi+g2^5*ixi*x2^5,g1^5*g2^5*ixi*g1^5*g2^5)
+k3*g1^5*g2^10*ixi*g1^5*g2^10
+k4*g1^4*g2^7*ixi*g1^4*g2^7
+k5*ttimes(ttimes(delu,delx2),g1^3*g2^4*ixi*g1^3*g2^4)
+k6*ttimes(ttimes(ttimes(ttimes(delu,delu),delx2),delx2)
           ,g1^2*g2*ixi*g1^2*g2)
+k7*g1^3*g2^4*ixi*g1^3*g2^4
+k8*ttimes(ttimes(delu,delx2),g1^2*g2*ixi*g1^2*g2)
+k9*g1^2*g2*ixi*g1^2*g2
,si)$

print("TEST IF w SKEW-PRIMITIVE")$
remainder(delw-w*ixi-g1^5*g2^10*ixi*w,si)$
print("TEST IF u^5 CENTRAL")$
remainder(u^5*x1-x1*u^5,si)$
remainder(u^5*x2-x2*u^5,si)$
bye$
\end{verbatim}
%\hrulefill

\section[Listing for self-linking of $G_2$]{Listing for self-linking of $\bm{G_2}$}\label{listG2}
%\hrulefill
\begin{verbatim}
link("tensor.mdl")$

select int(lam12,lam21,mu1,mu2)<w,v,u,z,x1,x2,g1,g2,ixi,q;
{{0,0,0,0,0,0,0,0,1,0},
 {5,4,3,2,1,1,0,0,0,0},
 {2,2,1,1,1,1,0,0,0,0},
 {0,0,0,0,1,1,0,0,0,0},
 {0,0,0,0,1,0,0,0,0,0},
 {0,0,0,1,0,0,0,0,0,0},
 {0,0,1,0,0,0,0,0,0,0},
 {0,1,0,0,0,0,0,0,0,0},
 {0,0,0,0,0,0,1,1,0,0},
 {0,0,0,0,0,0,1,0,0,0}}
>$

si:=ideal(
q*w-w*q, q*v-v*q, q*u-u*q, q*z-z*q,
q*x2-x2*q, q*x1-x1*q, q*g2-g2*q, q*g1-g1*q,
q^6+q^5+q^4+q^3+q^2+q+1, g1*g2-g2*g1, q*ixi-ixi*q,

g1*w-q^2*w*g1, g2*w-q^2*w*g2,
g1*v-q^6*v*g1, g2*v-q^6*v*g2,
g1*u-q^5*u*g1, g2*u-q^5*u*g2,
g1*z-q^4*z*g1, g2*z-q^4*z*g2,
g1*x2-q^3*x2*g1, g2*x2-q^3*x2*g2,
g1*x1-q*x1*g1, g2*x1-q*x1*g2,

z-x2*x1+q*x1*x2,
u-z*x1+q^2*x1*z,
v-u*x1+q^3*x1*u,
w-z*u+q^3*u*z,	
   
x1*(x1*(x1*(x1*x2-q^3*x2*x1)     -q^4*(x1*x2-q^3*x2*x1)*x1)
       -q^5*(x1*(x1*x2-q^3*x2*x1)-q^4*(x1*x2-q^3*x2*x1)*x1)*x1)
   -q^6*(x1*(x1*(x1*x2-q^3*x2*x1)-q^4*(x1*x2-q^3*x2*x1)*x1)
       -q^5*(x1*(x1*x2-q^3*x2*x1)-q^4*(x1*x2-q^3*x2*x1)*x1)*x1)*x1    
 -lam12*(1-g1^4*g2),
x2*z-q^4*z*x2-lam21*(1-g1*g2^2),

x1^7-mu1*(1-g1^7),
x2^7-mu2*(1-g2^7)
)$
 
si:=standard(si)$
%bye$%

delx1:=g1*ixi*x1+x1*ixi*1$
delx2:=g2*ixi*x2+x2*ixi*1$

delz:=remainder(ttimes(delx2,delx1)
            -ttimes(q*ixi,ttimes(delx1,delx2)),si)$
delz2 :=remainder(ttimes(delz,delz),si)$
delz3 :=remainder(ttimes(delz2,delz),si)$
delz4 :=remainder(ttimes(delz3,delz),si)$
delz5 :=remainder(ttimes(delz4,delz),si)$
delz6 :=remainder(ttimes(delz5,delz),si)$
delz7 :=remainder(ttimes(delz6,delz),si)$

f0:=(1-q^3)^7*mu1*mu2$
f1:=(-2*q^5-4*q^4+q^3-q^2+4*q+2)*lam12*lam21$
f2:=(6*q^5+8*q^4+6*q^3-3*q-3)*lam12*lam21^2$
f3:=(-q^4-3*q^3+q^2-3*q-1)*lam12*lam21^2$
f4:=(2*q^5+2*q^4+4*q^3+5*q^2+2*q-1)*lam12*lam21^3$
f5:=(q^5-2*q^4-4*q^3-7*q^2-6*q-3)*lam12*lam21^3$
f6:=(-4*q^5-2*q^4-2*q^3+2*q^2+2*q+4)*lam12*lam21^3$
f7:=(q^5+2*q^4+2*q^3+2*q)*lam12*lam21^3$

zz:=z^7
+f0*(1-g1^7)
+f1*z^2*x2^2*g1^4*g2
+f2*z*x2*g1^5*g2^3
+f3*z*x2*g1^4*g2
+f4*g1^7*g2^7
+f5*g1^6*g2^5
+f6*g1^5*g2^3
+f7*g1^4*g2$

delzz:=delz7
+f0*(ixi-g1^7*ixi*g1^7)
+f1*ttimes(ttimes(delz2,ttimes(delx2,delx2)),g1^4*g2*ixi*g1^4*g2)
+f2*ttimes(ttimes(delz,delx2),g1^5*g2^3*ixi*g1^5*g2^3)
+f3*ttimes(ttimes(delz,delx2),g1^4*g2*ixi*g1^4*g2)
+f4*g1^7*g2^7*ixi*g1^7*g2^7
+f5*g1^6*g2^5*ixi*g1^6*g2^5
+f6*g1^5*g2^3*ixi*g1^5*g2^3
+f7*g1^4*g2*ixi*g1^4*g2$

print("TEST IF zz SKEW_PRIMITIVE")$
remainder(delzz-zz*ixi-g1^7*g2^7*ixi*zz,si)$
print("TEST IF z^7 CENTRAL")$
remainder(z^7*x1-x1*z^7,si)$
remainder(z^7*x2-x2*z^7,si)$
bye$
\end{verbatim}
%\hrulefill

%%% Local Variables: 
%%% mode: latex
%%% TeX-master: "main"
%%% End: 

%\newcommand{\liti}[6]{#1, \emph{#2}, #3 {\bf #4} (#5), #6.}%1Name,2Titel,3Journal,4Number,5Year,6Pages

%%% Local Variables: 
%%% mode: latex
%%% TeX-master: "main"
%%% End: 

{\noindent \Large\bf Summary}\\
\addcontentsline{toc}{chapter}{Summary/Zusammenfassung} 

In this thesis we want to contribute to some classification results for pointed Hopf algebras with abelian coradical found recently by Andruskiewitsch and Schneider \cite{AS-p3,AS-p17,AS-Poin,AS-char}. Their lifting method produces new classes of Hopf algebras. These algebras are constructed from a linking datum consisting of a group, a Dynkin diagram, some linking parameters and a number of group elements and characters fulfilling certain compatibility conditions. These conditions are rather implicit and hence an explicit description of these Hopf algebras is often not easy. In this work we treat various aspects of such a description in detail.

 One of our main contributions is the clarification of the  concept of linking. Based on the original work \cite{AS-p17}, we first introduce some suitable terminology, Definitions \ref{link.diag}-\ref{link.alg}. Then we give an easily applicable criterion, Theorem \ref{mainthm}, that helps in deciding which linkings can produce finite dimensional Hopf algebras and what possible restrictions have to be imposed on the coradical. This involves simply counting certain objects in graphs and computing the so-called genus from this data. We extend this result to treat affine Dynkin diagrams as well, Theorem \ref{affinethm}. Examples of ``exotic'' linkings are given in Figure \ref{exotics}.  Some exceptional cases that usually have to be excluded from classification results come from setups we call self-linkings. We present the prototypes of Hopf algebras arising from such situations in Section \ref{rank2comps}. The new Hopf algebras derived from the diagram $B_2,$ which we compute using a Computer algebra program, are given in Figure \ref{B2rel}.

Another open question concerns the compatibility of the groups and the Dynkin diagrams in a linking datum. Although a general answer seems out of reach, we are able to contribute an answer for the groups $(\Z/(p))^2$ in Theorem \ref{groupthm}. We prove that apart from a few exceptions, all diagrams with at most four vertices can be used for the construction of finite dimensional pointed Hopf algebras with these groups as the coradical. %????????Furthermore we derive some necessary conditions that can be used to exclude some linkings 

Finally, the last major topic of this thesis is the investigation of the relation between the new Hopf algebras constructed by the lifting method. It turns out that different linking parameters lead to quasi-isomorphic Hopf algebras, Theorem \ref{linkthm}. All Hopf algebras that arise from the lifting method using only Dynkin diagrams of type $A_n$ display the same behaviour, Theorem \ref{mixedthm}. This means that all the finite dimensional pointed Hopf algebras constructed in this way, which only differ in their choice of parameters are 2-cocycle deformations of each other. Our proof should be easily adaptable to the Hopf algebras associated with the other types of finite Dynkin diagrams, once all parameters have been determined for these algebras explicitly. This raises the hope that Masuoka's conjecture in \cite{Mas} can be saved in spite of the counter-example in \cite{EG} by specializing it slightly (page \pageref{conjec}).\\

{\noindent \Large\bf Zusammenfassung}\\

In dieser Dissertation wollen wir zu Klassifizierungsresultaten f\"ur punktierte Hopfalgebren mit abelschem Koradikal beitragen, die vor kurzer Zeit von Andruskiewitsch und Schneider \cite{AS-p3,AS-p17,AS-Poin,AS-char} gewonnen wurden. Deren Liftingmethode erzeugt neue Klassen von Hopfalgebren. Diese Algebren werden ausgehend von einem Verbindungs-Datum konstruiert, welches aus einer Gruppe, einem Dynkin-Diagramm, einigen Parametern und einer Reihe von Gruppenelementen und -charakteren, die gewisse Kompatibilit\"atsbedingungen erf\"ullen, besteht. Diese Bedingungen sind ziemlich implizit gegeben und erlauben daher meist keine einfache explizite Beschreibung dieser Hopfalgebren. In dieser Arbeit behandeln wir ausf\"uhrlich verschiedene Aspekte solch einer Beschreibung.

Einer unserer wesentlichen Beitr\"age ist die explizite Ausarbeitung des Verbindungskonzepts. Aufbauend auf der Originalarbeit \cite{AS-p17} f\"uhren wir zun\"achst eine passende Terminologie ein, Definitionen \ref{link.diag}-\ref{link.alg}. Danach geben wir ein einfach anwendbares Kriterium an, Theorem \ref{mainthm}, das entscheiden hilft, welche Verbindungen zu endlichdimensionalen Hopfalgebren f\"uhren und was f\"ur Bedingungen an das Koradikal gestellt werden m\"ussen. Dazu mu\ss{}  man nur, gewisse Objekte in Graphen z\"ahlen und aus diesen Daten das sogenannte Geschlecht berechnen. Wir erweitern unser Resultat in Theorem \ref{affinethm} auch auf affine Dynkin-Diagramme. Beispiele exotischer Verbindungen gibt Abbildung \ref{exotics}. Einige Ausnahmef\"alle, welche normalerweise von den Klassifizierungsresultaten ausgenommen werden m\"ussen, entstehen in Situationen, die wir Selbstverbindungen nennen. Wir pr\"asentieren die Prototypen von Hopfalgebren in solch einer Situation in Abschnitt \ref{rank2comps}. Die neuen Hopfalgebren, die vom Diagramm $B_2$ stammen und welche wir mittels eines Computeralgebraprogramms berechnen, sind in Abbildung \ref{B2rel} angegeben.

Eine weitere offene Frage betrifft die Kompatibilit\"at der Gruppen mit den Dynkin-Diagrammen in einem Verbindungs-Datum. Obwohl eine generelle Antwort au\ss er Reichweite scheint, k\"onnen wir mit Theorem \ref{groupthm} eine Antwort f\"ur die Gruppen $(\Z/(p))^2$ geben. Wir beweisen, da\ss , bis auf wenige Ausnahmen, alle Diagramme mit maximal vier Ecken f\"ur die Konstruktion endlichdimensionaler punktierter Hopfalgebren mit diesen Gruppen als Koradikal benutzt werden k\"onnen.

Das letzte gro\ss e Thema dieser Arbeit ist die Untersuchung der Beziehungen zwischen den, mittels der Liftingmethode konstruierten, neuen Hopfalgebren. Es stellt sich heraus, da\ss{}  verschiedene Verbindungsparameter zu quasi-isomorphen Hopfalgebren f\"uhren, Theorem \ref{linkthm}. Alle Hopfalgebren, die mittels der Liftingmethode nur aus Diagrammen des Typs $A_n$ entstehen, zeigen das gleiche Verhalten, Theorem \ref{mixedthm}. Das hei\ss t, da\ss{}  alle so konstruierten, endlichdimensionalen punktierten Hopfalgebren, die sich nur durch die Wahl ihrer Parameter unterscheiden, 2-Kozyklus-Deformationen voneinander sind. Unser Beweis sollte sich einfach auf Hopfalgebren, die mit den anderen Typen endlicher Dynkin-Diagramme assoziiert sind, \"ubertragen lassen, sobald alle Parameter f\"ur diese Algebren explizit bestimmt worden sind. Dies best\"arkt die Hoffnung, da\ss{}  Masuokas Vermutung in \cite{Mas} trotz des Gegenbeispiels in \cite{EG} gerettet werden kann, indem man sie etwas spezialisiert (Seite \pageref{conjec}).

%%% Local Variables: 
%%% mode: latex
%%% TeX-master: "main"
%%% End: 

%\documentclass[a4paper,12pt]{letter}
%\usepackage{ngerman}
\addtolength{\oddsidemargin}{1cm}
\addtolength{\textwidth}{2pt}
\newcommand{\capti}[1]{{\begin{center}\LARGE\bf{#1}\end{center}}}
\newcommand{\fett}[1]{{\large\bf{#1}}}
\thispagestyle{empty}
\sloppy

%\begin{document}

\capti{Lebenslauf}
\vspace{2em}
\fett{Pers\"onliche Daten}
\begin{quote}
Daniel Didt\\

gebren am 3.~August 1974 in Rostow am Don (Ru\ss land)\\
verheiratet; 4 Kinder\\
\end{quote}
\fett{Schulausbildung}
\begin{quote}
\begin{description}
\item[1981-1989 ] Besuch einer Polytechnischen Oberschule in Leipzig und
ab 1983 einer Klasse mit erweitertem Russischunterricht
\item[1989-1993 ] Nach erfolgreicher Aufnahmepr\"ufung Besuch der
Spezialschule math.-naturwiss.-techn. Richtung  ``Wilhelm Ostwald'' in
Leipzig
\item[Juli 1993 ] Abitur
\end{description}
\end{quote}
\fett{Hochschulausbildung}
\begin{quote}
\begin{description}
\item[Okt.~1993 ] Beginn eines Doppelstudiums (Diplom-Mathematik und
Diplom-Physik) an der Universit\"at Leipzig
\item[Sept.~1995 ] Vordiplom in Mathematik (Nebenfach Informatik) und Physik
(Nebenfach Chemie)
\item[Feb.~1996 - Jan.~1999 ] Mitglied der Studienstiftung des Deutschen Volkes
\item[Sept. 1996 - Juli 1998 ] Leiter eines Mathematikzirkels f\"ur die Klassen 11/12 der ``Leipziger Sch\"ulergesellschaft f\"ur Mathematik''
\item[Jan.~1999 ] Abschluss des Physikstudiums mit Diplom%Thema der Diplomarbeit: "`Differentialkalk\"ule auf inhomogenen Quantengruppen"'
\item[M\"arz 2000 ] Abschluss des Mathematikstudiums mit Diplom
\item[Jan.~2000 - Dez.~2002] Promotion am Mathematischen Institut der LMU M\"unchen im Graduiertenkolleg ``Mathematik im Bereich ihrer Wechselwirkung mit der Physik''
\end{description}
\end{quote}

%\end{document}
\end{document}